\ifx \labelsONmargin\undefined

  \ifx\RequirePackage\undefined
    \expandafter\ifx\csname amsart.sty\endcsname\relax
        \documentstyle[12pt,amssymb,amscd]{amsart}
    \fi

    \def\atSign{@@}

    \def\mathbb{\Bbb}
    \def\mathfrak{\frak}
    \def\mathbf{\bold}
    \ifx\boldsymbol\undefined	
      \def\boldsymbol#1{{\bold #1}}
    \fi

    \def\mathbit{\boldsymbol}

    \makeatletter
    \newenvironment{proof}{%
         \@ifnextchar[{%
                       \expandafter\let\expandafter\end@proof
                         \csname endpf*\endcsname
                         \my@proof
                      }{\let\end@proof\endpf\pf}%
        }{\end@proof}
    \def\my@proof[#1]{\@nameuse{pf*}{#1}}
    \makeatother

    \def\xrightarrow[#1]#2{@>{#2}>{#1}>}
    \def\xleftarrow[#1]#2{@<{#2}<{#1}<}

    \def\providecommand#1{\def#1}
    \def\emph#1{{\em #1}}
    \def\textbf#1{{\bf #1}}

    \def\mathring{\overset{\,\,{}_\circ}}
  \else
      \ifx\usepackage\RequirePackage
      \else
        \documentclass[12pt]{amsart}
        \usepackage{amssymb,amscd}
      \fi
      \ifx \mathring\undefined		
        \DeclareMathAccent{\mathring}{\mathalpha}{operators}{"17}
      \fi

      \long\def\FAKEendPROOF{\endtrivlist}
      \ifx\endproof\FAKEendPROOF
	  \def\endproof{\qed\endtrivlist}
      \fi

      \ifx\mathbit\undefined	
        \DeclareMathAlphabet{\mathbit}{OML}{cmm}{b}{it}
      \fi

      \def\atSign{@}

      \advance\oddsidemargin -0.1\textwidth
      \evensidemargin\oddsidemargin
      \def\Sb#1\endSb{_{\substack{#1}}}
      \def\Sp#1\endSp{^{\substack{#1}}}

  \fi
\makeatletter
\@ifundefined{DeclareFontShape} 
     {\@ifundefined{selectfont}
             {
                \message{We need AmS-LaTeX, this version of LaTeX does not 
                        support it}\@@end
             }{
                \def\mathcal{\cal}
                
                \def\pcyr{%
                        \def\default@family{UWCyr}%
                        \let\oldSl@\sl
                        \def\sl{\def\default@shape{it}\oldSl@}%
                        \cyracc
                        \language\Russian\family{UWCyr}\selectfont
                }
             }}{
                \DeclareFontEncoding{OT2}{}{\noaccents@}
                \DeclareFontSubstitution{OT2}{cmr}{m}{n}
                \DeclareFontFamily{OT2}{cmr}{\hyphenchar\font45 }
                \DeclareFontShape{OT2}{cmr}{m}{n}{%
                     <5><6><7><8><9><10>gen*wncyr %
                     <10.95><12><14.4><17.28><20.74><24.88> wncyr10 %
                }{}
                \DeclareFontShape{OT2}{cmr}{m}{it}{%
                     <5><6><7><8><9><10> gen * wncyi%
                     <10.95><12><14.4><17.28><20.74><24.88> wncyi10%
                }{}
                \DeclareFontShape{OT2}{cmr}{bx}{n}{%
                     <5><6><7><8><9><10> gen * wncyb%
                     <10.95><12><14.4><17.28><20.74><24.88> wncyb10%
                }{}
                \DeclareFontShape{OT2}{cmr}{m}{sl}{%
                     <-> ssub * cmr/m/it%
                }{}
                \DeclareFontShape{OT2}{cmr}{m}{sc}{%
                     <5><6><7><8><9><10>%
                     <10.95><12><14.4><17.28><20.74><24.88> wncysc10%
                }{}
                \DeclareFontFamily{OT2}{cmss}{\hyphenchar\font45 }
                \DeclareFontShape{OT2}{cmss}{m}{n}{%
                     <8><9><10> gen * wncyss%
                     <10.95><12><14.4><17.28><20.74><24.88> wncyss10%
                }{}
                
                \def\cyrencodingdefault{OT2}
                \def\pcyr{%
                        \cyracc
                        \let\encodingdefault\cyrencodingdefault
                        \language\Russian\fontencoding{OT2}\selectfont
                }
             }   

\@ifundefined{theorembodyfont} 
     {
        \def\theorembodyfont#1{\relax}
        \makeatletter
          \let\@@th@plain\th@plain
          \def\th@plain{ \@@th@plain \slshape }
        \makeatother
        \let\normalshape\relax
     }{}   

\ifx\cprime\undefined
     \def\cprime{$'$}
\fi
\makeatletter
  \def\@sect@my#1#2#3#4#5#6[#7]#8{%
\ifnum #2>\c@secnumdepth
   \let\@svsec\@empty
 \else
   \refstepcounter{#1}%
\edef\@svsec{\ifnum#2<\@m
             \@ifundefined{#1name}{}{\csname #1name\endcsname\ }\fi
\noexpand\rom{\csname the#1\endcsname.}\enspace}\fi
 \@tempskipa #5\relax
 \ifdim \@tempskipa>\z@ 
   \begingroup #6\relax
   \@hangfrom{\hskip #3\relax\@svsec}{\interlinepenalty\@M #8\par}%
   \endgroup
   \if@article\else\csname #1mark\endcsname{%
        \ifnum \c@secnumdepth >#2\relax\csname the#1\endcsname. \fi#7}\fi
\ifnum#2>\@m \else
       \let\@tempf\\ \def\\{\protect\\}\addcontentsline{toc}{#1}%
{\ifnum #2>\c@secnumdepth \else
             \protect\numberline{%
               \ifnum#2<\@m
               \@ifundefined{#1name}{}{\csname #1name\endcsname\ }\fi
               \csname the#1\endcsname.}\fi
           #8}\let\\\@tempf
     \fi
 \else
  \def\@svsechd{#6\hskip #3\@svsec
    \@ifnotempty{#8}{\ignorespaces#8\unskip
       \ifnum\spacefactor<1001.\fi}%
        \ifnum#2>\@m \else
          \let\@tempf\\ \def\\{\protect\\}\addcontentsline{toc}{#1}%
            {\ifnum #2>\c@secnumdepth \else
              \protect\numberline{%
                \ifnum#2<\@m
                \@ifundefined{#1name}{}{\csname #1name\endcsname\ }\fi
                \csname the#1\endcsname.}\fi
             #8}\let\\\@tempf\fi}%
 \fi
\@xsect{#5}}
  \ifx\RequirePackage\undefined
  \let\@sect\@sect@my             
  \fi
  \ifx\RequirePackage\undefined	
  \def\th@remark@my{\theorempreskipamount6\p@\@plus6\p@
    \theorempostskipamount\theorempreskipamount
    \def\theorem@headerfont{\it}\normalshape}
    \let\th@remark\th@remark@my
  \else				
    \let\o@@remark\th@remark

    \ifx\theorempostskipamount\undefined		
    \else						
      \def\th@remark{\o@@remark
	\ifdim\theorempostskipamount < 2pt\relax
	  \theorempostskipamount\theorempreskipamount
	     \multiply\theorempostskipamount\tw@
	     \divide\theorempostskipamount\thr@@
	\fi
      }
    \fi
  \fi
\let\myLabel\@gobble
\def\labelsONmargin{\@mparswitchfalse\def\myLabel##1{\@bsphack\marginpar
                                  {\normalshape\tiny\rm Label ##1}\@esphack}}
\ifx \url\undefined
  \def\url#1{{\tt #1}}%
\fi
\def\cyracc{\def\u##1{
                \if \i##1\char"1A%
                \else \if I##1\char"12%
                \else \accent"24 ##1\fi\fi }%
\def\"##1{\if e##1{\char"1B}%
                \else \if E##1{\char"13}%
                \else \accent"7F ##1\fi\fi }%
\def\9##1{\if##1z\char"19 
\else\if##1Z\char"11 
\else\if##1E\char"03 
\else\if##1e\char"0B 
\else\if##1u\char"18 
\else\if##1U\char"10 
\else\if##1A\char"17 
\else\if##1a\char"1F 
\else\if##1p\char"7E 
\else\if##1P\char"5E 
\else\if##1Q\char"5F 
\else\if##1q\char"7F 
\else\if##1i\char"1A 
\else\if##1I\char"12 
\else\if##1N\char"7D 
\fi
\fi
\fi
\fi
\fi
\fi
\fi
\fi
\fi
\fi
\fi
\fi
\fi
\fi
\fi
}%
\def\cydot{{\kern0pt}}}%
\def\cydot{$\cdot$}

\ifx\Russian\undefined
        \def\Russian{0\relax
    \message{Don't know the hyphenation rules for Russian^^J
                        Please do INITeX with `input  russhyph' in the 
                        command line}%
                \gdef\Russian{0\relax}%
        }
\fi
\ifx\RequirePackage\undefined			
  \def\@putname#1#2#3#4{\def\@@ref{#3}\let\old@bf\bf
        \def\bf##1{\old@bf\if?\noexpand##1?{#4}\else##1\fi}%
	#1{#2}%
        \let\bf\old@bf}
\else						
  \def\@putname#1#2#3#4{\def\@@ref{#3}\let\old@bf\bf	
	\let\old@reset@font\reset@font			
        \def\bf##1{\old@bf\if?\noexpand##1?{#4}\else##1\fi}%
	\def\reset@font##1##2{\old@reset@font##1\if?\noexpand##2?{#4}\else##2\fi}#1{#2}%
        \let\bf\old@bf\let\reset@font\old@reset@font}
\fi
\let\my@ref=\ref
\def\ref#1{\@putname\my@ref{#1}{#1}{\tiny\rm\@@ref}}
\let\my@pageref=\pageref
\def\pageref#1{\@putname\my@pageref{#1}{#1}{\tiny\rm\@@ref}}
\let\my@cite=\cite
\def\cite#1{\@putname\my@cite{#1}{\@citeb}{\tiny\rm\@@ref}}
\makeatother
\fi
\relax 
\theoremstyle{plain} 
\theorembodyfont{\sl}

\numberwithin{equation}{section}
\theorembodyfont{\rm}
\theoremstyle{definition}

\newtheorem{definition}{Definition}[section]

\newtheorem{example}[definition]{Example}

\theoremstyle{remark}

\newtheorem{remark}[definition]{Remark} 

\theoremstyle{plain} 
\theorembodyfont{\sl}

\newtheorem{theorem}[definition]{Theorem}
\newtheorem{lemma}[definition]{Lemma}
\newtheorem{corollary}[definition]{Corollary}
\newtheorem{proposition}[definition]{Proposition}
\newtheorem{amplification}[definition]{Amplification}

\begin{document}
\bibliographystyle{amsplain}
\relax 

\title[Small degree Riemann--Roch theorem]{ Curves of infinite genus I \\
Riemann--Roch theorem for small degree }

\author{ Ilya Zakharevich }

\address{ Department of Mathematics, Ohio State University, 231 W.~18~Ave,
Columbus, OH, 43210 }

\email {ilya\atSign{}math.ohio-state.edu}

\date{ February 2002\quad Archived as \url{math.AG/0202057}\quad Printed: \today
}

\setcounter{section}{-1}

\maketitle
\begin{abstract}
The most useful and interesting line bundles over algebraic
curves of a very high genus have the ratio $ \delta $ of the degree to the genus
close to half-integer values, usually $ \delta\approx0 $, $ \delta\approx1/2 $, or $ \delta\approx1 $; the numeric
properties are very different in these three cases. This leads to three
different theories for curves of infinite genus.

For analytic curves of infinite genus, to get a theory parallel to
algebraic geometry one needs to restrict attention to holomorphic
sections satisfying some ``conditions on growth at infinity''. Each such
condition effectively attaches an ``ideal point'' to the curve; this
process is similar to compactification.

The theory of holomorphic functions on curves with such ``ideal
points'' is developed (the variant presented in the first part of the
series is tuned to the case $ \delta\approx0 $). Conditions on the ``lengths of handles''
of the curve are found which ensure the geometry to be parallel to
algebraic geometry.

It turns out that these conditions give no restriction on the
density of ideal points on the curve. In particular, such curves may have
a dense set of ideal points; these curves have no smooth points at all,
and have a purely fractal nature. (Such ``foam'' curves live near the
``periphery'' of the corresponding $ g=\infty $ moduli space; one needs to study
these curves too, since they may be included in the support of natural
measures arising on the moduli spaces.)

\end{abstract}
\tableofcontents

\section{Introduction }

\subsection{Motivations }\label{s02.20}\myLabel{s02.20}\relax  The desire to find a set of curves of infinite
genus for which ``most'' theorems of algebraic geometry hold stems from the
following:
\begin{enumerate}
\item
The existence of algebro-geometric description of a dense set of
special solutions of infinite-dimensional integrable systems; thus a hope
to describe other solutions by extending the algebraic geometry;
\item
The ability to describe the series of perturbation theory for string
amplitudes as integrals over moduli spaces of algebraic curves; thus a
hope to describe non-perturbative terms as integrals over some ambient
space which includes all the moduli spaces of algebraic curves;
\item
The hope that an appropriate completion of the union of compactified
moduli spaces may have a simpler geometry than the moduli spaces
themselves (stability phenomenon).
\end{enumerate}

In \cite{Zakh97Qua} we introduced several heuristics to describe a
possible candidate for such a completion; outline the principal
ingredients:
\begin{enumerate}
\item
one works with a pair of a curve $ C $ and a line bundle $ {\mathcal L} $ on $ C $;
\item
a curve of infinite genus is a small ``deformation'' of an algebraic
curve; small ``deformations'' of a pair should change the space of global
sections of $ {\mathcal L} $ (and $ H^{1}\left(C,{\mathcal L}\right) $) in a minimal possible way;
\item
the description should be as conformally invariant as possible;
\item
given a long tube $ A $ and a $ 1 $-form $ G $ such that $ \operatorname{Supp} G $ is ``deep
inside'' $ A $ and $ \int G\not=0 $, the set $ \left\{F|_{\partial A} \mid \bar{\partial}F=cG, c\in{\mathbb C}\right\} $ depends very
weakly\footnote{E.g., consider Fourier coefficients of $ F $ on two boundary circles for
$ A=\left\{1/N<|z|<N\right\} $, $ N\gg 1 $.} on $ G $ ({\em robustness\/});
\item
enumerative algebraic geometry ``requires'' compactification of
algebraic curves; one can substitute compactification by growth
conditions: the local sections of the line bundle $ {\mathcal O}\left(n\cdot P\right) $, $ P=\infty $, on $ {\mathbb P}^{1} $
coincide with local sections of $ {\mathcal O} $ on $ {\mathbb A}^{1} $ with the growth condition
$ O\left(|z|^{\beta}\right) $, $ n\leq\beta<n+1 $;
\item
Riemann--Roch theorem is a litmus test to check whether a procedure
to translate the growth conditions (such as $ \beta $) to degree $ n $ is ``correct'';
\item
duality theorem is a litmus test to check whether a growth condition
is ``reasonable'' (note the specialty of integer $ \beta $ w.r.t.~the change $ \beta \to
-\beta $ in the translation law $ n\leq\beta<n+1 $).
\end{enumerate}

In our theory the compactification process adds a very large
(possibly uncountable) set; thus to obtain a working theory it is crucial
to start with ``reasonable'' growth conditions. In this paper we introduce
a class of curves and growth conditions which satisfies the heuristics
above. We expect that (when defined) the moduli space of such curves may
lead to answers to the questions at the beginning of this sections.

\subsection{Finite-genus cases }\label{s02.30}\myLabel{s02.30}\relax  Recall which ``small deformations'' are
involved in the ``usual completion'' of the moduli spaces. If both the
initial and resulting curves are smooth, it is a two-steps procedure.
First, gluing: given a curve $ C $ with a pair of points $ x_{1,2} $, identify these
points; the resulting quotient curve $ C' $ has a double-point singularity $ x $
(locally isomorphic to solutions to $ uv=0 $). Given a line bundle $ {\mathcal L} $ over $ C $
and an identification of $ {\mathcal L}|_{x_{1}} $ and $ {\mathcal L}|_{x_{2}} $, one obtains a line bundle $ {\mathcal L}' $ over
$ C' $; sections of $ {\mathcal L}' $ can be identified with sections $ F $ of $ {\mathcal L} $ satisfying
$ F\left(x_{1}\right)=F\left(x_{2}\right) $. By robustness, far from $ x_{1,2} $ the condition $ F\left(x_{1}\right)=F\left(x_{2}\right) $ can
be ``compensated'' by allowing a pole at a point $ P\in C' $ near $ x $; thus
$ \left(C',{\mathcal L}'\left(P\right)\right) $ is a small deformation of $ \left(C,{\mathcal L}\right) $.

Second, ``smoothing'' of the singularity: cover $ C' $ by a neighborhood $ U $
of the double point $ x\in C' $ and ``the rest'' $ U' $; $ U $ may be identified with an
open subset $ U_{0} $ of $ V_{0} $, here $ V_{\varepsilon}=\left\{\left(u,v\right)\in{\mathbb C}^{2} \mid uv=\varepsilon\right\} $; glue a curve $ C^{\varepsilon} $ by
replacing $ U $ by a ``deformation'' $ U_{\varepsilon}\subset V_{\varepsilon} $ of $ U_{0} $. Similarly, given
trivialization of $ {\mathcal L} $ near $ x_{1,2} $, one can deform $ {\mathcal L}' $ into a line bundle $ {\mathcal L}^{\varepsilon} $ on
$ C^{\varepsilon} $; thus $ {\mathcal L}'\left(P\right) $ into $ {\mathcal L}^{{\mathcal E}}\left(Q\right) $, $ Q\in V_{\varepsilon} $. Since coordinates $ \left(u_{0},v_{0}\right) $ of $ Q $ determine
$ \varepsilon=u_{0}v_{0} $, this gives a two-parametric deformation of the pair $ \left(C,{\mathcal L}\right) $.

\subsection{The dust }\label{s02.40}\myLabel{s02.40}\relax  Since $ g\left(C^{\varepsilon}\right)=g\left(C\right)+1 $, one needs infinitely many
operations of Section~\ref{s02.30} to obtain a curve of infinite genus. Each
such operation is equivalent to cutting two small disks $ R_{i} $, $ R_{i}' $ out of $ C $
(with radii satisfying $ r_{i}r_{i}'=|\varepsilon| $), then identifying $ \gamma_{i}=\partial R_{i} $ with $ \gamma_{i}'=\partial R_{i}' $
(the angle of rotation is determined by $ \operatorname{Arg} \varepsilon $) on $ D\buildrel{\text{def}}\over{=}C\smallsetminus\left(\bigcup_{I}R_{i}\cup R_{i}'\right) $.
Denote the result of these identifications by $ C^{*} $.

The description above has a gaping hole. Indeed, if $ I $ is infinite,
$ \partial D $ is strictly larger than $ \bigcup_{I}\left(\gamma_{i}\cup\gamma_{i}'\right) $; call the remaining part of $ \partial D $ the
{\em dust\/} of $ C^{*} $. There are two different possibilities: glue $ C^{*} $ out of $ D $, or
out of the interior $ \mathring{D} $ of $ D $. In the first case $ C^{*} $ is {\em very\/} non-smooth at
points of the dust, and we need to {\em define\/} which functions are
``holomorphic'' at these points. In the other case $ C^{*} $ is not compact;
compactification of $ C^{*} $ requires addition of the dust. By heuristics of
Section~\ref{s02.20}, to define line bundles on $ C^{*} $ one needs a growth
condition near each point of the dust. Thus either way leads to growth
conditions.

\subsection{Growth conditions and representations of $ \protect \operatorname{SL}_{2}\left({\mathbb C}\right) $ }\label{s02.50}\myLabel{s02.50}\relax  As the example
of $ {\mathbb P}^{1}={\mathbb A}^{1}\cup\left\{\infty\right\} $ of Section~\ref{s02.20} shows, the definition of ``appropriate''
growth conditions and degree is delicate even if dust consists of
isolated points; the situation should be much harder for ``massive'' dust.
The principal tool of our approach comes from the following observation:
given a typical Banach norm on the space of functions, the functions of
finite norm have a restricted growth rate near {\em every point}.

Suppose that the norm has good properties w.r.t.~gluing a function
from its restriction on open subsets. Then it makes sense to ask whether
a function has a finite norm ``where it is analytic''; most of the commonly
used norms satisfy this condition. If so, then finiteness of the norm is
a growth condition at each point of the set of non-analyticity of the
function. This gives us some growth conditions at the dust.

A selection principle for ``reasonable'' norms now includes conformal
invariance (or at least some strong bounds on how the norm can change
under conformal transformations). It also helps if the norm is Hilbert.
Recall the action of $ \operatorname{SL}\left(2,{\mathbb C}\right) $ on $ {\mathbb P}^{1} $ allows (see \cite{GelNai50Uni,Vog87Uni})
a family of {\em unitary\/} representations of $ \operatorname{SL}\left(2,{\mathbb C}\right) $ on sections of the line
bundles\footnote{For a real oriented line bundle $ {\mathcal L} $, $ {\mathcal L}^{\alpha} $ is well-defined for $ \alpha\in{\mathbb C} $;
same for $ {\mathcal L}\otimes_{{\mathbb R}}{\mathbb C} $.} $ {\mathcal L}_{\alpha}=\left(\omega\otimes\bar{\omega}\right)^{\alpha/2} $, $ \alpha\in\left(0,1\right) $; this gives a family of
conformally-invariant Hilbert norms $ \|\|_{\alpha} $. We call a section $ F $ on
$ D\subset C $ {}$ \alpha $-{\em acceptable\/} if $ F=G|_{D} $, here $ \|G\|_{\alpha}<\infty $. Now we can glue $ C^{*} $ out of
$ \mathring{D} $, as in Section~\ref{s02.40}, taking $ \alpha $-acceptability (on $ \mathring{D} $) as a growth
condition on the dust; this defines a line bundle\footnote{I.e., a ``usual'' line bundle on the complement to dust, plus growth
conditions on the dust.} $ \widetilde{{\mathcal L}}_{\alpha} $ on $ C^{*} $.

Similarly, one can define $ \widetilde{{\mathcal L}}_{\alpha}\left(n\cdot P\right) $ for $ P $ in the dust by requiring
that $ \left(z-P\right)^{n}F $ is acceptable instead of $ F $; such modified growth conditions
define, in our conventions, a different line bundle. Alternatively, $ n $ is
the ``multiplicity of a divisor'' at the dust points, or the contribution
of infinity to the degree of a divisor.

The Hilbert norms on sections of $ {\mathcal L}_{\alpha} $ are equivalent to the Sobolev
norms on spaces $ H^{1-\alpha}\left({\mathbb P}^{1},{\mathcal L}_{\alpha}\right) $; thus one can work with analogues of these
norms on any compact curve $ C $. Being geometrically-defined, bundles $ {\mathcal L}_{\alpha} $
have an added convenience of {\em auto-gluing\/}: any identification of curves $ \gamma $
and $ \gamma' $ on $ {\mathbb P}^{1} $ leads to identification of $ {\mathcal L}_{\alpha}|_{\gamma} $ and $ {\mathcal L}_{\alpha}|_{\gamma'} $. This provides a
very convenient ``base point'' on the Jacobian. As we will see it later, a
modification of the gluing between $ {\mathcal L}_{\alpha}|_{\gamma_{i}} $ and $ {\mathcal L}_{\alpha}|_{\gamma_{i}'} $ for infinite number
of indices can lead to a failure of the Riemann--Roch theorem; thus the
validity of the Riemann--Roch for this particular point of the Jacobian
is a property of the curve itself.

Now one can {\em define\/} a line bundle on $ C^{*} $ as a modification of the line
bundle $ \widetilde{{\mathcal L}}_{\alpha} $ by divisors on $ \mathring{D} $, or by modifications of the gluings of $ {\mathcal L}_{\alpha}|_{\gamma_{i}} $
and $ {\mathcal L}_{\alpha}|_{\gamma_{i}'} $, or by modifications of growth conditions at several points of
the dust (thus ``a divisor on the dust'').

\begin{remark} In fact, in the body of this paper we do not use the language
of divisors at all. As it is easy to see, one can replace a divisor on $ \mathring{D} $
by a change of the gluings of $ {\mathcal L}_{\alpha}|_{\gamma_{i}} $ and $ {\mathcal L}_{\alpha}|_{\gamma_{i}'} $. While divisors at
infinity cannot be translated to a similar change of the gluing, the
necessary generalizations of our results are trivial, as far as the
divisor at infinity is has finite support. \end{remark}

\subsection{Three theories }\label{s01.60}\myLabel{s01.60}\relax  Unfortunately, for $ \alpha\not=0 $ there is no naturally
defined notion of a {\em complex-analytic\/} section of $ {\mathcal L}_{\alpha} $. However, for
$ \alpha\in\frac{1}{2}{\mathbb Z} $, one can define similar growth conditions on sections of the
line bundle $ \omega^{\alpha} $, which {\em is\/} complex-analytic, and the action of $ \operatorname{SL}\left(2,{\mathbb C}\right) $ on
this bundle is ``very similar'' to the action on $ \left(\omega\otimes\bar{\omega}\right)^{\alpha/2} $. Additionally,
the Hilbert norms allow limits when $ \alpha $ goes to 0 or 1; while the limits
are not positive-definite, they induce Hilbert norms on an appropriate
subspace (or a quotient space) of codimension 1. This gives 3
satisfactory theories: $ \alpha\in\left\{0,1/2,1\right\} $.

In the case $ g\left(C^{*}\right)<\infty $, the defined above ``reference bundles'' $ \widetilde{{\mathcal L}}_{\alpha} $ on $ C^{*} $
in these three cases are $ {\mathcal O} $, $ \omega^{1/2} $, and $ \omega $; the degrees are 0, $ g-1 $, and
$ 2\left(g-1\right) $. By our definition of a line bundle on $ C^{*} $, there is a well-defined
notion of ``its relative degree'' w.r.t.~$ \widetilde{{\mathcal L}}_{\alpha} $. Thus in the case $ g\left(C^{*}\right)=\infty $ we
get three theories: one $ \left(\alpha=0\right) $ with a well-defined degree $ d $ of a line
bundle; another $ \left(\alpha=1/2\right) $ with a well-defined difference $ d-g $; the third
$ \left(\alpha=1\right) $ with a well-defined difference $ d-2g $.

Out of these three, $ \alpha=1/2 $ gives the most interesting, self-dual
theory. Moreover, in this case the geometrically defined gluing of $ {\mathcal L}_{\alpha}|_{\gamma} $
and $ {\mathcal L}_{\alpha}|_{\gamma'} $ automatically ``adds a pole'' for each glued pair of circles
(thus each added handle), similar to pole at $ Q $ in Section~\ref{s02.30}.
However, the corresponding norm on $ H^{1/2}\left(C,\omega^{1/2}\right) $ is non-local (as all
Sobolev norms with fractional index), and is conformally-invariant only
``approximately''.

\subsection{Case $ \alpha=0 $ }\label{s01.70}\myLabel{s01.70}\relax  In this paper of the series, we consider the case $ \alpha=0 $
only. This simplifies the discussion, since the principal objects are
functions, not (fractional-degree) differential forms. Moreover, the
norms we consider are manifestly conformally-invariant, and defined by
local formulae.

The drawbacks are: first, the duality theorem needs to be postponed
until we consider the case $ \alpha=1 $. Second, to get a positive-definite
invariant norm, we need to consider quotient-spaces of functions by
constants (the constants are in the kernel of the ``norm''). Third, the
auto-gluing in the case $ \alpha=0 $ differs a lot from the ``small deformation'' of
Section~\ref{s02.30}. Essentially, it corresponds to consideration of $ {\mathcal L}^{\varepsilon} $
without the added pole at $ Q $; to get a ``small deformation'' (thus a hope to
get a finite dimension of ``global sections'' after an infinite number of
such steps), one needs to add some one-dimensional ``slack'', something
similar to allowing a pole at $ Q $ in Section~\ref{s02.30}.

By robustness, it is not very important {\em which\/} slack we allow instead
of allowing a pole. Allowing a pole at $ Q $ is equivalent to replacing $ \bar{\partial}F=0 $
by $ \bar{\partial}F = c\delta_{Q} $; here $ c\in{\mathbb C} $, $ \delta_{Q} $ is the $ \delta $-function at $ Q $. Since we need to
consider functions up to a constant anyway, it makes sense to allow a
jump by an additive constant when we glue $ \partial R_{i} $ and $ \partial R_{i}' $. Thus the last two
inconveniences of the case $ \alpha=0 $ partially compensate each other.

This makes starting our consideration with the case $ \alpha=0 $ very
convenient: we can introduce principal concepts without unnecessary
complications. Moreover, this case is needed anyway for the general
theory: it is dual to the consideration of the {\em partial period mapping\/}
$ \Gamma\left(M,\omega\right) \to {\mathbb C}^{g} $ of taking periods of global holomorphic forms along
$ A $-cycles on the Riemann surface. Until the remaining papers of the
series appear, an interested reader can refer to \cite{Zakh97Qua}, where all
three cases $ \alpha\in\left\{0,1/2,1\right\} $ are considered (though with much stricter
assumptions than in this paper).

\subsection{The principal results } In this paper of the series we define the
principal notions: an infinite-genus curve and a line bundle on it (with
$ \alpha=0 $ only), and discuss only the simplest possible properties of these
objects: the Riemann--Roch theorem. The principal result is Theorem~%
\ref{th33.35}, (and the amplifications in Sections~\ref{s3.45},~\ref{s3.55} and~\ref{s6.6})
which show that for validity of the Riemann--Roch theorem the only
condition is that the removed disks $ R_{i} $, $ R_{i}' $ (notations of Section~%
\ref{s02.40}) are ``small'' enough. (The formalization of the latter notion is
the notion of conformal distance, see Definition~\ref{def31.40}.)

The striking corollary of this fact is that there is no restriction
on the {\em position\/} of the disks, only on their ``sizes''. In particular, dust
can take arbitrary large proportion of the whole curve $ C^{*} $; not excluding
the case when the $ C^{*} $ consists of dust only (no smooth points on $ C^{*} $ at all).

This is the principal difference of the approach of this paper to
one in \cite{Zakh97Qua}: the much stricter conditions of \cite{Zakh97Qua} required
an annulus of smooth points around each cycle $ \gamma_{i} $, $ \gamma_{i}' $. The other
significant difference is that we allow gluing the curve $ C^{*} $ out of a
(possibly infinite) collection of curves $ C_{\left[k\right]} $ (with removed regions
$ R_{k,i} $). This allows, e.g., a uniform consideration of curves of Section~%
\ref{s02.40} together with\footnote{Note that while we expect our conditions of Riemann--Roch theorem to be
close to optimal when we glue $ C^{*} $ out of one curve $ C $, they must be very
non-optimal in the case of gluing of pants.} (more traditional) curves glued out of an infinite
collection of {\em pants\/} (spheres with 3 disks removed); see also Section~%
\ref{s6.70} for the example of yet another useful type of curves.

Let us sketch the principal ingredients of our presentation: curves
we consider are glued of {\em model domains\/} (Section~\ref{s35.30}); analytic
functions are replaced by {\em Sobolev-holomorphic functions\/} on model domains
(Section~\ref{s2.45}); the gluing rules are introduced in Section~\ref{s30.10}.
Sections~\ref{s10.70} and~\ref{s10.80} introduce the translation rules from the
``usual'' $ \left(g<\infty\right) $ Riemann--Roch theorem to the language of
Sobolev-holomorphic functions and gluing data.

Sections~\ref{s4.11},~\ref{s40.20} and~\ref{s40.30} reduce the (translated)
Riemann--Roch theorem to a condition of almost-transversality for two
appropriately defined subspaces in the space of functions on the {\em smooth
part of the boundary\/} of the model domains. (This is one of the key
conceptual ingredients: since this smooth part of the boundary is
enumerated by a discrete collection of indices, this should be considered
as a kind of {\em discretization\/} of the initial problem.) The remaining part
of Section~\ref{h42} introduces the translation of this almost-transversality
condition to the estimates of the norms (or of the essential spectrum) of
certain infinite matrices. This section concludes by the simplest
possible effective form of the Riemann--Roch theorem.

Section~\ref{h50} introduces generalizations of this simplest form which
are needed to study models of curves appearing in integrable systems, as
well as those needed for general divisor--line-bundle correspondence; in
Section~\ref{s6.3} we show that the curves which satisfy the Riemann--Roch
theorem may be of purely fractal nature. We also discuss special
properties enjoyed by the bundle $ {\mathcal O} $.

Finally, in the appendix (Section~\ref{h4}) we prove (and discuss the
motivations) for the particular form of Fredholm theorem used in
our treatment of almost transversality.

Note that the most of the statements of this paper are technically
straightforward; thus the motivations and heuristics may be as important
as the particular formulations of statements. Let us list less
straightforward technical statements forming the foundation of our
methods: Theorem~\ref{th35.149} (and Theorem~\ref{th35.150}) allow the
``discretization'' mentioned above; Lemma~\ref{lm131.32} allows application of
the above statements to real curves inside complex curves (as well as the
generalization of principal results to the case of many-Jordan-curves
boundary in
Section~\ref{s6.6}); Lemma~\ref{lm30.70} translates bundle gluing data to settings
of Section~\ref{h4} (in the context similar to Segal--Wilson's Universal
Grassmannian \cite{SegWil85Loo,PreSeg86Loo}). Applying this translation for
infinitely many curves requires a stronger ``discretization'' condition:
{\em fatness\/} (see Section~\ref{s3.55}); it is achieved in Theorem~\ref{th35.157} in
quasi-circular case, estimates of Section~\ref{s3.55} show that in ``most of
the cases'' fatness follows from the other assumptions of this paper.
Theorem~\ref{th41.20} introduces examples of foam curves.

\subsection{Moduli spaces and the Universal Grassmannian } Let us also mention the
natural problems which we {\em do not\/} discuss in this part of the series.
First of all, the duality theorem requires consideration of two values $ \alpha $,
$ \alpha' $ with $ \alpha+\alpha'=1 $; in this paper $ \alpha=0 $, so we do not consider the duality here
(see \cite{Zakh97Qua} instead---with more assumptions than we require in this
paper).

Second, the principal unit we consider here is a {\em model}, i.e., a
curve together with its representation via gluing of finite-genus pieces.
We do not consider the question when two models define ``the same'' curve.\footnote{There are two possible reasons why different gluing data may define the
same curve: to construct a model, one needs to {\em cut\/} the curve along a
collection of cycles, and choose identification of the pieces with
subsets of compacts curves. It is relatively easy to describe the
possible ambiguities {\em given a fixed choice\/} of the homotopy classes of cuts.
\endgraf
However, choosing different homotopy classes of cuts may lead to a
``significantly different'' model of a curve. However, recall the
conjecture of \cite{Zakh97Qua}: if two different choices of cuts both lead to
the gluing data satisfying the conditions of Theorem~\ref{th33.35}, then they
differ only for a finite number of cuts. The heuristic for this
conjecture is that the cuts which are {\em small cycles across\/} thin long
handles. Given infinitely many handles which are thinner and thinner, there
is essentially no choice: all the ``other'' cycles are going to be too
long.}
In other words, here we treat the question how big is the collection of
curves corresponding to points of the moduli space, not what is exactly a
{\em point\/} of the moduli space. Recall, however, that \cite{Zakh97Qua} introduces
the mapping of curves with a distinguished ``quasi-smooth'' point to the
Sato's Universal Grassmannian (\cite{SatSat82Sol,PreSeg86Loo}), which leads
to the notion of ``sameness'' for two models. The same mapping is still
defined for the curves we consider in this paper, this this approach
works for the foam curves too.

Another topic missing in this paper is the divisor--line-bundle
correspondence. It is more or less trivial to define $ {\mathcal L}\left(D\right) $ for a line
bundle $ {\mathcal L} $ and a finite divisor $ D $; similarly, one can do the same for
divisors with infinite support, as far as points with multiplicity $ -1 $ are
close enough to points of multiplicity 1. (This is similar to how
\cite{Zakh97Qua} uses the results of Section~\ref{s3.60} to show that any ``bounded''
line bundle of degree 0 may be realized using constant gluing function $ \psi_{j} $
of Section~\ref{s40.20}.) However, the ``interesting'' theory would work with
divisors $ D $ of infinite degree; in such cases $ \alpha $ for $ {\mathcal L}\left(D\right) $ is different from
$ \alpha $ for $ {\mathcal L} $. Again, we cannot discuss {\em this\/} topic until we have theories
suitable for different values of $ \alpha $.

\subsection{Historic remarks } The roots of this paper go back to Yu.~I.~Manin's
seminars of the spring of 1981 (see \cite{Zakh97Qua} for details), as well as
McKean and Trubowitz work \cite{McKTru76Hil} on the hyper-elliptic case.
During the last several years Feldman, Kn\"orrer and Trubowitz made a major
breakthrough using an unrelated approach (cf.~\cite{FelKnoTru96Inf}). The
relation of these curves to what we discuss here is explained in Section~%
\ref{s6.70}.

Numerous alternative approaches to curves of infinite genus exist,
both in rigorous settings, and in papers written using physical level of
arguments. One can break them into two different categories: one is
restricted to curves which allow a well-behaved finite-sheet covering
over $ {\mathbb C}{\mathbb P}^{1} $ (ramified at infinitely many points, and with a ``significant
singularity'' over $ \infty\in{\mathbb C}{\mathbb P}^{1} $); these curves are similar to hyperelliptic
curves. Another deals with curves similar to those in
our settings, but with severe restriction on the dust (e.g., with the
dust which consists of one point, or a finite number of points, only).
Such curves appear in study of, e.g., double-periodic solutions of KP
equation.

The first approach is developed in \cite{Schmi96Int,McKVan97Act,%
MulSchmiSchra98Hyp,Bik00Rie}, and \cite{GesHol00Dar}. The paper dealing
with second approach are \cite{Kosh89Cur}, which describes the period matrix
(this question is dual to what we investigate in this paper), \cite{Mer00Rie}
and \cite{Mer99Asy}, which investigate the Riemann--Roch problem, the
Jacobian, and the divisor-bundle correspondence; the book
\cite{GiesKnoTru93Geo} contains comprehensive studies of the geometry of
infinite-genus curves which appear in particular problems of mathematical
physics. As in the case of curves in \cite{FelKnoTru96Inf}, these particular
curves are very special cases of the curves we consider in this paper (as
well as in \cite{Zakh97Qua}). The paper \cite{Jin96Max} describes the class of
curves which cannot be extended (so are analogues of compact curves). In
the paper \cite{Ichi97Schot} holomorphic forms on a curve with a Schottky
model are studied as Taylor series of parameters of this model. (The
relation of curves studied in \cite{Jin96Max} and \cite{Ichi97Schot} and the curves
we consider in this paper is not yet clear.)

The author is most grateful to I.~M.~Gelfand, A.~Givental,
A.~Goncharov, D.~Kazhdan, M.~Kontsevich, Yu.~I.~Manin, H.~McKean,
V.~Serganova, A.~Tyurina and members of A.~Morozov's seminar for
discussions which directed many approaches applied here.

\section{Preliminaries }\label{h35}\myLabel{h35}\relax 

\subsection{Notations }\label{s35.24}\myLabel{s35.24}\relax  Consider a collection of topological vector spaces
$ V_{i} $, $ i\in I $. Then $ \prod_{i\in I}V_{i} $ denotes the space of all collections $ \left(v_{i}\in V_{i}\right)_{i\in I} $ with
the projective limit topology; $ \bigoplus_{i\in I}V_{i} $ denotes the space of collections
with a finite number of non-0 terms with the inductive limit topology.
For a collection with $ V_{i}\subset V $, $ \sum V_{i} $ denotes the closure of the image of $ \bigoplus V_{i} $
in $ V $. If all $ V_{i} $ are Hilbert spaces, $ \bigoplus_{l_{2},i\in I}V_{i} $ denotes the subspace of
$ \prod_{i\in I}V_{i} $ consisting of collections $ v=\left(v_{i}\in V_{i}\right)_{i\in I} $ with $ \|v\|^{2}\buildrel{\text{def}}\over{=}\sum_{i}\|v_{i}\|^{2}<\infty $.

We say that a topological vector space $ V $ has a {\em Hilbert topology}, if
there is a Hilbert norm on $ V $ which induces the topology of $ V $. We say that
two subspaces $ S_{1},S_{2}\subset V $ are {\em comparable\/} if $ S_{1}\cap S_{2} $ has a finite codimension in
both $ S_{1} $ and $ S_{2} $; then $ \operatorname{reldim}\left(S_{1},S_{2}\right) $ is the difference of these
codimensions. Two subspaces $ S_{1},S_{2}\subset V $ are {\em quasi-complementary\/} if $ \dim 
S_{1}\cap S_{2}<\infty $ and $ \operatorname{codim} S_{1}+S_{2}<\infty $. The {\em excess\/} of a quasi-complementary pair of
subspaces is $ \dim  \left(S_{1}\cap S_{2}\right) - \operatorname{codim}\left(S_{1}+S_{2}\right) $.

Since we need to consider both Sobolev spaces and spaces of
cohomology, we reserve the letter $ H $ for the former, and $ {\mathbit H} $ for the latter.

In what follows we define several flavors of {\em distortions}. List them
here for reference purposes: $ \Delta\left(E\right) $ in Section~\ref{s2.50}, $ \Delta\left(\gamma\right) $ in Section~%
\ref{s2.3} (used in Definition~\ref{def27.20} of quasi-circularity),
$ \Delta\left(E,E',\gamma,\gamma',\varphi,\psi\right) $ in Section~\ref{s4.11}, and $ \Delta\left(\varphi\right) $ in Section~\ref{s27}. To simplify
the discussion, most of the time we assume that these numbers are
uniformly bounded (otherwise we would need to incorporate them into the
estimates). Note that in most important particular cases these numbers
are all 1.

\subsection{The model domains }\label{s35.30}\myLabel{s35.30}\relax  The complex curves (or their
generalizations) we consider here are going to be glued of several pieces
$ D_{k} $, $ k\in K $, each piece being a closed subset of a compact complex curve $ C_{k} $.
The subsets $ D_{\bullet} $ we consider here have the boundary consisting of the
smooth part,\footnote{Section~\ref{s6.6} introduces modifications allowing Jordan curves as
components of the boundary.} which may have infinitely many connected components, and of
the accumulation points of these components. Enumerate the smooth
components of the boundaries of all the pieces $ D_{k} $ as $ \gamma_{j} $, $ j\in J $. Given one
such component $ \gamma_{j} $, the gluing process associates to it another component
$ \gamma_{j'} $, and a smooth orientation-inverting identification $ \varphi_{j}\colon \gamma_{j} \to \gamma_{j'} $.

Describe the effect of these modifications when both sets $ K $ and $ J $
are finite. Each gluing either decreases the number of connected
components by 1 (if $ \gamma_{j,j'} $ were on different components), or increase the
genus of one of the components by 1. As a result, out of $ p=|K| $ connected
pieces of genera $ g_{k} $, with $ d_{k} $ components of boundary each, one can glue
one compact connected curve of genus $ g=\sum g_{k}-p+1+\sum d_{k}/2 $, $ |J|=\sum d_{k} $. Call the
collection $ D_{\bullet} $ the {\em model\/} of the resulting curve.

Assume for a moment that all $ g_{k}=0 $. Then there are two different ways
to increase $ g $: either by increasing $ p $ (assuming $ d_{k}\geq3 $), or by increasing
$ d_{k} $. On the other hand, gluing two pieces of genus 0 along a pair of
components of boundary gives a piece of genus 0 too; this decreases $ p $ by
1, and increases $ \sum d_{k} $ by 2. Consequently, it is possible to substitute an
increment of $ p $ by an increment of $ \sum d_{k} $; eventually, one can replace a
model by one with $ p=1 $ (and possibly large $ |J| $). As the motivations in
\cite{Zakh97Qua} show, this substitution can be made to work also when one has
much more information (``growth conditions at infinity'') attached to the
pieces too.

In some sense, what we do in this paper is the formalization of the
last remark, if ``the growth conditions at infinity'' are understood as
``the growth compatible with $ H^{1} $-Sobolev smoothness''.

\begin{remark} One can make similar arguments that one can cut a piece with
$ d_{1}+d_{2}-2 $ components of boundary into two pieces with $ d_{1} $ and $ d_{2} $ components
of boundary correspondingly. However, as the examples of Section~\ref{s6.3}
show, this argument works only if $ d_{k} $ is finite: some pieces with an
infinite number of ``holes'' cannot be cut into an infinite number of
pieces with a finite number of ``holes'' each. Thus consideration of pieces
with infinitely many ``holes'' leads to new effects which cannot be
described by gluing together simpler pieces.

Due to these new effects, and a possibility to replace many ``simple''
pieces by one ``complicated'' piece, the pieces with an infinite number of
holes are especially important for us. Up to Section~\ref{s2.45}, we
describe what are ``holomorphic functions'' on such
pieces.\footnote{Section~\ref{s54} describes additional conditions on the pieces to imply the
Riemann--Roch theorem.} \end{remark}

\subsection{Generalized Sobolev spaces }\label{s2.40}\myLabel{s2.40}\relax  For our purposes we need to slightly
extend some standard notions of the theory of Sobolev spaces (compare
with \cite{Fried82Int,Rau91Par}). First, recall the notions which we use
without any modification.

The {\em Sobolev\/} $ s $-{\em norm}, $ s\in{\mathbb R} $, on smooth rapidly decreasing functions $ f\left(x\right) $
on $ {\mathbb R}^{n} $ is defined by $ \|f\|_{s}^{2} = \int|\widehat{f}\left(\xi\right)\left(1+|\xi|^{s}\right)|^{2}d\xi $ (here $ \widehat{f} $ is the Fourier
transform of $ f $). The completion w.r.t.~this norm gives a Hilbert space
$ H^{s}\left({\mathbb R}^{n}\right) $ called the {\em Sobolev space}. Any element of $ H^{s}\left({\mathbb R}^{n}\right) $ may be identified
with a generalized function $ f $ on $ {\mathbb R}^{n} $; these generalized functions are
those for which the Fourier transform $ \widehat{f} $ is locally-$ L_{2} $, and the integral
of the Sobolev $ s $-norm converges. Since multiplication by a smooth
function with a compact support is a continuous operator in $ H^{s}\left({\mathbb R}^{n}\right) $, it
makes sense to consider the generalized functions which are {\em locally
Sobolev}, i.e., become Sobolev after a multiplication by any smooth
function with a compact support. While the finiteness of the Sobolev norm
reflects both the degree of smoothness of the function, and its decay at
infinity, the property of being locally Sobolev reflects the smoothness
only.

Locally Sobolev functions form a vector space $ H_{\text{loc}}^{s}\left({\mathbb R}^{n}\right) $; it has a
natural topology (of the appropriate inverse limit). Moreover, this
topological vector space is invariant w.r.t.~diffeomorphisms of $ {\mathbb R}^{n} $. This
makes it possible to define the topological vector space $ H_{\text{loc}}^{s}\left(M\right) $ for any
manifold $ M $; a generalized function $ f $ on $ M $ belongs to $ H_{\text{loc}}^{s}\left(M\right) $ if for any
coordinate chart $ {\mathbb R}^{n}\supset V \xrightarrow[]{\varphi} U\subset M $ and any smooth function $ \psi $ on $ U $ the
(generalized) function $ \varphi^{*}\left(\psi f\right) $ is locally Sobolev on $ {\mathbb R}^{n} $. Similarly, one
can define the topological vector space\footnote{See Section~\ref{s35.24} on how we avoid the conflict with cohomology.} $ H_{\text{loc}}^{s}\left(M,{\mathcal L}\right) $ of locally Sobolev
sections of an arbitrary finite-dimensional vector bundle $ {\mathcal L} $ on $ M $.

The vector space $ H_{\text{loc}}^{s}\left(M\right) $ has a natural topology of the inverse
limit. Moreover, if $ M $ is compact, then, as it is easy to see, this
topology is equivalent to the topology given by a Hilbert norm. In such a
case we use notation $ H^{s}\left(M\right) $ instead of $ H_{\text{loc}}^{s}\left(M\right) $.

Multiplications by smooth functions and diffeomorphisms induce
continuous operators in $ H^{s}\left(M\right) $. If $ \psi f=0 $ for an appropriate smooth function
$ \psi $ such that $ \psi\left(m\right)\not=0 $, one says that $ f $ {\em vanishes near\/} $ m\in M $. The {\em support\/} $ \operatorname{Supp}
f\subset M $ consists of points $ m\in M $ such that $ f $ does not vanish near $ m $; it is a
closed subset of $ M $.

The only deviation from the ``classical'' terminology is in the following

\begin{definition} Consider a subset $ U $ of the manifold $ M $. Let $ \mathring H^{s}\left(U\right) $ denote
the closure of the vector subspace $ \left\{f\in H^{s}\left(M\right) \mid \operatorname{Supp} f\subset U\right\} $.

Consider $ V\subset M $. Let $ H^{s}\left(V\subset M\right)=H^{s}\left(M\right)/\mathring H^{s}\left(M\smallsetminus V\right) $. \end{definition}

The deviation from the standard definition is that we do not require
that $ U $ is closed, and $ V $ is open. Obviously, $ \mathring H^{s}\left(U\right)\subset\mathring H^{s}\left(\bar{U}\right) $, but this
inclusion may be proper, as the example below shows. Call spaces $ \mathring{H}\left(U\right) $
and $ H\left(V\subset M\right) $ the {\em generalized Sobolev spaces}.

\begin{example} Consider a disjoint\footnote{There are simpler examples, but this one gives a domain we are
going to deal with; see Section~\ref{s6.3}.} family of open subsets $ V_{i}\subset M $, Let
$ {\mathcal V}=\overline{\bigcup V_{i}}\smallsetminus\bigcup V_{i} $. Suppose that
\begin{enumerate}
\item
The natural mapping $ \bigoplus_{l_{2}}\mathring H^{s}\left(V_{i}\right) \xrightarrow[]{\iota} H^{s}\left(M\right) $ is a (continuous)
monomorphism\footnote{I.e., the image is closed, and the mapping is an isomorphism onto the
image.}.
\item
There is a function $ g\in H^{s}\left(M\right) $ with $ \operatorname{Supp} g\subset{\mathcal V} $.
\end{enumerate}
The first condition insures that the space $ \mathring H^{s}\left(\bigcup V_{i}\right) $ is the image of
the mapping $ \iota $. Hence any non-zero function $ f\in\mathring H^{s}\left(\bigcup V_{i}\right) $ satisfies the
condition $ \operatorname{Supp} f \cap \bigcup V_{i}\not=\varnothing $.

Thus the function $ g\in H^{s}\left(M\right) $ with $ \operatorname{Supp} g\subset{\mathcal V} $ satisfies $ f\in\mathring H^{s}\left(\overline{\bigcup V_{i}}\right) $,
but $ f\notin\mathring H^{s}\left(\bigcup V_{i}\right) $. One of the standard facts of the theory of Hausdorff
dimension is that the second condition is satisfied if $ \dim _{\text{Hausdorff}}{\mathcal V}>\dim 
M+2s $. Moreover, it works if $ s=0 $ and $ {\mathcal V} $ has a positive measure. \end{example}

In Section~\ref{s6.3} we show how to construct a family of disks $ V_{i} $ which
satisfy the first condition. The centers of these disks may be an
arbitrary locally discrete set $ \Delta $. Moreover, one can find $ \Delta $ such that the
corresponding set $ {\mathcal V} $ does not depend on radii; any set with an empty
interior can be obtained as such $ {\mathcal V} $. Thus the above construction works for
any $ s\leq0 $ (in this part of the series we are most interested in the case
$ s=0 $).

On the other hand, if $ U $ has smooth boundary, then $ \mathring H^{s}\left(U\right)=\mathring H^{s}\left(\bar{U}\right) $.

\subsection{Norms on Sobolev spaces }\label{s1.30}\myLabel{s1.30}\relax  As explained above, the Sobolev space of
sections of a line bundle is defined up to topological equivalence only,
it has no canonical Hilbert norm. However, in application we will need to
consider Hilbert direct sums, which require a specification of the
{\em Hilbert norm\/} (as opposed to {\em Hilbert topology}, see Section~\ref{s35.24}). As
explained in Section~\ref{s02.50}, the spaces $ H^{1-\alpha}\left(C,\omega^{\alpha}\right) $ have a canonically
defined norm (or an appropriate approximation); gluing $ \omega^{\alpha}|_{D_{k}} $ for pieces
$ D_{k}\subset C_{k} $ allows a definition of a norm on the space of global sections via
the Hilbert direct sum in $ k $. In the case we consider here, $ \alpha\in\left\{0,1\right\} $, and
the norm {\em is\/} canonically defined.

Given a compact oriented surface $ C $ with a conformal structure,
denote by $ \omega_{\perp} $ the $ 90^{\circ} $-counterclockwise rotation of a section $ \omega $ of $ \Omega^{1}\left(C\right) $.
Then $ \|\omega\|^{2}\buildrel{\text{def}}\over{=}\int_{C}\omega\omega_{\perp} $ gives a canonically defined Sobolev norm on $ H^{0}\left(C,\Omega^{1}\right) $.
Restricting to the components $ \omega $, $ \bar{\omega} $ of $ \Omega^{1}\otimes{\mathbb C}=\omega\oplus\bar{\omega} $, one obtains the norms
$ \|\alpha\|^{2}=-\frac{i}{2}\int_{C}\bar{\alpha}\alpha $ on $ H^{0}\left(C,\bar{\omega}\right) $ and $ \|\alpha\|^{2}=\frac{i}{2}\int_{C}\bar{\alpha}\alpha $ on $ H^{0}\left(C,\omega\right) $.

While $ H^{1}\left(C,{\mathcal O}\right) $ does not carry a natural norm, for a compact $ C $ the
mapping $ \bar{\partial}\colon H^{1}\left(C,{\mathcal O}\right) \to H^{0}\left(C,\bar{\omega}\right) $ identifies $ H^{1}\left(C,{\mathcal O}\right)/\operatorname{const} $ with a closed
subspace of $ H^{0}\left(C,\bar{\omega}\right) $. This provides a canonically defined norm on
$ H^{1}\left(C,{\mathcal O}\right)/\operatorname{const} $.

\begin{remark} \label{rem1.305}\myLabel{rem1.305}\relax  One could define a similar norm using the operator $ \partial $
instead; however, the result is going to be the same due to
\begin{equation}
\int\overline{\partial f}\partial f=\int\bar{\partial}\bar{f}\partial f= -\int\bar{f}\bar{\partial}\partial f=\int\bar{f}\partial\bar{\partial}f=-\int\partial\bar{f}\bar{\partial}f=-\int\overline{\bar{\partial}f}\bar{\partial}f;
\notag\end{equation}
this identity comes very handy in Section~\ref{s5.7}. In particular, the norm of
$ f\in H^{1}\left(C,{\mathcal O}\right)/\operatorname{const} $ is $ \|df\|/\sqrt{2} $.

The norm on $ H^{1}\left(C,{\mathcal O}\right)/\operatorname{const} $ induces canonically defined norms on
$ H^{1}\left(D\subset C,{\mathcal O}\right)/\operatorname{const} $ and on $ H^{0}\left(D\subset C,\bar{\omega}\right) $. Note that the Hilbert norm on $ H^{0}\left(D\subset C,\bar{\omega}\right) $
coincides with $ \|\alpha\|^{2}=-\frac{i}{2}\int_{D}\bar{\alpha}\alpha $. \end{remark}

Consider a subset $ D\subset C $ such that $ 1\notin\mathring{H}^{s}\left(D\right)\subset H^{s}\left(C\right) $, here 1 is considered
as an element of $ H^{s}\left(C\right) $. Since $ \mathring{H}^{s}\left(D\right) $ is closed in $ H^{s}\left(C\right) $, the mapping $ \mathring{H}^{s}\left(D\right)
\hookrightarrow H^{s}\left(C\right)/\operatorname{const} $ is a monomorphism, thus a norm on $ H^{s}\left(C\right)/\operatorname{const} $ induces a
norm on $ \mathring{H}^{s}\left(D\right) $.

\begin{lemma} Consider a connected compact complex curve $ C $. Then $ 1\notin\mathring{H}^{s}\left(D\right) $ if
$ C\smallsetminus D $ has a non-empty interior; or if $ C\smallsetminus D $ contains a smooth curve and
$ s>1/2 $; or if $ C\smallsetminus D $ contains a connected component which is not a point, and
$ s=1 $; or if $ C\smallsetminus D $ is non-empty and $ s>1 $. \end{lemma}

\begin{proof} All the statements except the last but one follow from the
continuity properties of the restriction to submanifolds. The remaining
statement is equivalent to the following statement: Let $ \gamma $ be a connected
subset of a complex curve $ C $, and $ \gamma $ is not a point. Consider a sequence
$ \left(\psi_{k}\right) $ of smooth functions on $ C $, and a sequence $ U_{k} $ of neighborhoods of $ \gamma $
such that $ \psi_{k}|_{U_{k}}=0 $. It is enough get a contradiction with $ \psi_{k} \to 1 $ in
$ H^{1}\left(C\right) $.

Since multiplication by a smooth function is continuous in $ H^{1}\left(C\right) $, it
is enough to find $ \varepsilon $ such that $ \|\Psi\left(1-\psi_{k}\right)\|_{H^{1}\left(C\right)}>\varepsilon $ for any $ k $; here $ \Psi $ is an
appropriate cut-off function. This makes the question local on $ C $, so we
may assume $ C={\mathbb C}{\mathbb P}^{1} $; moreover, it is enough to show that
$ \int_{C\smallsetminus\gamma}\|d\left(\Psi\left(1-\psi_{k}\right)\right)\|^{2}>\varepsilon $; this is formulated completely in terms of
holomorphic geometry of $ C\smallsetminus\gamma $. Thus we may assume that $ C\smallsetminus\gamma $ is a unit disk.
Now one can apply the first part of the lemma. \end{proof}

This provides a canonically defined norm on $ \mathring{H}^{1}\left(D\right) $ if $ C\smallsetminus D $ contains a
Jordan curve. In what follows we use the canonically defined norms above
unless specified otherwise.

\section{Sobolev holomorphic functions }

\subsection{Sobolev-holomorphic functions and decomposition for $ g=0 $ }\label{s2.45}\myLabel{s2.45}\relax 

Differential operators act on Sobolev spaces decreasing $ s $ by the
degree of the operator, and do not increase the support. Thus given an
element $ f $ of $ H^{1}\left(D\subset C,{\mathcal O}\right) $, $ \bar{\partial}f $ is a correctly defined element of $ H^{0}\left(D\subset C,\bar{\omega}\right) $.

\begin{definition} Given a closed subset $ D $ of a compact complex curve $ C $, an
$ H^{1} $-{\em holomorphic function\/} on $ D $ is an element $ f $ of $ H^{1}\left(D\subset C,{\mathcal O}\right) $ which satisfies
the condition $ \bar{\partial}f=0\in H^{0}\left(D\subset C,\bar{\omega}\right) $. Denote the the space of $ H^{1} $-holomorphic
functions on $ D\subset C $ by $ {\mathcal H}^{1}\left(D\subset C\right) $ (or just $ {\mathcal H}^{1} $). \end{definition}

Note that the Sobolev spaces in this definition are generalized
ones. The norm on $ H^{1}\left(D\subset C\right)/\operatorname{const} $ induces a canonically defined norm on
$ {\mathcal H}^{1}\left(D\subset C\right)/\operatorname{const} $.

\begin{remark} The heuristic on the use of ``generalized Sobolev'' space is that
$ f\in H^{1}\left(D\subset C,{\mathcal O}\right) $ is an equivalence class modulo functions with support
``inside'' $ C\smallsetminus D $; in other words, we ``keep'' the information about
$ f|_{\partial D} $. The equation $ \bar{\partial}f=0 $ is again satisfied only up to functions with
support ``inside'' $ C\smallsetminus D $; in other words, the equations should be satisfied
``also'' on $ \partial D $. \end{remark}

Consider $ D $ and $ C $ as in the definition above. From now on assume that
$ D\not=C $. Since $ C\smallsetminus D $ is open, it may be represented as a disjoint union of open
connected sets $ R_{j}\subset C $, $ j\in J $; here $ J $ is an appropriate set of indices. For
each $ j\in J $ let $ D_{j}=C\smallsetminus R_{j} $, $ D\subset D_{j}\subset C $. Then $ D=\bigcap_{j}D_{j} $. Consider the restriction
mapping $ {\mathcal H}^{1}\left(D_{j}\subset C\right) \to {\mathcal H}^{1}\left(D\subset C\right) $, and the induced mapping of quotient spaces
$ {\mathcal H}^{1}\left(D_{j}\subset C\right)/\operatorname{const} \to {\mathcal H}^{1}\left(D\subset C\right)/\operatorname{const} $. Taken together for every $ j\in J $, these
mappings define a mapping $ \widetilde{\rho}\colon \bigoplus{\mathcal H}^{1}\left(D_{j}\subset C\right)/\operatorname{const} \to {\mathcal H}^{1}\left(D\subset C\right)/\operatorname{const} $.

\begin{theorem} \label{th35.149}\myLabel{th35.149}\relax  If $ C $ is of genus 0, the described above mapping $ \widetilde{\rho} $
extends continuously to a Fredholm mapping $ \rho\colon \bigoplus_{l_{2}}{\mathcal H}^{1}\left(D_{j}\subset C\right)/\operatorname{const} \to
{\mathcal H}^{1}\left(D\subset C\right)/\operatorname{const} $.
In fact $ \rho $ is a natural unitary mapping of Hilbert spaces. \end{theorem}

\begin{proof} Consider a complex curve $ C $ of arbitrary genus. Consider the
mapping $ \bar{\partial}\colon H^{1}\left(C,{\mathcal O}\right) \to H^{0}\left(C,\bar{\omega}\right) $. Since $ \bar{\partial} $ is an elliptic operator of degree
1, it is Fredholm, and any function in $ \operatorname{Ker}\bar{\partial} $ is smooth, similarly for
$ \operatorname{Ker}\bar{\partial}^{*} $. As a corollary, $ \operatorname{Ker} \bar{\partial} $ is spanned by 1, and $ \operatorname{Ker} \bar{\partial}^{*} $ consists of
global holomorphic $ 1 $-forms. From now on assume $ C={\mathbb C}{\mathbb P}^{1} $. In particular,
$ \operatorname{Coker}\bar{\partial}=0 $.

Due to the conventions on norms from Section~\ref{s1.30}, $ \bar{\partial} $ induces a
unitary isomorphism $ H^{1}\left(C,{\mathcal O}\right)/\operatorname{const} \to H^{0}\left(C,\bar{\omega}\right) $. Given $ f\in{\mathcal H}^{1}\left(D\subset C\right) $, consider
two different liftings $ \widetilde{f}_{1} $, $ \widetilde{f}_{2} $ of $ f $ to elements of $ H^{1}\left(C,{\mathcal O}\right) $. By definition,
$ \widetilde{f}_{1}-\widetilde{f}_{2}\in\mathring{H}^{1}\left(C\smallsetminus D,{\mathcal O}\right) $. Moreover, $ \bar{\partial} \widetilde{f}_{1,2}\in\mathring{H}^{0}\left(C\smallsetminus D,\bar{\omega}\right) $. This implies that $ \bar{\partial}\widetilde{f}_{1} $ is a
canonically defined element of
\begin{equation}
{\mathcal R}_{D\subset C} \buildrel{\text{def}}\over{=} \mathring{H}^{0}\left(C\smallsetminus D,\bar{\omega}\right)/\bar{\partial} \mathring{H}^{1}\left(C\smallsetminus D,{\mathcal O}\right).
\notag\end{equation}
Since $ \bar{\partial} $ is surjective, the mapping $ {\mathcal H}^{1}\left(D\subset C\right)/\operatorname{const} \to {\mathcal R}_{D\subset C}\colon f+\operatorname{const} \mapsto \bar{\partial}\widetilde{f}_{1} $
is a unitary mapping of Hilbert spaces.

\begin{lemma} Consider disjoint open subsets $ R_{j} $, $ j\in J $, of a compact complex
curve $ C $. Let $ R=\coprod_{j}R_{j} $. Then $ \mathring{H}^{0}\left(R,\bar{\omega}\right) \simeq \bigoplus_{l_{2}}\mathring{H}^{0}\left(R_{j},\bar{\omega}\right) $, $ \mathring{H}^{1}\left(R,{\mathcal O}\right) \simeq
\bigoplus_{l_{2}}\mathring{H}^{1}\left(R_{j},{\mathcal O}\right) $, the former natural isomorphism is unitary, the latter is an
isomorphism of topological vector spaces. \end{lemma}

\begin{proof} By definition, an element $ f\in\mathring{H}^{s}\left(R\right) $ may be approximated by an
element $ f' $ of $ H^{s}\left(C\right) $ with $ \operatorname{Supp} f'\subset R $. Since $ \operatorname{Supp} f' $ is closed, $ \operatorname{Supp} f' $ is
compact, thus $ \operatorname{Supp} f' $ is contained in a finite union of several domains
$ R_{s} $. Thus $ \mathring{H}^{s}\left(R\right)\subset\sum_{j}\mathring{H}^{s}\left(R_{j}\right) $, which implies $ \mathring{H}^{s}\left(R\right)=\sum_{j}\mathring{H}^{s}\left(R_{j}\right) $. For $ s=0 $ the
subspaces $ \mathring{H}^{0}\left(R_{j}\right)\subset H^{0}\left(C\right) $ are obviously orthogonal, which proves one
statement of the lemma.

Similarly, $ \bar{\partial} $-images of $ \mathring{H}^{1}\left(R_{j},{\mathcal O}\right) $ in $ H^{0}\left(C,\bar{\omega}\right) $ are orthogonal; thus
images of $ \mathring{H}^{1}\left(R_{j},{\mathcal O}\right) $ in $ H^{1}\left(C,{\mathcal O}\right)/\operatorname{const} $ are orthogonal. One may assume that $ C $
is connected, and $ |J|>1 $, thus $ \mathring{H}^{1}\left(R,{\mathcal O}\right) $ projects monomorphically to
$ H^{1}\left(C,{\mathcal O}\right)/\operatorname{const} $. Since $ \mathring{H}^{1}\left(R_{j},{\mathcal O}\right) $ lie in $ \mathring{H}^{1}\left(R,{\mathcal O}\right) $, this proves the remaining
statement of the lemma. \end{proof}

As a corollary, we can see that $ {\mathcal R}_{D\subset C} $ can be naturally identified
with $ \bigoplus_{l_{2}}{\mathcal R}_{D_{j}\subset C} $, and that this identification is unitary. Applying the
same arguments to $ D_{j} $ instead of $ D $, one can see that $ {\mathcal H}^{1}\left(D\subset C\right)/\operatorname{const} $ is
isomorphic to $ \bigoplus_{l_{2}}{\mathcal H}^{1}\left(D_{j}\subset C\right)/\operatorname{const} $. Obviously, the restriction mapping
$ {\mathcal H}^{1}\left(D_{j}\subset C\right)/\operatorname{const} \to {\mathcal H}^{1}\left(D\subset C\right)/\operatorname{const} $ is compatible with this identification.
This finishes the proof of the theorem. \end{proof}

\subsection{Decomposition for an arbitrary genus }\label{s2.50}\myLabel{s2.50}\relax  To generalize the theorem
above to the case of $ C $ of arbitrary genus $ g\left(C\right) $, we need to compensate for
$ \bar{\partial}\colon H^{1}\left(C,{\mathcal O}\right) \to H^{0}\left(C,\bar{\omega}\right) $ being not surjective. Results of this section allow
to make the statements of this paper slightly more general, the price
being slightly more cumbersome formulations. The results of this paper
can be weakened by assuming that all the pieces we glue our curve from
are of genus 0; since these weaker results are still interesting, one can
skip this section on the first reading, assuming that during the gluing
of Section~\ref{s3.2} all the finite-genus pieces are of genus 0, and ignoring
all lower indices $ E $.

The arguments of Section~\ref{s2.45} presumed solvability of $ \bar{\partial}f=\alpha $. Since
$ \operatorname{Im}\bar{\partial}\subset H^{0}\left(C,\bar{\omega}\right) $ has codimension $ g $, for $ g>0 $ some compensation is needed.
Consider an arbitrary projection $ p\colon H^{0}\left(C,\bar{\omega}\right) \to \operatorname{Im}\bar{\partial} $; then for any
$ \alpha\in H^{0}\left(C,\bar{\omega}\right) $ the expression $ f=\bar{\partial}^{-1}\left(p\alpha\right) $ is a correctly defined element of
$ H^{1}\left(C,{\mathcal O}\right)/\operatorname{const} $. If $ E=\operatorname{Ker} p $, then $ f $ is a solution of $ \bar{\partial}f\equiv\alpha \mod E $.

This leads to the following definition:

\begin{definition} A subspace $ E\subset H^{0}\left(C,\bar{\omega}\right) $ is an {\em excess space\/} if the natural
pairing between $ E $ and the space $ \Gamma_{\text{an}}\left(C,\omega\right) $ of global holomorphic $ 1 $-forms is
non-degenerate. Given an excess space $ E $, let $ {\mathcal H}_{E}^{1}\left(D\subset C\right) $ consists of
$ f\in H^{1}\left(D\subset C\right) $ with $ \bar{\partial}f\in\pi E $, here $ \pi $ is the projection from $ H^{0}\left(C,\bar{\omega}\right) $ to $ H^{0}\left(D\subset C,\bar{\omega}\right) $.

For a subspace $ V\subset H^{0}\left(C,\bar{\omega}\right) $ and $ D\subset C $, let $ V_{D}=\left\{f\in V \mid \operatorname{Supp} f\subset D\right\} $. An excess
space $ E $ is $ D $-{\em supported\/} if $ E=E_{D} $. Given a collection $ R_{i} $, $ i\in I $, of disjoint
subsets of $ C $, and $ R=\overline{\bigcup R_{i}} $, an excess space $ E $ {\em is\/} $ R_{\bullet} $-{\em split\/} if $ E_{R}=\sum E_{\bar{R}_{i}} $.

Given an excess space $ E $, the {\em distortion\/} $ \Delta\left(E\right) $ is the norm of the
projector to $ \operatorname{Im}\bar{\partial} $ along $ E $. \end{definition}

Obviously, $ \dim  E=g\left(C\right) $. Since the particular choice of $ E $ is not
important for the following arguments, we use notation $ {\mathcal H}_{E}^{1} $ without
mentioning $ E $ otherwise. If $ C={\mathbb C}{\mathbb P}^{1} $, then $ {\mathcal H}_{E}^{1}={\mathcal H}^{1} $. Note that the spaces
$ {\mathcal H}_{E}^{1}/\operatorname{const} $ are equipped with natural Hilbert norms induced from $ H^{1}/\operatorname{const} $.

\begin{remark} \label{rem35.170}\myLabel{rem35.170}\relax  Obviously, excess subspaces exist. It is easy to find a
$ D_{\bullet} $-split one for any collection $ \left\{D_{\bullet}\right\} $ of disjoint open subsets. In fact,
given any closed subset $ D\subset C $ of non-zero measure, one can find a
$ D $-supported excess space $ E $. \end{remark}

\begin{theorem} \label{th35.150}\myLabel{th35.150}\relax  In notations of Section~\ref{s2.45}, suppose that $ E $ is
$ R_{\bullet} $-split. Then the mapping $ \widetilde{\rho}\colon \bigoplus{\mathcal H}_{E}^{1}\left(D_{j}\subset C\right)/\operatorname{const} \to {\mathcal H}_{E}^{1}\left(D\subset C\right)/\operatorname{const} $ can be
extended to an invertible mapping $ \rho\colon \bigoplus_{l_{2}}{\mathcal H}_{E}^{1}\left(D_{j}\subset C\right)/\operatorname{const} \to
{\mathcal H}_{E}^{1}\left(D\subset C\right)/\operatorname{const} $. Moreover, $ \rho $ is unitary if $ g\left(C\right)=0 $. \end{theorem}

\begin{proof} For a subspace $ V $ of $ H^{0}\left(C,\bar{\omega}\right) $ put $ V^{\vartheta}\buildrel{\text{def}}\over{=}V\cap\operatorname{Im}\bar{\partial} $; use the same
notation for subspaces of quotients of $ H^{0}\left(C,\bar{\omega}\right) $. Obviously, $ \bar{\partial} $ sends
$ {\mathcal H}^{1}\left(D\subset C\right)/\operatorname{const} $ to $ {\mathcal R}_{D\subset C}^{\vartheta} $. Similarly, $ {\mathcal H}_{E}^{1}\left(D\subset C\right)/\operatorname{const} $ is identified with
$ \left({\mathcal R}_{D\subset C}+\bar{E}\right)^{\vartheta} $, here $ \bar{E}\subset H^{0}\left(C,\bar{\omega}\right)/\bar{\partial}\mathring{H}^{1}\left(C\smallsetminus D,{\mathcal O}\right) $ is the projection of $ E\subset H^{0}\left(C,\bar{\omega}\right) $. It
is clear that $ \dim  \bar{E} = \dim  E $.

Since $ V^{\vartheta}\subset V\subset H^{0}\left(C,\bar{\omega}\right) $ is defined by equations $ \left< \alpha,v \right>=0 $, $ \alpha\in\Gamma_{\text{an}}\left(C,\omega\right) $,
the projection of $ \left(V+E\right)^{\vartheta} $ to $ V $ along $ E $ is an isomorphism if $ E\cap V=0 $. Thus
$ \left(V+E\right)^{\vartheta} $ is naturally identified with $ V/E^{V} $ via taking the quotient by $ E $;
here $ E^{V}=V\cap E $. This identification preserves the topology.

If $ W\subset H^{1}\left(C,{\mathcal O}\right) $, the same argument works for $ V\subset H^{0}\left(C,\bar{\omega}\right)/\bar{\partial}W $ (substituting
$ \bar{E} $ for $ E $), since the pairing with $ \Gamma_{\text{an}}\left(C,\omega\right) $ vanishes on $ \bar{\partial}W $. Fix $ R\subset C $; put
$ V_{R}=H^{0}\left(R,\bar{\omega}\right)/\bar{\partial}\mathring{H}^{1}\left(C\smallsetminus D,{\mathcal O}\right) $. Then $ \bar{E}^{V_{R}}=\bar{E}_{R} $. Thus $ \left({\mathcal R}_{D\subset C}+\bar{E}\right)^{\vartheta} $ is identified with
$ {\mathcal R}_{D\subset C}/\bar{E}_{C\smallsetminus D} $.

The splitness property can be restated as $ E_{C\smallsetminus D}=\bigoplus E_{\bar{R}_{i}} $. Now $ {\mathcal R}_{D\subset C}
=\bigoplus_{l_{2}}{\mathcal R}_{D_{i}\subset C} $ implies $ {\mathcal R}_{D\subset C}/E_{C\smallsetminus D}=\bigoplus_{l_{2}}{\mathcal R}_{D_{i}\subset C}/E_{\bar{R}_{i}} $, which finishes the proof. \end{proof}

\begin{remark} \label{rem35.180}\myLabel{rem35.180}\relax  The mapping $ \rho $ is not necessarily unitary if $ g\left(C\right)>0 $. The
non-unitary component of the identifications of the theorem is the
projection along $ \bar{E} $ from $ \left(V+\bar{E}\right)^{\vartheta} $ to $ V/\bar{E}^{V} $. Consequently, it is enough to
estimate the ``distortion'' of this projection; or the angles between $ \bar{E} $ and
$ \left({\mathcal R}_{D\subset C}+\bar{E}\right)^{\theta} $, and between $ \bar{E}/E_{C\smallsetminus D} $ and $ {\mathcal R}_{D\subset C}/E_{C\smallsetminus D} $.

If $ E $ is split w.r.t.~the collection $ \left(\left(R_{i}\right),D\right) $, then the latter
subspaces are orthogonal. Thus the degree of non-unitarity of $ \rho $ is
majorated by a function of the distortion $ \Delta\left(E\right) $. \end{remark}

\subsection{Riemann problem decomposition }\label{s2.3}\myLabel{s2.3}\relax  Recall that for a submanifold $ N\subset M $
of codimension $ d $ the restriction mapping $ H^{s}\left(M\right) \to H^{s-d/2}\left(M\right) $ is continuous
as long as $ s>d/2 $.

\begin{lemma} \label{lm131.32}\myLabel{lm131.32}\relax  Consider a compact smooth real curve $ \gamma $ which is a
submanifold in a real surface $ C $. Then the restriction mapping $ \widetilde{\rho}\colon H^{1}\left(C\right) \to
H^{1/2}\left(\gamma\right) $ induces an isomorphism $ H^{1}\left(\gamma\subset C\right)\simeq H^{1/2}\left(\gamma\right) $ of topological vector
spaces. \end{lemma}

\begin{proof} Since $ \widetilde{\rho}f=0 $ if $ \operatorname{Supp} f\cap\gamma=0 $; by continuity, $ \widetilde{\rho} $ vanishes on $ \mathring{H}^{1}\left(C\smallsetminus\gamma\right) $,
thus induces a continuous mapping $ \rho\colon H^{1}\left(\gamma\subset C\right) \to H^{1/2}\left(\gamma\right) $. Obviously, the
surjectivity and injectivity of $ \rho $ are local properties; by the invariance
of the Sobolev topology w.r.t.~diffeomorphisms, it is enough to
consider one particular curve $ \gamma $.

Assume $ \gamma=\partial D $, $ D $ being the unit disk in $ S^{2}={\mathbb C}{\mathbb P}^{1} $. Since $ \rho $ commutes with
the action of the group $ {\mathbb T} $ of rotations, it is enough to show that $ \rho $
induces uniformly bounded (from above and from below) isomorphisms
between the isotypical components of $ {\mathbb T} $. For $ k\geq0 $ let $ f_{k}\left(z\right)=z^{k} $ if $ |z|\leq1 $,
$ f_{k}\left(z\right)=\bar{z}^{-k} $ if $ |z|\geq1 $. For $ k<0 $ define $ f_{k}\left(z\right)=\bar{f}_{-k}\left(z\right) $. A simple calculation
shows that the $ H^{1} $-norms of $ f_{k} $ grow as $ |k|^{1/2} $, $ k\not=0 $; similarly for
$ \|f_{k}|_{\partial D}\|_{H^{1/2}} $. Thus $ \rho|_{\left<f_{k}\right>} $ is an isomorphism of topological vector spaces
(here $ \left<f_{k}\right> $ is the vector subspace spanned by $ f_{k} $), thus $ \rho $ is surjective.

To show injectivity, it is enough to consider isotypical
components one-by-one. A function from such a component can be written as
$ z^{n}\psi\left(z\right) $ near $ \partial D $, here $ \psi $ is rotation-invariant, thus $ \psi\left(z\right)=\psi\left(|z|\right) $. Thus
reduces the problem to the following $ 1 $-dimensional problem: show that a
smooth function $ \psi $ on $ \left[0.5,2\right] $ such that $ \psi\left(1\right)=0 $ can be $ H^{1} $-approximated by a
function which vanishes near 1. In turn, this is obvious. \end{proof}

\begin{definition} \label{def131.35}\myLabel{def131.35}\relax  Given $ f\in H^{1}\left(D\subset C,{\mathcal O}\right)/\operatorname{const} $, define $ \|f\|_{1,\text{int}} $ as
$ \|df\|_{L_{2}\left(D,\Omega^{1}\right)} $. Here $ d $ is de Rham differential. \end{definition}

It is clear that $ \|f\|_{1,\text{int}}\leq\sqrt{2}\|f\|_{H^{1}\left(D\subset C,{\mathcal O}\right)} $ and the coefficient
can be reduced to 1 if $ f\in{\mathcal H}^{1}\left(D\subset C,{\mathcal O}\right)/\operatorname{const} $.

\begin{lemma} \label{lm131.40}\myLabel{lm131.40}\relax  Consider a closed subset $ D\subset C $ of a compact real surface $ C $
with a smooth boundary. Then the norm $ \|\bullet\|_{1,\text{int}} $ induces a Hilbert space
structure on $ H^{1}\left(D\subset C,{\mathcal O}\right)/\operatorname{const} $ compatible with the natural Hilbert topology
on this space. \end{lemma}

\begin{proof} Obviously, the norm $ \|f\|_{1,\text{int}}\leq\|f\|_{H^{1}\left(D\subset C\right)/\operatorname{const}} $, thus it is
enough to show that one can majorate $ \|f\|_{H^{1}\left(D\subset C\right)/\operatorname{const}} $ given $ \|f\|_{1,\text{int}} $. In
other words, given a function $ f $ defined in $ D $ with $ \int_{D}|df|^{2}d\mu\leq1 $, it is
enough to construct a continuation $ g $ of this function to $ C $ so that
$ \int_{C}|dg|^{2}d\mu\leq M $ (for an appropriate $ M $ which does not depend on $ f $).

It is clear that the existence of such a continuation depends on the
local properties of $ f $ near $ \partial D $, thus one may assume that $ C $ is $ {\mathbb C}P^{1} $, $ D $ is
the unit disk. Now the proof can proceed as in the previous lemma. \end{proof}

\begin{lemma} \label{lm131.44}\myLabel{lm131.44}\relax  Consider a closed measure-0 subset $ \gamma $ of a complex curve $ C $.
Then $ {\mathcal H}_{E}^{1}\left(\gamma\subset C\right)={\mathcal H}^{1}\left(\gamma\subset C\right)=H^{1}\left(\gamma\subset C\right) $. \end{lemma}

\begin{proof} By definition of $ {\mathcal H}^{1} $, it is enough to show that $ H^{0}\left(\gamma\subset C\right)=0 $, or
that any $ L_{2} $-function on $ C $ can be approximated by an $ L_{2} $-function with (the
closure of) the support inside $ C\smallsetminus\gamma $. In turn, this follows from the fact
that there exists an open subset $ U $, $ \gamma\subset U\subset C $, with an arbitrary small
measure. \end{proof}

\begin{corollary} If $ \gamma $ is a smooth real curve inside a complex curve $ C $,
then $ {\mathcal H}^{1}\left(\gamma\subset C\right)\simeq H^{1/2}\left(\gamma\right) $. Suppose that $ \gamma $ breaks $ C $ into two pieces $ D_{\pm} $. Let an
excess space $ E $ be $ D_{\pm} $-split, and let $ H_{\pm}^{1/2}\left(\gamma\right) $ be images of $ {\mathcal H}_{E}^{1}\left(\bar{D}_{\pm}\right) $ inside
$ H^{1/2}\left(\gamma\right) $. Then $ H^{1/2}\left(\gamma\right)/\operatorname{const} $ is a direct sum of $ H_{+}^{1/2}\left(\gamma\right)/\operatorname{const} $ and
$ H_{-}^{1/2}\left(\gamma\right)/\operatorname{const} $ as a topological vector space. The mappings $ {\mathcal H}_{E}^{1}\left(\bar{D}_{\pm}\right)/\operatorname{const}
\to H_{\pm}^{1/2}\left(\gamma\right)/\operatorname{const} $ are isomorphisms of topological vector spaces. \end{corollary}

\begin{proof} Apply Theorem~\ref{th35.150} to $ D=\gamma $, and $ D_{1,2}=D_{\pm} $. Now the statement
follows from the lemmas above. \end{proof}

\begin{remark} Heuristically, the first part of the corollary is similar to
the following statements: any function from $ H^{1/2}\left(\gamma\right) $ can be approximated
by a function which is analytic near $ \gamma $; if $ C={\mathbb C}{\mathbb P}^{1} $, any analytic near $ \gamma $
function can be represented as a sum of two functions analytic in
neighborhoods of $ D_{+} $ and $ D_{-} $ correspondingly. In other words, it is similar
to a statement about density of analytic functions inside Sobolev spaces.
\end{remark}

\begin{remark} Since $ {\mathcal H}_{E}^{1}\left(\bar{D}_{\pm}\right)/\operatorname{const} $ carries a naturally defined Hilbert norm,
so do $ H_{\pm}^{1/2}\left(\gamma\right)/\operatorname{const} $, thus $ H^{1/2}\left(\gamma\right)/\operatorname{const} $. In this norm the subspaces
$ H_{\pm}^{1/2}\left(\gamma\right)/\operatorname{const} $ are orthogonal. Call this norm on $ H^{1/2}\left(\gamma\right)/\operatorname{const} $ the
{\em embedding norm}. In Definition~\ref{def131.85} we define a different norm on
$ H^{1/2}\left(\gamma\right)/\operatorname{const} $. \end{remark}

\begin{remark} The norm on $ H^{1/2}\left(\gamma\right)/\operatorname{const} $ and the decomposition
\begin{equation}
H^{1/2}\left(\gamma\right)/\operatorname{const}=H_{+}^{1/2}\left(\gamma\right)/\operatorname{const}\oplus H_{-}^{1/2}\left(\gamma\right)/\operatorname{const}
\notag\end{equation}
depend on the inclusion $ \gamma \hookrightarrow C $. Clearly, the subspaces $ H_{\pm}^{1/2}\left(\gamma\right)/\operatorname{const} $
depend only on inclusions $ \gamma \hookrightarrow \bar{D}_{\pm} $. However, the norms on $ H_{\pm}^{1/2}\left(\gamma\right)/\operatorname{const} $
depend also on the inclusions $ D_{\pm} \hookrightarrow C $, since the norms on $ {\mathcal H}^{1}\left(D_{\pm}\subset C\right)/\operatorname{const} $
depend on these inclusions. \end{remark}

From now on we assign indices $ + $ and $ - $ to the parts $ D_{\pm} $ into which an
oriented real curve $ \gamma\subset C $ breaks $ C $ so that the orientation of $ \gamma $ coincides
with the orientation of the boundary of $ D_{-} $. Typically, we will have
several clockwise circles $ \gamma_{j} $ bounding disjoint disks $ D_{j+} $; the complement
to these disks is $ D=\bigcap_{j}D_{j-} $; it is this complement we are interested this,
and $ \bigcup\gamma_{j} $ is the properly oriented boundary of $ D $.

\begin{remark} In the applications the parts $ D_{+} $ and $ D_{-} $ do not play symmetrical
roles. Typically, $ D_{+} $ is a ``small'' domain; moreover, if $ C={\mathbb C}{\mathbb P}^{1} $, $ D_{+} $ is often
a small disk.

Below we use a trick to postpone one complicated calculation until
Section~\ref{s3.55}. The trick boils down to changing the norm on the subspace
$ H_{+}^{1/2}\left(\gamma\right)/\operatorname{const} $ in the following way: \end{remark}

\begin{definition} \label{def131.85}\myLabel{def131.85}\relax  Consider a real oriented connected curve $ \gamma $ inside a
complex curve $ C $. Suppose that $ \gamma $ breaks $ C $ into two parts $ D_{\pm} $. Consider the
decomposition $ H^{1/2}\left(\gamma\right)/\operatorname{const} = H_{+}^{1/2}\left(\gamma\right)/\operatorname{const}\oplus H_{-}^{1/2}\left(\gamma\right)/\operatorname{const} $. Consider the
norm on $ H_{-}^{1/2}\left(\gamma\right)/\operatorname{const} $ defined by its identification with $ {\mathcal H}_{E}^{1}\left(C_{-}\subset C\right)/\operatorname{const} $
(with the natural norm), consider the norm on $ H_{+}^{1/2}\left(\gamma\right)/\operatorname{const} $ defined by
its identification with $ {\mathcal H}_{E}^{1}\left(C_{-}\subset C\right)/\operatorname{const} $ with the norm $ \|\bullet\|_{1,\text{int}} $ from
Definition~\ref{def131.35}. Call the induced direct sum norm on $ H^{1/2}\left(\gamma\right)/\operatorname{const} $
the $ + $-{\em skewed norm}. \end{definition}

The reason to consider the $ + $-skewed norm is Theorem~\ref{th35.157}.

\begin{remark} Section~\ref{s3.55} provides an alternative version of the theory
which does not use the $ + $-skewed norms. This allows dropping one of the
conditions on the pieces (quasi-circularity), the price being slightly
more complicated conditions on the ``distance'' between components of the
boundary of the pieces.

In other words, consideration of the $ + $-skewed norm provides some
shortcuts in the discussion which follows, but should not
significantly influence the class of curves allowed by these discussions. \end{remark}

\begin{remark} \label{rem131.90}\myLabel{rem131.90}\relax  Obviously, $ \|\alpha\|_{1,\text{int}}\leq\|\alpha\|_{1} $ for $ \alpha\in{\mathcal H}^{1} $. If $ D $ is a disk, it is
easy to show that $ \|\alpha\|_{1,\text{int}}=\|\alpha\|_{1}/\sqrt{2} $. Due to Lemma~\ref{lm131.40}, for a domain
$ D $ with a smooth boundary $ \gamma $, $ \|\alpha\|_{1}\leq c\|\alpha\|_{1,\text{int}} $ for an appropriate number $ c $.
Call the minimal such number the {\em distortion\/} $ \Delta\left(\gamma\right) $ {\em of the curve\/} $ \gamma $. \end{remark}

\subsection{Restriction to boundary curves } Consider a complex curve $ C $ with a
closed subset $ D $. Let $ R_{j} $, $ j\in J $, be the collection of connected components
of $ C\smallsetminus D $. Let $ \mathring{D}_{j} $ be the interior of $ C\smallsetminus R_{j} $. An excess space $ E $ is $ D $-{\em adjusted\/}
if it is $ R_{\bullet} $-split, and is $ \left\{\mathring{D}_{j},R_{j}\right\} $-split for any $ j\in J $. Starting from this
section, we consider only $ D $-adjusted excess spaces.

\begin{definition} Consider a complex curve $ C $ with a closed subset $ D $. Suppose
that one of the connected components $ R_{j}\subset C $ of $ C\smallsetminus D $ is bounded by an
oriented smooth curve $ \gamma_{j} $. Denote the well-defined restriction mappings
$ H^{1}\left(D\subset C\right)/\operatorname{const} \to H^{1/2}\left(\gamma_{j}\right)/\operatorname{const} $ by $ \beta_{j} $, the compositions of $ \beta_{j} $ with the
projections $ H^{1/2}\left(\gamma_{j}\right)/\operatorname{const} \to H_{\pm}^{1/2}\left(\gamma_{j}\right)/\operatorname{const} $ by $ \beta_{j\pm} $. \end{definition}

\begin{definition} Call a closed subset $ R $ of a compact complex curve $ C $
{\em pseudo-smooth\/} if all the connected components $ D_{j} $, $ j\in J $, of $ C\smallsetminus R $ have smooth
boundaries $ \gamma_{j} $. Call curves $ \gamma_{j} $ the {\em smooth boundaries\/} of $ R $. \end{definition}

\begin{definition} \label{def35.153}\myLabel{def35.153}\relax  Suppose that $ D $ is pseudo-smooth with smooth
boundaries $ \gamma_{j} $, $ j\in J $. Denote by $ \widetilde{\beta} $ the mapping $ \prod_{j}\beta_{j}\colon H^{1}\left(D\subset C\right)/\operatorname{const} \to
\prod_{j}H^{1/2}\left(\gamma_{j}\right)/\operatorname{const} $, and by $ \widetilde{\beta}_{\pm} $ the mappings $ \prod_{j}\beta_{j\pm}\colon H^{1}\left(D\subset C\right)/\operatorname{const} \to
\prod_{j}H_{\pm}^{1/2}\left(\gamma_{j}\right)/\operatorname{const} $. \end{definition}

Since $ H^{1/2}\left(\gamma_{j}\right)/\operatorname{const} $ is equipped with a natural Hilbert norm, it
makes sense to consider $ \bigoplus_{l_{2}}H^{1/2}\left(\gamma_{j}\right)/\operatorname{const} $. (This is the first place
where the distinction between Hilbert norms and Hilbert topologies
becomes important.)

\begin{theorem} \label{th35.155}\myLabel{th35.155}\relax  The mapping $ \widetilde{\beta}_{-} $ defined above sends $ {\mathcal H}_{E}^{1}\left(D\subset C\right)/\operatorname{const} $ into
$ V_{-}=\bigoplus_{l_{2}}H_{-}^{1/2}\left(\gamma_{j}\right)/\operatorname{const} $. Then the induced mapping $ \beta_{-} $ into $ V_{-} $ is an
invertible continuous mapping. Moreover, $ \beta_{-} $ is unitary if $ g\left(C\right)=0 $. \end{theorem}

\begin{proof} Indeed, $ \beta $ is a composition of the mapping $ \rho^{-1} $ of Theorem~%
\ref{th35.150} with a direct sum of the mapping $ {\mathcal H}^{1}\left(D_{j}\subset C\right)/\operatorname{const} \to
H_{-}^{1/2}\left(\gamma_{j}\right)/\operatorname{const} $; here $ D_{j}=C\smallsetminus R_{j} $. By definition of $ H_{-}^{1/2} $, the latter mapping
is an isomorphism of Hilbert spaces. \end{proof}

\begin{remark} \label{rem35.156}\myLabel{rem35.156}\relax  What is important to us in this result is that though we
assume that the boundary of each of $ R_{j} $ is smooth, we do not assume that
the {\em whole\/} boundary of $ D $ is smooth. Indeed, if the number $ |J| $ of connected
components of $ C\smallsetminus D $ is infinite, then in addition to $ \bigcup_{j\in J}\partial R_{j} $, $ \partial D $ contains
also the {\em dust\/}: all the accumulation points of the curves $ \partial R_{j} $. It is easy
to construct examples when the dust is very massive; see Section~\ref{s6.3}.
For example, it may have a positive measure. It may also coincide with $ D $.
\end{remark}

\begin{theorem} \label{th35.157}\myLabel{th35.157}\relax  The mapping $ \widetilde{\beta} $ defined above sends $ {\mathcal H}_{E}^{1}\left(D\subset C\right) $ into
$ \bigoplus_{l_{2}}H^{1/2}\left(\gamma_{j}\right)/\operatorname{const} $, here each $ H^{1/2}\left(\gamma_{j}\right)/\operatorname{const} $ is equipped with the
$ + $-skewed norm. \end{theorem}

\begin{proof} It is enough to show that $ \operatorname{Im}\widetilde{\beta}_{+}\subset\bigoplus_{l_{2}}H_{+}^{1/2}\left(\gamma_{j}\right)/\operatorname{const} $.
In other words, given an $ H^{1} $-function $ f $ on $ C $, the
sequence $ \left(n_{j}\right) $ is in $ l_{2} $, here $ n_{j}\buildrel{\text{def}}\over{=}\|\left(f|_{\gamma_{j}}\right)_{+}\|_{H_{+}^{1/2}\left(\gamma_{j}\right)/\operatorname{const}} $, and $ g_{\pm} $ are
the $ \pm $-components of $ g\in H^{1/2}\left(\gamma_{j}\right) $. In turn, by Lemma~\ref{lm131.40} this
follows from the restriction mapping $ H^{1}\left(C\right)/\operatorname{const} \to \bigoplus_{l_{2}}H^{1}\left(D_{j}\right)/\operatorname{const} $
having a norm $ \leq1 $ if we consider norms $ \|\|_{1,\text{int}} $ on $ H^{1}\left(D_{j}\right)/\operatorname{const} $. In turn,
the latter statement follows from the definition of the norm $ \|\|_{1,\text{int}} $. \end{proof}

\begin{definition} Let $ \beta_{E} $ be the mapping $ {\mathcal H}_{E}^{1}\left(D\subset C\right) \to \bigoplus_{l_{2}}H^{1/2}\left(\gamma_{j}\right)/\operatorname{const} $ induced
by $ \widetilde{\beta} $, let $ \beta=\beta_{E}|_{{\mathcal H}^{1}\left(D\subset C\right)} $. \end{definition}

\section{Gluing the curve from the pieces }

\subsection{Gluing data and mismatch }\label{s30.10}\myLabel{s30.10}\relax 

\begin{definition} Consider two connected oriented closed real curves $ \gamma_{1} $ and
$ \gamma_{2} $. The {\em curve gluing data\/} for the pair $ \left(\gamma_{1},\gamma_{2}\right) $ is a pair of mutually
inverse diffeomorphisms $ \varphi_{1}\colon \gamma_{1} \to \gamma_{2} $ and $ \varphi_{2}\colon \gamma_{2} \to \gamma_{1} $ which reverse the
orientations. The {\em bundle gluing data\/} for the pair $ \left(\gamma_{1,2}\right) $ and the curve
gluing data $ \left(\varphi_{1,2}\right) $ is a pair of smooth complex-valued functions $ \psi_{i} $ on
$ \gamma_{i} $, $ i=1,2 $, such that $ \psi_{1}\cdot\varphi_{1}^{*}\left(\psi_{2}\right)=1 $.

The {\em degree\/} of the bundle gluing data is\footnote{Recall that $ \operatorname{ind}\psi\buildrel{\text{def}}\over{=}\frac{1}{2\pi i}\int d \log \psi $; here $ \psi $ is a nowhere-0 function on
$ S^{1} $.} $ \operatorname{ind}\psi_{1}=\operatorname{ind}\psi_{2} $. \end{definition}

Clearly, the gluing data for a pair $ \left(\gamma_{1},\gamma_{2}\right) $ induces gluing data for
a pair $ \left(\gamma_{2},\gamma_{1}\right) $.

\begin{definition} Given the curve and bundle gluing data $ \left(\varphi_{1,2},\psi_{1,2}\right) $ for a pair
$ \left(\gamma_{1,2}\right) $, call a pair of functions $ f_{j}\in H^{s}\left(\gamma_{j}\right) $, $ j=1,2 $, {\em compatible with the
gluing data\/} if $ f_{1}=\psi_{1}\cdot\varphi_{1}^{*}\left(f_{2}\right) $. The {\em mismatch\/} of a pair $ \left(f_{1,2}\right) $ is a pair
$ \left(\delta_{1,2}\right) $, $ \delta_{j}\in H^{s}\left(\gamma_{j}\right) $, $ j=1,2 $, given by $ \delta_{j}=f_{j}-\psi_{j}\cdot\varphi_{j}^{*}\left(f_{k}\right) $, for $ \left(j,k\right)=\left(1,2\right) $ or
(2,1). \end{definition}

Clearly, if $ \left(f_{1},f_{2}\right) $ is compatible with gluing data for $ \left(\gamma_{1},\gamma_{2}\right) $, then
$ \left(f_{2},f_{1}\right) $ is compatible with the corresponding gluing data for $ \left(\gamma_{2},\gamma_{1}\right) $.
Similarly, if $ \left(\delta_{1},\delta_{2}\right) $ is the mismatch of $ \left(f_{1},f_{2}\right) $, then $ \left(\delta_{2},\delta_{1}\right) $ is the
mismatch of $ \left(f_{2},f_{1}\right) $. Moreover, $ \delta_{2}=-\psi_{2}\cdot\varphi_{2}^{*}\left(\delta_{1}\right) $.

\begin{definition} Consider a collection $ {\mathcal D}=\left(D_{k}\subset C_{k}\right)_{k\in K} $ of pseudo-smooth closed
subsets of compact complex curves. Let $ J_{k} $ be the set of connected
components of $ C_{k}\smallsetminus D_{k} $, $ J=\coprod_{k}J_{k} $. Let $ \gamma_{j} $, $ j\in J $, be the boundary of the
connected component which corresponds to $ j $. The {\em gluing decomposition\/} for
$ {\mathcal D} $ is a decomposition of $ J $ into a disjoint union of pairs $ \left\{j,j'\right\} $. Define a
mapping $ '\colon J \to J\colon j \mapsto j'\colon j' \mapsto j $. The {\em curve and bundle gluing data\/} for
such a gluing is a gluing decomposition together with curve and bundle
gluing data $ \left(\varphi_{j},\varphi_{j'}\right) $, $ \left(\psi_{j},\psi_{j'}\right) $ for the pair of curves $ \left(\gamma_{j},\gamma_{j'}\right) $ for each
pair $ \left\{j,j'\right\} $.

Suppose that $ \operatorname{ind}\psi_{j}=0 $ for all but a finite number of $ j\in J $. Then the
{\em degree\/} of the bundle gluing data is the sum of degrees over all pairs
$ \left(\gamma_{j},\gamma_{j'}\right) $.

Call a collection of functions $ F_{k}\in H^{1}\left(D_{k}\subset C_{k}\right) $, $ k\in K $, {\em compatible with
the gluing data}, if the pair $ \left(f_{j},f_{j'}\right) $ is compatible with the gluing data
for $ \left(\gamma_{j},\gamma_{j'}\right) $ for each pair $ \left\{j,j'\right\} $; here $ f_{j}=F_{k_{j}}|_{\gamma_{j}} $ (assuming that $ \gamma_{j}\subset C_{k_{j}} $).
Define similarly the {\em mismatch\/} $ \left(\delta_{j}\right) $, $ \delta_{j}\in H^{1/2}\left(\gamma_{j}\right) $, $ j\in J $, of such a
collection. \end{definition}

It is clear that mismatches $ \delta_{j} $, $ \delta_{j'} $ satisfy $ \delta_{j'}=-\psi_{j'}\cdot\varphi_{j'}^{*}\left(\delta_{j}\right) $.

\subsection{The curve and the line bundle }\label{s3.2}\myLabel{s3.2}\relax  Given the curve and bundle gluing
data $ \left(D_{k}\subset C_{k},',\varphi_{\bullet},\psi_{\bullet}\right) $, one can associate to it some more or less familiar
objects. The {\em associated curve\/} $ C $ is the set obtained from $ \coprod_{k}D_{k} $ by
identifying points on the smooth parts of the boundaries via $ \varphi_{\bullet} $. As
quotients do, $ C $ is equipped with a natural topology. A point $ m\in C $ has one
or two preimages. A point $ m\in C $ is a {\em dust point\/} if one of its preimages on
$ \coprod_{k}D_{k} $ lies in the dust of the corresponding component $ D_{k} $. The {\em dust\/} $ C_{\infty}\subset C $
consists of dust points; it is a closed subset of $ C $. Obviously, $ C_{\infty} $ is
empty unless a complement to one of $ D_{k} $ has infinity many connected
components. It is clear that $ C_{\text{fin}}\buildrel{\text{def}}\over{=}C\smallsetminus C_{\infty} $ has a natural structure of a
complex curve. However, $ C_{\text{fin}} $ may be empty.

If $ \gamma_{j} $ and $ \gamma_{j'} $ have no dust points on them, the common image of these
curves on $ C_{\text{fin}} $ is a smooth cycle on $ C_{\text{fin}} $.

Similarly to $ C $, one can glue a ``set-theoretic\footnote{I.e., a set with a projection $ \pi\colon {\mathcal L} \to C $; fibers of $ \pi $ are one-dimensional
vector spaces.} line bundle'' $ {\mathcal L} $ over $ C $
starting from $ \coprod_{k}D_{k}\times{\mathbb C} $ and gluing via $ \left(\varphi_{\bullet},\psi_{\bullet}\right) $. It may be not a topological
line bundle; however, it is an analytic line bundle over $ C_{\text{fin}} $. Usual
definitions of the dual bundle, of the tensor product of bundles, of the
line bundles $ \omega $ and $ \bar{\omega} $ work without any change in this situation.

Let $ H_{\text{loc}}^{1}\left(C,{\mathcal L}\right) $ consist of $ H_{\text{loc}}^{1} $-{\em sections\/} of $ {\mathcal L} $ on $ C $: an $ H_{\text{loc}}^{1} $-{\em section\/}
is a collection of elements of $ H^{1}\left(D_{k}\subset C_{k}\right) $ which are compatible with the
gluing data. (One can naturally define the {\em support\/} of an $ H_{\text{loc}}^{1} $-section,
thus one can also define what is a {\em section\/} of $ {\mathcal L} $ on $ U\subset C $.) Similarly one
can define $ H_{\text{loc}}^{0}\left(C,{\mathcal L}\right) $ (without any compatibility conditions on $ \gamma_{j} $) and
the operator $ \bar{\partial}\colon H_{\text{loc}}^{1}\left(C,{\mathcal L}\right) \to H_{\text{loc}}^{0}\left(C,{\mathcal L}\otimes\bar{\omega}\right) $. The vector space $ {\mathcal H}_{\text{loc}}^{1}\left(C,{\mathcal L}\right) $
of {\em locally-\/}$ H^{1} $-{\em holomorphic\/} sections consists of $ f\in H_{\text{loc}}^{1}\left(C,{\mathcal L}\right) $ such that
$ \bar{\partial}f=0 $; here $ \bar{\partial}f\in H_{\text{loc}}^{0}\left(C,{\mathcal L}\otimes\bar{\omega}\right) $. Define similarly the space $ {\mathcal H}_{E,\text{loc}}^{1}\left(C,{\mathcal L}\right) $.

The motivation for this definition is the following:

\begin{lemma} \label{lm351.20}\myLabel{lm351.20}\relax  Consider a complex curve $ C $ and a compact smooth real
curve $ \gamma\subset C $ which breaks $ C $ into two parts $ D_{\pm} $. Consider two functions
$ f_{\pm}\in H_{\text{loc}}^{1}\left(D_{\pm}\subset C\right) $ such that $ f_{+}|_{\gamma}=f_{-}|_{\gamma} $ and $ \bar{\partial}f_{\pm}=0\in H_{\text{loc}}^{0}\left(D_{\pm}\subset C\right) $. Then there is a
unique function $ f\in H_{\text{loc}}^{1}\left(C\right) $ such that $ f_{\pm}=f|_{D_{\pm}} $. Moreover, $ \bar{\partial}f=0\in H_{\text{loc}}^{0}\left(C\right) $. In
particular, $ f $ is analytic near $ \gamma $. \end{lemma}

\begin{proof} If $ f $ exists, then obviously $ \bar{\partial}f\in L_{2,\text{loc}} $ should vanish. The
uniqueness of $ f $ is also obvious. Show the existence of $ f $; we may
drop the conditions $ \bar{\partial}f_{\pm}=0 $.

We know already that the mapping of restriction to $ \gamma $ is surjective,
thus we may suppose $ f_{+}|_{\gamma}=f_{-}|_{\gamma}=0 $. It is enough to show that $ f_{+} $
allows a continuation-by-0 without loosing its
smoothness class $ H^{1} $. This is a local statement, so we may assume that $ C $
is a neighborhood of the unit circle $ \gamma $ in $ {\mathbb C} $.

Again, consider the action of the group of rotations. One may
restrict the attention to one isotypical component in $ H^{1}\left(C\right) $. This reduces
the problem to one-dimensional: given an $ H^{1} $-function $ g\left(t\right) $ on $ \left[0.5,2\right] $ such
that $ g\left(1\right)=0 $, the extension-by-0 from $ \left[0.5,1\right] $ to $ \left[0.5,2\right] $ obviously has the
same norm. \end{proof}

In this paper we are most interested in line bundles $ {\mathcal L} $ of small
degree. As explained in Section~\ref{s01.60}, this leads to consideration of
the operator $ \bar{\partial} $ sending $ H^{1} $ to $ H^{0} $. This is the reason for our interest in
$ H^{1} $-holomorphic functions.

\begin{remark} Above, the index loc relates to having no restriction on how
the sequence $ \|f_{k}\| $ grows; here $ f_{k} $ is the restriction of $ f\in H^{s}\left(C,{\mathcal L}\right) $ to $ D_{k}\subset C $.
Such an approach is sufficient if $ K $ is finite (which is the most
interesting case in our approach). However, it is easy to modify this to
work with infinite collections $ K $, see Section~\ref{s10.80}. \end{remark}

Recall (see Section~\ref{s01.70}) that to expect Riemann--Roch theorem to
hold, one needs to add some slack, allowing some non-strictly holomorphic
sections. By robustness, it is not very important which non-holomorphic
functions are allowed; we add slack by allowing a finite-dimensional
mismatch at each gluing:

\begin{definition} \label{def351.15}\myLabel{def351.15}\relax  Given a curve and bundle gluing data $ \left(\varphi_{j},\psi_{j}\right)_{j\in J} $ for
a collection $ \left(D_{k}\subset C_{k}\right)_{k\in K} $ with boundary curves $ \left(\gamma_{j}\right)_{j\in J} $, the {\em mismatch
allowance\/} is a collection $ \left(V_{j}\right)_{j\in J} $ consisting of vector subspaces
$ V_{j}\subset H^{1/2}\left(\gamma_{j}\right) $ such that $ V_{j'}=\psi_{j'}\cdot\varphi_{j'}^{*}\left(V_{j}\right) $.

An $ H_{\text{loc}}^{1} $-{\em section\/} $ F $ {\em modulo\/} $ \left(V_{j}\right)_{j\in J} $ is a collection $ \left(F_{k}\right)_{k\in K} $ such that
$ F_{k}\in H^{1}\left(D_{k}\subset C_{k}\right) $ and the mismatch $ \left(\delta_{j}\right)_{\psi\in J} $ of $ \left(F_{k}\right) $ satisfies $ \delta_{j}\in V_{j} $.
Define similarly $ {\mathcal H}_{\text{loc}}^{1} $-{\em sections\/} and $ {\mathcal H}_{E,\text{loc}}^{1} $-{\em sections modulo\/} $ \left(V_{j}\right) $. \end{definition}

\begin{definition} Consider an involution ' of a set $ J $. Given a quantity $ t_{j} $,
$ j\in J $, such that $ t_{j}=t_{j'} $, let $ \sum_{\left\{j,j'\right\}}t_{j}=\frac{1}{2}\sum_{j\in J}t_{j} $. Similarly, if $ V_{j} $ is a
vector space with a fixed isomorphism $ \varphi_{i} $ between $ V_{j} $ and $ V_{j'} $ such that
$ \varphi_{i'}=\varphi_{i}^{-1} $, let $ \bigoplus_{\left\{j,j'\right\}}V_{j} $ is the subspace of $ \bigoplus_{j}V_{j} $ formed by sequences
$ \left(v_{j}\right)_{j\in J} $ such that $ v_{j'}=\varphi_{j}v_{j} $. \end{definition}

If $ J=J_{0}\coprod J_{0}' $, then $ \bigoplus_{\left\{j,j'\right\}}V_{j} $ is canonically isomorphic to $ \bigoplus_{j\in J_{0}}V_{j} $.
Similarly, define $ \prod_{\left\{j,j'\right\}} $ etc.

\subsection{Finite-genus Riemann--Roch theorem via gluing data }\label{s10.70}\myLabel{s10.70}\relax  This theorem
relates the dimension of two vector spaces. One is the space of global
sections of a line bundle on a curve. Another is the first homology. The
first step to formulate the infinite-genus variant is the
translation of the case $ g<\infty $ to our notations.

\begin{theorem} In the conditions of the previous section suppose that the
set $ K $ is finite, and each subset $ D_{k} $ has a smooth boundary (thus its
complement $ C_{k}\smallsetminus D_{k} $ has finitely many connected components). Then the vector
subspace of $ \prod_{\left\{j,j'\right\}}H^{1/2}\left(\gamma_{j}\right)=\bigoplus_{\left\{j,j'\right\}}H^{1/2}\left(\gamma_{j}\right) $ formed by mismatches of
elements of $ \prod_{k}{\mathcal H}^{1}\left(D_{k}\subset C_{k}\right) $ is a closed subspace of finite codimension.
Denote this codimension by $ h^{1} $.

The vector space of global analytic sections is finite dimensional,
denote its dimension by $ h^{0} $. Then $ h^{0}-h^{1}=d-g+1 $, here $ d $ is the degree of the
bundle gluing data, and $ g=|J|/2-|K|+1+\sum_{k}g\left(C_{k}\right) $. \end{theorem}

\begin{proof}[Sketch of the proof ] Since in what follows we are going to prove
significantly more general results, let us show only that this statement
is a generalization of the ``usual'' Riemann--Roch theorem, which is
formulated using the language of analytic sections, not Sobolev-class
section. To simplify the discussion, assume that the gluing data $ \left(\varphi_{j}\right) $ and
$ \left(\psi_{j}\right) $ consists of real-analytic functions, that no connected components of
the boundary of the same piece $ D_{k} $ are glued together, and that the
complex analytic curve $ C $ obtained after gluing is connected.

Since functions $ \varphi_{j} $ can be analytically extended to neighborhoods of
$ \gamma_{j} $, we can glue $ C $ of neighborhoods $ \widetilde{D}_{k} $ of $ D_{k} $. Images $ U_{k} $ of $ \widetilde{D}_{k} $ in $ C $ form a
covering of $ C $, $ U_{k}\cap U_{k'} $ is a union of annuli (thus Stein), and $ U_{k} $ is Stein
(unless $ |K|=1 $, $ |J|=0 $, when the theorem is obvious).

Moreover, $ D_{k} $ is identified with a closed subset of $ C $, so we may
assume that $ D_{k}\subset C $ and $ \gamma_{j}\subset C $, and $ \gamma_{j}=\gamma_{j'} $ up to orientation change. We may
assume that the bundle gluing functions $ \psi_{j} $ on $ \gamma_{j} $ can be extended to the
corresponding connected component $ V_{\left\{j,j'\right\}} $ of $ U_{k}\cap U_{k'} $, thus define a line
bundle $ {\mathcal L} $ over $ C $; sections of $ {\mathcal L} $ are represented by functions on $ U_{k} $ with
vanishing
mismatches of boundary values. A simple calculation shows that $ \deg {\mathcal L}=d $,
$ g\left(C\right)=g $. Everything being Stein, we can calculate cohomology by the \v Cech
complex $ {\mathcal C}_{\text{an}} $
\begin{equation}
0 \to \bigoplus_{k\in K}\Gamma_{\text{an}}\left(U_{k},{\mathcal L}\right) \to \bigoplus_{\left\{j,j'\right\}}\Gamma_{\text{an}}\left(V_{j,j'},{\mathcal L}\right) \to\text{ 0.}
\notag\end{equation}
Moreover, $ {\mathcal L}|_{U_{k}} $ is already trivialized, so we may substitute $ {\mathcal O} $ instead of
$ {\mathcal L} $, with an appropriate modification of the differential of the complex.
After this change the differential becomes the operator of taking
the mismatch. We want to show that the cohomology of the complex above
coincides with the cohomology of the complex $ {\mathcal C}_{H} $
\begin{equation}
0 \to \bigoplus_{k\in K}{\mathcal H}^{1}\left(D_{k}\right) \xrightarrow[]{\widetilde{\mu}} \bigoplus_{\left\{j,j'\right\}}H^{1/2}\left(\gamma_{j}\right) \to\text{ 0.}
\notag\end{equation}
Call the differential $ \widetilde{\mu} $ the {\em operator of taking the mismatch}.

There is natural inclusion $ {\mathcal C}_{\text{an}} \hookrightarrow {\mathcal C}_{H} $, consider the induced mapping
of cohomology. On the level of $ {\mathbit H}^{0} $ (here $ {\mathbit H} $ denotes cohomology) it is
automatically an injection. On the other hand, as Lemma~\ref{lm351.20} shows,
any function compatible with the gluing data $ \left(\varphi_{\bullet},\psi_{\bullet}\right) $ induces an
$ {\mathcal H}^{1} $-section of $ {\mathcal L} $. Thus the mapping of $ {\mathbit H}^{0} $ is an isomorphism.

Consider the spaces $ {\mathbit H}^{1} $. Due to the duality theorem,
$ {\mathbit H}^{1}\left(C,{\mathcal L}\right)^{*}={\mathbit H}^{0}\left(C,\omega\otimes{\mathcal L}^{-1}\right) $. Given a section $ \alpha $ of $ \omega\otimes{\mathcal L}^{-1} $, the pairing with a
$ 1 $-cocycle $ c\in\bigoplus_{\left\{j,j'\right\}}\Gamma_{\text{an}}\left(V_{j,j'},{\mathcal L}\right) $ is $ \int_{\gamma}\alpha c $, here $ \gamma $ is a suitably oriented
curve with connected components generating $ 1 $-homology of annuli $ V_{j,j'} $.
Taking $ \gamma=\bigcup_{\left\{j,j'\right\}}\gamma_{j} $ shows that the linear functional on $ {\mathbit H}^{1}\left({\mathcal C}_{\text{an}}\right) $ induced
by $ \alpha $ can be passed through the mapping $ {\mathcal C}_{\text{an}} \to {\mathcal C}_{H} $, thus $ {\mathbit H}^{1}\left({\mathcal C}_{\text{an}}\right) \to {\mathbit H}^{1}\left({\mathcal C}_{H}\right) $
is an injection. Since $ \operatorname{Im}{\mathcal C}_{\text{an}} $ is dense in $ {\mathcal C}_{H} $, to show the surjectivity it
is enough to show that the image of each component $ {\mathcal H}^{1}\left(D_{k}\right) \to \bigoplus H^{1/2}\left(\gamma_{j}\right) $
of the differential is closed, which follows from Theorem~\ref{th35.155}. \end{proof}

\subsection{Plan of the campaign }\label{s10.80}\myLabel{s10.80}\relax  We have shown that the finite-genus case
of Riemann--Roch theorem coincides with the calculation of the index of
the operator $ \widetilde{\mu} $ of taking the mismatch. The target of this paper is to
investigate the mismatch operator in the more general case of arbitrary
(possibly infinite) genus---assuming that the degree remains finite.
There are two obstacles to restate the above theorem in the infinite
genus case: first, $ g-d $ becomes infinite. Second, by Theorem~\ref{th35.155}, the
natural topology on $ \operatorname{Im}\widetilde{\mu} $ is the Hilbert space topology, which is very far
from both the topology on $ \prod_{\left\{j,j'\right\}}H^{1/2}\left(\gamma_{j}\right) $ and on $ \bigoplus_{\left\{j,j'\right\}}H^{1/2}\left(\gamma_{j}\right) $, so
there is no hope to get a finite-dimensional cokernel of $ \widetilde{\mu} $.

The trick to tackle the first problem is the allowance subspaces we
introduced above. Instead of considering the mapping $ \widetilde{\mu} $ to
$ \bigoplus_{\left\{j,j'\right\}}H^{1/2}\left(\gamma_{j}\right) $, consider the induced mappings into
$ \bigoplus_{\left\{j,j'\right\}}\left(H^{1/2}\left(\gamma_{j}\right)/V_{j}\right) $; here $ V_{j} $ is an arbitrary finite-dimensional
subspace of $ H^{1/2}\left(\gamma_{j}\right) $. This would change the right-hand side of
Riemann--Roch theorem to $ d-g+1+\sum_{\left\{j,j'\right\}}\dim  V_{j} $. Now if $ \dim  V_{j} $ ``compensates''
the contribution of $ j\in J $ into $ g $, then $ d-g+1+\sum_{\left\{j,j'\right\}}\dim  V_{j} $ makes sense as a
finite number. This is so, for example, if $ |K|<\infty $, and $ \dim  V_{j}=1 $ for all
but a finite number of $ j\in J $.

Similarly, to compensate for an infinite $ K $ with $ g\left(C_{k}\right)=0 $ for almost
any $ k\in K $, it is enough to consider $ \bigoplus_{k\in K}{\mathcal H}^{1}\left(D_{k}\right)/W_{k} $ instead of $ \bigoplus_{k\in K}{\mathcal H}^{1}\left(D_{k}\right) $;
here $ W_{k} $ is an arbitrary $ 1 $-dimensional subspace of $ {\mathcal H}^{1}\left(D_{k}\right) $ (it is
convenient to assume that restrictions of functions from $ W_{k} $ to any
boundary component $ \gamma_{j} $ of $ C_{k} $ are in $ V_{j} $). As we will see, it is most
convenient to take $ W_{k} $ consisting of constant functions.

Finally, one can also put correction terms if infinitely many curves
$ C_{k} $ are not rational. Suppose that $ {\mathcal H}^{1}\left(D_{k}\right)\subset\widetilde{{\mathcal H}}^{1}\left(D_{k}\right)\subset H^{1}\left(D_{k}\subset C_{k}\right) $, and $ \dim 
\widetilde{{\mathcal H}}^{1}\left(D_{k}\right)/{\mathcal H}^{1}\left(D_{k}\right)=g\left(C_{k}\right) $. Then substitution of $ \widetilde{{\mathcal H}} $ instead of $ {\mathcal H} $ leads to the
index formula for $ \widetilde{\mu} $ which does not include $ g\left(C_{k}\right) $. Moreover, one can take
$ \widetilde{{\mathcal H}}^{1}={\mathcal H}_{E}^{1} $ for an appropriate $ D $-adjusted and $ D $-supported excess space $ E $.

As a result, we obtain the following reformulation of the
finite-genus Riemann--Roch theorem:

\begin{theorem} In the conditions of the previous section suppose that the
set $ K $ is finite, and each subset $ D_{k} $ has a smooth boundary (thus the
complement $ C\smallsetminus D_{k} $ has finitely many connected components). Suppose also
that $ 1\in V_{j} $ for any $ j $. Then the operator of taking the mismatch
modulo $ V_{j} $
\begin{equation}
\bigoplus_{k\in K}{\mathcal H}_{E}^{1}\left(D_{k}\right)/\operatorname{const} \xrightarrow[]{\mu'} \bigoplus_{\left\{j,j'\right\}}H^{1/2}\left(\gamma_{j}\right)/V_{j}
\notag\end{equation}
is Fredholm, and its index is equal to $ d+\sum_{\left\{j,j'\right\}}\left(\dim  V_{j}-1\right) $. \end{theorem}

The way to handle the second problem is suggested by Theorem~%
\ref{th35.155}: the space which appears in the theory in a natural way is
$ \bigoplus_{l_{2}}\left(H^{1/2}\left(\gamma_{j}\right)/\operatorname{const}\right) $; thus it is natural to replace $ \bigoplus_{\left\{j,j'\right\}}H^{1/2}\left(\gamma_{j}\right)/V_{j} $
by $ \bigoplus_{l_{2},\left\{j,j'\right\}}H^{1/2}\left(\gamma_{j}\right)/V_{j} $, provided $ V_{j}=\left<\operatorname{const}\right> $ for all but a finite
number of indices. Similarly, one should replace $ \bigoplus $ by $ \bigoplus_{l_{2}} $ as the target
of $ \mu' $ as well.

Combining these two arguments, we need to investigate the operator
\begin{equation}
\bigoplus_{l_{2},k\in K}{\mathcal H}_{E_{k}}^{1}\left(D_{k}\right)/\operatorname{const} \xrightarrow[]{\mu} \bigoplus_{l_{2},\left\{j,j'\right\}}H^{1/2}\left(\gamma_{j}\right)/V_{j}
\notag\end{equation}
of taking the mismatch; here we suppose that $ 1\in V_{j} $, $ E_{k} $ is $ D_{k} $-supported,
and $ \dim  V_{j}=1 $ for all but a finite number of $ j\in J $. (Note that these
conditions imply that $ \psi_{j}=\operatorname{const} $ for all but a finite number of $ j\in J $.) We
want to find the cases when this operator is continuous, Fredholm, and of
the index prescribed by Riemann--Roch theorem.

One part of this question is trivial to answer: by Theorem~%
\ref{th35.155}, $ \mu $ is continuous as far as we consider the $ + $-skewed norms on
$ H^{1/2}\left(\gamma_{j}\right)/V_{j} $. (Later we will see that this skewing may be replaces by
appropriate assumptions about curves $ \gamma_{j} $, but for the time being restrict
our attention to the $ + $-skewed norm.) The other questions require more
information about properties of the operator $ \mu $.

\section{Riemann--Roch theorem }\label{h42}\myLabel{h42}\relax 

\subsection{The mapping of gluing }\label{s4.11}\myLabel{s4.11}\relax  Here we show that the gluing condition on one
particular pair of curves $ \gamma_{j} $, $ \gamma_{j'} $ may be reformulated in terms of an
operators sending the $ + $-parts of a pair of functions to the $ - $-part.

\begin{definition} Consider a complex curve $ C $ with an excess space $ E $ and a real
smooth oriented curve $ \gamma\subset C $ splitting $ C $ into two domains $ D_{\pm} $ (as usual, the
orientation of $ \gamma $ is compatible with the orientation of $ D_{-} $). If $ D_{+} $ has no
handles, and $ E $ is $ D_{-} $-supported, call the inclusion $ \gamma\hookrightarrow C $ {\em admissible}.

Given two such pairs $ \gamma\subset C $, $ \gamma'\subset C' $ with excess spaces $ E $, $ E' $ which are
$ D_{-}- $ and $ D_{-}' $-supported, say that an orientation-inverting gluing $ \varphi\colon \gamma \to \gamma' $
is {\em compatible\/} with $ E $, $ E' $ if $ E\oplus E' $ is an excess space for the complex curve
$ C^{*}=D_{-}\cup_{\varphi}D'_{-} $. \end{definition}

Consider bundle gluing data $ \psi $, $ \psi' $ for the identification $ \varphi $. It
determines a line bundle $ {\mathcal L} $ over $ C^{*} $. Construct 3 subspaces $ \left(V_{+},V_{-},L\right) $ of
the vector space $ V=H^{1/2}\left(\gamma\right)\oplus H^{1/2}\left(\gamma'\right)\colon V_{\pm} $ correspond to $ \pm $-parts of $ H^{1/2} $,
and $ L $ consists of pairs $ \left(f,f'\right)\in V $ which are compatible with the curve and
bundle gluing data: $ f=\psi\cdot\varphi^{*}\left(f'\right) $.

The following lemma technical lemma is a key tool for translating
properties of the operator $ \mu $ of mismatch (from Section~\ref{s10.80}) to the
usual Fredholm theorem. Heuristically, it states that these 3 subspaces
are in general position, and calculates the corresponding relative
dimensions:

\begin{lemma} \label{lm30.70}\myLabel{lm30.70}\relax  Let $ \bar{V}=H^{1/2}\left(\gamma\right)/\operatorname{const}\oplus H^{1/2}\left(\gamma'\right)/\operatorname{const} $, $ \bar{L} $ be the image of $ L $
w.r.t.~the natural projection $ \pi\colon V \to \bar{V} $. Let $ \bar{V}=\bar{V}_{+}\oplus\bar{V}_{-} $ be the decomposition
of $ \bar{V} $ corresponding to the decompositions $ H^{1/2}/\operatorname{const} =H_{+}^{1/2}/\operatorname{const}\oplus
H_{-}^{1/2}/\operatorname{const} $. Then there is a continuous mapping $ A\colon \bar{V}_{+} \to \bar{V}_{-} $ such that
\begin{enumerate}
\item
the graph $ L_{A} $ of $ A $ is comparable with $ \bar{L} $;
\item
$ \operatorname{reldim}\left(\bar{L},L_{A}\right)=\operatorname{ind}\psi+\delta $; here $ \delta $ is 0 if $ \psi\equiv \operatorname{const} $, and is 1 otherwise.
\end{enumerate}

If $ \varphi $ is compatible with $ E $, $ E' $ and $ \psi\equiv \operatorname{const} $, one can chose $ A $ so that
$ L_{A}=\bar{L} $. In particular, this is so if $ C $, $ C' $ are of genus 0. \end{lemma}

\begin{proof} By the closed graph theorem, it is enough to show that $ \bar{L} $ and $ \bar{V}_{-} $
are quasi-complementary, and calculate the excess. In turn, it is enough
to do the same with $ L $ and the preimage $ V_{-} $ of $ \bar{V}_{-} $ w.r.t.~the projection $ V
\to \bar{V} $. (It is this step what introduces $ \delta $ into the statement.)

We defined a line bundle $ {\mathcal L} $ on $ C^{*} $; its (continuous) global section $ f $
of $ {\mathcal L} $ corresponds to two continuous functions $ F $, $ F' $ on $ D_{-} $, $ D_{-}' $, which
satisfy the gluing relationship $ f|_{\gamma} = \psi\cdot\varphi^{*}\left(f'|_{\gamma'}\right) $. Obviously, $ \deg {\mathcal L}=\operatorname{ind}\psi $.

By Lemma~\ref{lm351.20}, there is a natural holomorphic structure on $ {\mathcal L} $.
Let $ \widehat{C}=C\cup C' $, $ \widehat{D}_{-}=D_{-}\cup D_{-}' $, $ \widehat{\gamma}=\gamma\cup\gamma' $, and $ \gamma^{*} $ be the common image of $ \gamma $, $ \gamma' $ in $ C^{*} $.
Consider the complexes
\begin{align} {\mathcal C}_{1}\colon H^{1}\left(C^{*},{\mathcal L}\right) & \xrightarrow[]{\bar{\partial}} H^{0}\left(C^{*},{\mathcal L}\otimes\bar{\omega}\right)/E
\notag\\
{\mathcal C}_{2}\colon H^{1}\left(\widehat{D}_{-}\subset\widehat{C},{\mathcal O}\right) & \xrightarrow[]{\bar{\partial}} H^{0}\left(\widehat{D}_{-}\subset\widehat{C},\bar{\omega}\right)/\widehat{E};
\notag\end{align}
here $ \widehat{E}=E\oplus E' $; clearly, $ \widehat{E} $ may be considered as a subspace of both the
$ H^{0} $-spaces. By the (finite-genus) Riemann--Roch theorem $ \bar{\partial} $ in $ {\mathcal C}_{1} $ is
Fredholm of index $ \deg {\mathcal L}+1 $; additionally, $ \operatorname{Im}\bar{\partial} $ is given by a finite
number of independent equations: $ \beta\in\operatorname{Im}\bar{\partial} $ iff $ \left<\alpha,\beta\right>=0 $, here $ \left<\alpha,\beta\right>=\int_{C}\alpha\beta $,
and $ \beta $ is a holomorphic section of $ {\mathcal L}^{*}\otimes\omega $.

By the definition of $ {\mathcal L} $, there is an inclusion $ \iota $ of the first complex
into the second one; denote components of $ \iota $ by $ \iota_{{\mathcal O}} $, $ \iota_{\bar{\omega}} $. Obviously, $ \iota_{\bar{\omega}} $ an
isomorphism, and $ \iota_{{\mathcal O}} $ is a closed inclusion. Since $ \bar{\partial}_{{\mathcal C}_{2}} $ is surjective, $ \operatorname{Ker}
\bar{\partial}_{{\mathcal C}_{2}} $ and the image $ L^{*} $ of $ H^{1}\left(C^{*},{\mathcal L}\right) $ in $ H^{1}\left(\widehat{D}_{-}\subset\widehat{C},{\mathcal O}\right) $ are quasi-complementary of
excess $ \deg {\mathcal L}+1 $.

Let $ W_{{\mathcal L}} $ and $ W_{{\mathcal O}} $ be the subspaces of $ H^{1}\left(C^{*},{\mathcal L}\right) $ and $ H^{1}\left(\widehat{D}_{-}\subset\widehat{C},{\mathcal O}\right) $ consisting
of functions vanishing at $ \gamma^{*} $ and $ \widehat{\gamma} $ correspondingly. By Lemma~\ref{lm351.20} $
\iota_{{\mathcal O}} $ identifies $ W_{{\mathcal L}} $ and $ W_{{\mathcal O}} $. Since $ L^{*}\supset W_{{\mathcal O}} $, the subspaces $ \left(W_{{\mathcal O}}+\operatorname{Ker} \bar{\partial}_{{\mathcal C}_{2}}\right)/W_{{\mathcal O}} $ and
$ L^{*}/W_{{\mathcal O}} $ of $ H^{1}\left(D_{-}\subset C,{\mathcal O}\right)/W_{{\mathcal O}} $ are quasi-complementary of excess $ \deg {\mathcal L}-g\left(C^{*}\right)+1 $.

On the other hand, restriction to $ \widehat{\gamma} $ identifies $ H^{1}\left(\widehat{D}_{-}\subset\widehat{C},{\mathcal O}\right)/W_{{\mathcal O}} $ with
$ H^{1/2}\left(\widehat{\gamma}\right) $. Since this identification sends $ \left(W_{{\mathcal O}}+\operatorname{Ker} \bar{\partial}_{{\mathcal C}_{2}}\right)/W_{{\mathcal O}} $ to $ H_{-}^{1/2}\left(\widehat{\gamma}\right) $, and
$ L^{*}/W_{{\mathcal O}} $ to $ L $, this finishes the proof of the first part of the lemma.

It is clear that $ L\cap V_{-} $ consists of global solutions of $ \bar{\partial}\varphi\in\widehat{E} $ in
sections of $ {\mathcal L} $. Thus if $ \varphi $ is compatible with $ E $, $ E' $ and $ \psi\equiv \operatorname{const} $, $ L\cap V_{-} $
consists of constants, and $ \bar{L}\cap\bar{V}_{-}=0 $. In other words, $ \bar{L} $ and $ \bar{V}_{-} $ are
complementary, hence $ \bar{L} $ is a graph of a continuous mapping $ \bar{V}_{+} \to \bar{V}_{-} $. \end{proof}

\begin{definition} If $ \varphi $ is compatible with $ E $, $ E' $ and $ \psi\equiv \operatorname{const} $, the {\em distortion\/}
$ \Delta\left(E,E',\gamma,\gamma',\varphi,\psi\right) $ is the norm of the operator $ A $ from the last statement of the
lemma. Otherwise put $ \Delta\left(E,E',\gamma,\gamma',\varphi,\psi\right) $ to be 1. \end{definition}

In this definition the curves $ \gamma $, $ \gamma' $ are considered together with
inclusions into complex curves $ C $, $ C' $ with a choice of excess spaces $ E $,
$ E' $.

\begin{remark} There exist excess spaces which are compatible with arbitrary
diffeomorphisms $ \varphi $. Indeed, given a curve $ C $ with an admissible real curve
$ \gamma\subset C $, perform $ g=g\left(C\right) $ cuts on $ D_{-} $ along real disjoint curves $ B_{1},\dots ,B_{g} $ so
that the resulting curve is of genus 0. Let $ W $ be the vector space of
linear functionals on $ L_{2}\left(D_{-},\omega\right) $ spanned by integrals along these curves.
Call $ E $ {\em cyclic\/} if pairing with $ E $ induces the same space $ W $ of linear
functionals on $ {\mathcal H}^{0}\left(D_{-}\subset E,\omega\right) $ (here $ {\mathcal H}^{0} $ is $ \operatorname{Ker}\bar{\partial} $ in $ H^{0} $). Obviously, if $ E $, $ E' $
are cyclic, they are compatible with arbitrary gluing $ \varphi $.

It is possible to construct cyclic excess spaces in any situation we
deal with in this paper. Indeed, if there is a annulus $ U $ around the cycle
$ B_{s} $ which lies completely inside $ D\subset C $, then the form $ \bar{\partial}\arg \left(z\right) $ on $ U $
represents the integration along the cycle.

Later we see that in the foam situation, as in Section~\ref{s6.3}, such
an annulus may not exist. Sketch how to deal with foam curves glued of
curves of non-0 genus; since one needs no such arguments if $ g\left(C\right)=0 $, and
with additional cuts one can always achieve this, we do not discuss these
arguments in more details. First, one can find an annulus $ U\subset C $ around a
cycle $ B $ such that it contains any connected component of $ C\smallsetminus D $ which
intersects $ U $. One may assume \cite{CourHurv} that $ U\cap D $ is conformally
equivalent to an annulus with cuts along concentric arcs. Now $ \bar{\partial}\arg \left(z\right)|_{U} $
is the suitable element of $ E $. \end{remark}

\subsection{Global sections }\label{s40.20}\myLabel{s40.20}\relax  Recall that $ H_{\text{loc}}^{1}\left(C,{\mathcal L}\right) $ was defined as the
vector
subspace of $ \prod_{k}H^{1}\left(D_{k}\subset C_{k},{\mathcal O}\right) $ consisting of collections of functions $ \left(F_{k}\right) $
which are compatible with the gluing data. As suggested by arguments in
Section~\ref{s10.80}, introduce the following refinement of Definition~%
\ref{def351.15}:

\begin{definition} Given a curve and bundle data $ \left(C,{\mathcal L}\right) $, and a mismatch
allowance $ \left(V_{j}\right) $ for these data such that $ \operatorname{const}\in V_{j} $ for any $ j $, define a norm
$ \|F\| $ of an $ H_{\text{loc}}^{1} $-section $ F=\left(F_{k}\right) $ modulo $ \left(V_{j}\right) $ as $ \|F\|^{2}\buildrel{\text{def}}\over{=}\sum_{k}\|F_{k}\|^{2} $; here $ \|F_{k}\| $
is taken in $ H^{1}\left(D_{k}\subset C_{k},{\mathcal O}\right)/\operatorname{const} $; call $ F $ with $ \|F\|<\infty $ a {\em global section modulo\/}
$ \left(V_{j}\right) $. Denote the space of such global sections by $ H_{\left[V^{\bullet}\right]}^{1}\left(C,{\mathcal L}\right) $. Define
similarly $ H^{0}\left(C,{\mathcal L}\otimes\bar{\omega}\right)\buildrel{\text{def}}\over{=}\bigoplus_{l_{2}}H^{0}\left(D_{k}\subset C_{k},{\mathcal O}\right) $ and
$ H_{E}^{0}\left(C,{\mathcal L}\right)\buildrel{\text{def}}\over{=}\bigoplus_{l_{2}}H^{0}\left(D_{k}\subset C_{k},{\mathcal O}\right)/E_{k} $. \end{definition}

Obviously, the operator $ \bar{\partial} $ induces an operator $ H_{\left[V^{\bullet}\right]}^{1}\left(C,{\mathcal L}\right) \to
H^{0}\left(C,{\mathcal L}\otimes\bar{\omega}\right) $, and $ H_{\left[V^{\bullet}\right]}^{1}\left(C,{\mathcal L}\right)\cap{\mathcal H}_{\text{loc}}^{1}\left(C,{\mathcal L}\right) $ coincides with the kernel of this
operator. We denote this operator by the same symbol $ \bar{\partial} $; however, note
that this operator has a slightly different semantic than the operator $ \bar{\partial} $
acting on generalized functions. If $ F\in H_{\left[V^{\bullet}\right]}^{1}\left(C,{\mathcal L}\right) $, then ``honest'' $ \bar{\partial}F $ can
be represented as a sum of an $ L_{2} $-section of $ {\mathcal L}\otimes\bar{\omega} $, and of generalized
functions with support on the curves $ \gamma_{j} $, $ j\in J $; each of these generalized
functions corresponds to an element of $ V_{j} $. When we consider $ \bar{\partial}:
H_{\left[V^{\bullet}\right]}^{1}\left(C,{\mathcal L}\right) \to H^{0}\left(C,{\mathcal L}\otimes\bar{\omega}\right) $, we keep only the $ L_{2} $-component.

Similarly, $ H_{\left[V^{\bullet}\right]}^{1}\left(C,{\mathcal L}\otimes\bar{\omega}\right)\cap{\mathcal H}_{\text{loc,}E}^{1}\left(C,{\mathcal L}\right) $ is $ \operatorname{Ker} \bar{\partial}\colon H_{\left[V^{\bullet}\right]}^{1}\left(C,{\mathcal L}\right) \to
H_{E}^{0}\left(C,{\mathcal L}\otimes\bar{\omega}\right) $. In this flavor of the $ \bar{\partial} $-operator we forget not only about the
$ \delta $-function components on the curves $ \gamma_{j} $, but also about the ``components''
of $ \bar{\partial}F $ in $ \oplus_{l_{2}}E_{k} $.

\begin{remark} The consideration of the index of the latter flavor of the
$ \bar{\partial} $-operator is the principal target of this paper. Is is forgetting about
$ \delta $-function components and $ E $-components which allows $ \bar{\partial} $ to have a finite
index (in ``good'' situations we are going to describe soon). Indeed, we
expect the ``honest'' operator $ \bar{\partial} $ to have an infinite negative index; taking
a quotient of the target space increases this index by exactly an amount
needed for it to become finite. \end{remark}

\subsection{Reduction to boundary }\label{s40.30}\myLabel{s40.30}\relax  Suppose that the excess subspaces $ E_{k} $ for
curves
$ C_{k} $ are $ D_{k} $-supported. Here we introduce two subspaces $ W_{E,\text{an}} $ and $ W_{\varphi\psi V} $ of an
appropriate Hilbert space $ W $; the relative position of these subspaces is
going to encode all the information about the operator $ \bar{\partial} $ we need.

By Theorem~\ref{th35.157}, there is a natural bounded operator $ \beta $ of
taking the boundary value modulo constants,
\begin{equation}
\beta\colon \bigoplus_{l_{2},k}H^{1}\left(D_{k}\subset C_{k},{\mathcal O}\right)/\operatorname{const} \to \bigoplus_{l_{2},j}H^{1/2}\left(\gamma_{j}\right)/\operatorname{const},
\notag\end{equation}
here each summand $ H^{1/2}\left(\gamma_{j}\right)/\operatorname{const} $ is equipped with the $ + $-skewed norm. Let
$ W=\bigoplus_{l_{2},j}H^{1/2}\left(\gamma_{j}\right)/\operatorname{const} $, let $ \beta_{E,\text{an}} $ be the restriction of $ \beta $ to
$ \bigoplus_{l_{2},k}{\mathcal H}_{E}^{1}\left(D_{k}\subset C_{k},{\mathcal O}\right)/\operatorname{const} $, $ W_{E,\text{an}}=\operatorname{Im}\beta_{E,\text{an}} $.

The vector space $ W $ is naturally decomposed into a sum of vector
subspaces $ W_{\pm} $, as in Section~\ref{s2.3}. If the distortions $ \Delta\left(E_{k}\right) $ of the excess
spaces are bounded, the component $ \beta_{-} $ of $ \beta_{\text{an}} $ corresponding to this
decomposition is invertible; thus $ \beta $ and $ \beta_{\text{an}} $ are monomorphisms; in
particular, $ W_{E,\text{an}} $ is a closed subspace. Consider the vector subspace
$ W_{\varphi\psi V}\subset W $ consisting of functions compatible with the gluing conditions up
to elements of allowance spaces. $ W_{\varphi\psi V} $ is always closed as a direct sum of
closed subspaces.

Note that a component $ H^{1/2}\left(\gamma_{j}\right)/V_{j} $ of the target space of the
mismatch operator $ \mu $ from Section~\ref{s10.80} can be identified with the
quotient of $ H^{1/2}\left(\gamma_{j}\right)\oplus H^{1/2}\left(\gamma_{j'}\right) $ by $ W_{\varphi\psi V}^{\left[j\right]} $, which is $ V_{j}\oplus V_{j'} $ summed with
the graph of the operator $ f \mapsto \psi_{j}\cdot\varphi_{j}^{*}f $. This defines a new Hilbert norm
on $ H^{1/2}\left(\gamma_{j}\right)/V_{j} $. In what follows we are going to use {\em this norm\/} when we
consider the $ l_{2} $-sum $ \bigoplus_{l_{2},\left\{j,j'\right\}}H^{1/2}\left(\gamma_{j}\right)/V_{j} $. Now this sum can be
identified with $ W/W_{\varphi\psi V} $.

\begin{remark} Another reason to introduce the new norm is dictated by
Theorem~\ref{th4.90} (see also the discussion which follows the theorem). Note
that if all the distortions $ \Delta\left(E,E',\gamma,\gamma',\varphi,\psi\right) $ (defined in Section~\ref{s4.11}) are
bounded, then this norm is equivalent to the old one. \end{remark}

\begin{theorem} \label{th52.40}\myLabel{th52.40}\relax  Suppose that all the spaces $ V_{j} $ but a finite number are
spanned by constants, and the distortions $ \Delta\left(E_{k}\right) $ of the excess spaces are
bounded. Then the following conditions are equivalent:
\begin{enumerate}
\item
the operator $ \bar{\partial}\colon H_{\left[V^{\bullet}\right]}^{1}\left(C,{\mathcal L}\right) \to H_{E}^{0}\left(C,{\mathcal L}\otimes\bar{\omega}\right) $ is Fredholm of index $ d $;
\item
the mismatch operator $ \mu\colon \bigoplus_{l_{2},k}{\mathcal H}_{E}^{1}\left(D_{k}\subset C_{k},{\mathcal O}\right)/\operatorname{const} \to
\bigoplus_{l_{2},\left\{j,j'\right\}}H^{1/2}\left(\gamma_{j}\right)/V_{j} $ is Fredholm of index $ d $;
\item
the vector subspaces $ W_{E,\text{an}}\subset W $ and $ W_{V\varphi\psi}\subset W $ defined above are
quasi-complementary of excess $ d $.
\end{enumerate}
\end{theorem}

\begin{proof} We know that $ W_{E,\text{an}} $ is identified with $ \bigoplus_{l_{2},k}{\mathcal H}_{E}^{1}\left(D_{k}\subset C_{k},{\mathcal O}\right)/\operatorname{const} $.
This identifies the operator $ \mu $ with the projection $ W_{E,\text{an}} \to W/W_{\varphi\psi V} $, which
proves the equivalence of the last two conditions.

On the other hand, the operator
\begin{equation}
\bar{\partial}_{Y}\colon Y=\bigoplus_{l_{2},k}H^{1}\left(D_{k}\subset C_{k},{\mathcal O}\right)/\operatorname{const} \to \bigoplus_{l_{2}}H^{0}\left(D_{k}\subset C_{k},{\mathcal O}\right)/E_{k}=H_{E}^{0}\left(C,{\mathcal L}\otimes\bar{\omega}\right)
\notag\end{equation}
is an epimorphism. This the information about the operator $ \bar{\partial}\colon H_{\left[V^{\bullet}\right]}^{1}\left(C,{\mathcal L}\right)
\to H_{E}^{0}\left(C,{\mathcal L}\otimes\bar{\omega}\right) $ is encoded in the relative position of the subspaces
$ H_{\left[V^{\bullet}\right]}^{1}\left(C,{\mathcal L}\right) $ and $ \operatorname{Ker}\bar{\partial}_{Y} $ in $ Y $; note that $ \operatorname{Ker}\bar{\partial}_{Y}= \bigoplus_{l_{2},k}{\mathcal H}_{E}^{1}\left(D_{k}\subset C_{k},{\mathcal O}\right)/\operatorname{const} $.
Using $ \beta_{-} $, the intersection of these subspaces is identified with
$ W_{E,\text{an}}\cap W_{V\varphi\psi} $. Since $ H^{1}\left(C,{\mathcal L}\right)\subset Y $ coincides with $ \beta^{-1}W_{\varphi\psi V} $, the sum of $ H^{1}\left(C,{\mathcal L}\right) $
and $ \operatorname{Ker}\bar{\partial}_{Y} $ is closed iff $ W_{\varphi\psi V}+\beta\left(\operatorname{Ker}\bar{\partial}_{Y}\right) =W_{\varphi\psi V}+W_{E,\text{an}} $ is closed; moreover,
the corresponding codimensions coincide. \end{proof}

Consider the case when $ V_{j} $ is spanned by 1 and $ \psi_{j} $. If $ \psi_{j}\equiv \operatorname{const} $, $ \dim 
V_{j}=1 $, otherwise $ \dim  V_{j}=2 $. In the notations above $ \beta $ identifies the vector
space of global analytic sections of $ {\mathcal L} $ modulo $ \left(V_{j}\right) $ with $ W_{E,\text{an}}\cap W_{V\varphi\psi} $. Since
the mismatch vanishes exactly on collections of functions compatible with
the gluing data, the vector space of possible mismatches modulo $ V_{j} $ of
collections of functions from $ \bigoplus_{l_{2},k}H^{1}\left(D_{k}\subset C_{k},{\mathcal O}\right)/\operatorname{const} $ may be identified
with $ W/W_{V\varphi\psi} $.

Use the notations of Section~\ref{s10.80}. By arguments of Section~%
\ref{s10.80}, the finite genus Riemann--Roch theorem may be rewritten in the
following way: the vector subspaces $ W_{E,\text{an}} $ and $ W_{V\varphi\psi} $ are quasi-complementary
of excess $ d+\sum_{\left\{j,j'\right\}}\left(\dim  V_{j}-1\right) $. Codify this observation in the
arbitrary-genus case by

\begin{definition} Say that curve and bundle gluing data
$ \left(D_{k}\subset C_{k},\gamma_{j},',\varphi_{j},\psi_{j}\right)_{k\in K,j\in J} $ with excess spaces $ \left(E_{k}\right) $ {\em satisfy Riemann\/}--{\em Roch
theorem\/} if $ \psi_{j}\equiv \operatorname{const} $ for all but a finite number of $ j\in J $, the distortions
$ \Delta\left(E_{k}\right) $ are uniformly bounded, and the conditions of Theorem~\ref{th52.40} are
satisfied with $ d=\sum_{\left\{j,j'\right\}}\left(\operatorname{ind}\psi_{j}+\dim  V_{j}-1\right) $; here $ V_{j} $ is spanned by 1 and $ \psi_{j} $.
\end{definition}

\subsection{Riemann--Roch theorem and operators $ {\protect \mathcal C} $ and $ {\protect \mathcal R} $ }\label{s352}\myLabel{s352}\relax  Now, when we
formulated the requirements of Riemann--Roch theorem, it makes sense to
investigate when these conditions hold. It is the last alternative of
Theorem~\ref{th52.40} that we are going to check.

Recall that the vector space $ W $ is naturally decomposed into a sum of
vector subspaces $ W_{\pm} $. Since $ \beta_{-} $ is an isomorphism provided the distortions
$ \Delta\left(E_{k}\right) $ are bounded, $ W_{\text{an}} $ is a graph of continuous mapping $ {\mathcal C}=\beta_{+}\circ\beta_{-}^{-1}\colon W_{-} \to
W_{+} $. By Lemma~\ref{lm30.70}, if the excess spaces are compatible with the
gluing mappings,\footnote{In fact, it is enough if the excess spaces are compatible with all but a
finite number of gluing mappings. We will not repeat this remark in what
follows.} $ W_{V\varphi\psi} $ is comparable with a graph of a closed mapping $ {\mathcal R}:
W_{+} \to W_{-} $, with the relative dimension being $ \sum_{\left\{j,j'\right\}}\left(\operatorname{ind}\psi_{j}+\dim  V_{j}-1\right) $. If
distortions $ \Delta\left(E_{k,k'},\gamma_{j,j'},\varphi_{j},\psi_{j,j'}\right) $ are uniformly bounded, the operator $ {\mathcal R} $
is continuous.

\begin{theorem}[abstract Riemann--Roch theorem]  \label{th4.90}\myLabel{th4.90}\relax  Consider two vector
subspaces $ V_{1,2}\subset H $ of a Hilbert space $ H=H_{1}\oplus H_{2} $. Suppose that $ V_{1} $ is
comparable with the graph of a closed mapping $ A_{1}\colon H_{1} \to H_{2} $ and the
relative dimension of $ V_{1} $ and this graph is $ d_{1} $. Suppose $ V_{2} $ is comparable
with the graph of a bounded mapping $ A_{2}\colon H_{2} \to H_{1} $ and the relative
dimension of $ V_{2} $ and this graph is $ d_{2} $. If $ A_{1}\circ A_{2} $ is defined everywhere, and
either $ A_{1}\circ A_{2} $ is compact, or $ A_{2}\circ A_{1}|_{\text{Dom }A_{1}} $ is compact, then $ V_{1} $ and $ V_{2} $ are
quasi-complementary {\em in the space\/} $ \operatorname{Dom}\left(A_{1}\right)\oplus H_{2} $ with the excess $ d_{1}+d_{2} $. \end{theorem}

We postpone the discussion of this (more or less trivial)
generalization of Fredholm theorem until the Appendix (Section~\ref{h4}). What
is important for us now is the fact that this statement is not
invertible, but is very close to be such (see Remark~\ref{rem6.60}). Moreover,
note that for a bounded operator $ A_{1} $ this is just a reformulation of
Fredholm theorem. Until Section~\ref{s3.45}, we concentrate our attention on the
cases when $ A_{1} $ is bounded. (If $ A_{1} $ is not bounded, then $ \operatorname{Dom}\left(A_{1}\right) $ is
considered as a Hilbert space, the norm being $ \|v\|_{A_{1}}^{2}\buildrel{\text{def}}\over{=}\|v\|^{2}+\|Av\|^{2} $.)

Now, to show that a Riemann--Roch theorem holds, it is enough to
investigate the compositions $ {\mathcal R}\circ{\mathcal C} $ and $ {\mathcal C}\circ{\mathcal R} $ for compactness. Consider the
operators $ {\mathcal R} $ and $ {\mathcal C} $ separately.

\subsection{Properties of the operator $ {\protect \mathcal R} $ }\label{s27}\myLabel{s27}\relax  From now on, we suppose that the
excess spaces are compatible with the gluing data, as defined in Section~%
\ref{s4.11}. In such a case the operator $ {\mathcal R} $ is defined up to addition of a
finite-dimensional operator.

With respect to $ j $-decomposition, the operator $ {\mathcal R} $ has a very simple
block structure: with a suitable numeration of blocks (separately in $ W_{+} $
and $ W_{-} $) and merging of 2 blocks corresponding to $ j $ and $ j' $, the operator $ {\mathcal R} $
becomes block-diagonal. Each block of this operator depends only on $ \gamma_{j} $,
$ \gamma_{j'} $ (considered with inclusions into the corresponding curves $ C_{k} $, $ C_{\widetilde{k}} $) and
$ \varphi_{j} $, $ \psi_{j} $ for one pair $ \left(j,j'\right) $ of matching indices. By Lemma~\ref{lm30.70}, each
block is bounded. Recall that the norm of this block is denoted by
$ \Delta\left(\gamma,\varphi,\psi\right) $.

In this section we are interested in describing when the operator $ {\mathcal R} $
is bounded. Since considering this condition one can drop an arbitrary
finite collection of blocks of $ {\mathcal R} $, it is enough to consider the case when
$ \psi_{j}\equiv \operatorname{const} $.

It is clear that $ \Delta\left(\gamma,\varphi,\psi\right)\leq\max \left(\psi,\psi^{-1}\right)\Delta\left(\gamma,\varphi,1\right) $. Thus if all the
constants $ \psi_{j} $ are bounded, it is enough to consider the case when $ \psi_{j}=1 $ for
any $ j $. There is one case when it is easy to estimate $ \Delta\left(\gamma,\varphi,1\right) $: when for
any $ j $ the identification $ \varphi_{j} $ of curves $ \gamma_{j} $ and $ \gamma_{j'} $ can be extended to a
conformal mapping of the corresponding complex curves $ C_{k} $ and $ C_{k'} $, both of
genus 0. Obviously, in such a case the corresponding block of the
operator $ {\mathcal R} $ can be written as a composition of two operators: the
inclusion $ \iota_{j}^{\text{skew}} $ of $ H_{+}^{1/2}\left(\gamma_{j}\right)/\operatorname{const} $ (with the $ + $-skewed norm) into
$ H^{1/2}\left(\gamma_{j}\right)/\operatorname{const} $ (with non-skewed norm); and the unitary identification $ \varphi_{j*} $
of $ H^{1/2}\left(\gamma_{j}\right)/\operatorname{const} $ and $ H^{1/2}\left(\gamma_{j'}\right)/\operatorname{const} $ (with non-skewed norms) via $ \varphi_{j} $.

By Remark~\ref{rem131.90}, the operators $ \iota_{j}^{\text{skew}} $ are proportional to
unitary operators if all the curves $ \gamma_{j} $ are circles. In this case, $ {\mathcal R} $ is
bounded iff $ |\psi_{j}| $ are uniformly bounded.

\begin{definition} Curve gluing data is {\em Schottky\/} if all diffeomorphisms $ \varphi_{j} $ may
be extended to a conformal transformation from $ C_{k} $ to $ C_{k'} $; here $ \gamma_{j}\subset C_{k} $,
$ \gamma_{j'}\subset C_{k'} $. Curve gluing data is {\em circular\/} if all the curves $ C_{k} $ are
isomorphic to $ {\mathbb C}{\mathbb P}^{1} $ and all the curves $ \gamma_{j} $ are circles. \end{definition}

Give modifications of these notions more suitable to our situation:

\begin{definition} \label{def27.20}\myLabel{def27.20}\relax  Call the norm of the operator $ \varphi_{j*} $ the {\em distortion\/} $ \Delta\left(\varphi_{j}\right) $
{\em of the diffeomorphism\/} $ \varphi_{j} $. Curve gluing data is {\em quasi-Schottky\/} if the set
$ \left\{\Delta\left(\varphi_{j}\right) \mid j\in J\right\} $ is bounded. Curve gluing data is {\em quasi-circular\/} if the set
$ \left\{\Delta\left(\gamma_{j}\right) \mid j\in J\right\} $ is bounded. The bundle gluing data is {\em bounded\/} if all the
functions $ \psi_{j} $ are uniformly bounded. \end{definition}

We conclude that the operator $ {\mathcal R} $ is always a direct sum of bounded
operators, thus is a closed operator. Moreover, for a bounded bundle
gluing data for quasi-circular quasi-Schottky curve gluing data this
operator is a bounded invertible operator. In such a case, the question
of compactness of $ {\mathcal R}\circ{\mathcal C} $ and $ {\mathcal C}\circ{\mathcal R} $ is reduced to the question of compactness
of $ {\mathcal C} $.

One expects that curve gluing data being quasi-Schottky (or even
Schottky) is not a very significant restriction: any ``reasonable'' curve
should have such a ``representation''. Later we show that quasi-circularity
condition may be circumvented. On the other hand, unbounded bundle gluing
data comes very naturally in the theory of divisors; we postpone
discussion of such gluing data until Section~\ref{s3.45}.

\subsection{Strong Riemann--Roch conditions } By definition, $ {\mathcal C}=\beta_{+}\circ\beta_{-}^{-1} $. The
operator $ \beta_{-} $ is unitary if all the curves $ C_{k} $ are of genus 0; if $ K=\left\{1\right\} $ (so
there is only one curve $ C_{1} $) then one can bound the operators $ \beta_{-} $ and $ \beta_{-}^{-1} $
by a constant depending only on the distortion $ \Delta\left(E\right) $ of the excess spaces
$ E $ for the curve $ C_{1} $.

\begin{definition} A collection of compact curves $ C_{k} $, $ k\in K $, with excess spaces $ E_{k} $ is
{\em tame\/} if the set $ \left\{\Delta\left(E_{k}\right) \mid k\in K\right\} $ is bounded. \end{definition}

Since the operator $ \beta_{+} $ is always a contraction (due to our choice of
the $ + $-skewed norm), we obtain:

\begin{proposition} If the collection $ \left(C_{k}\right) $ is tame, then the operator $ {\mathcal C} $ is
continuous. \end{proposition}

\begin{corollary} \label{cor28.35}\myLabel{cor28.35}\relax  Suppose the curves $ \left(C_{k}\right) $ form a tame quasi-Schottky
quasi-circular collection, and the bundle gluing data is bounded. Then $ {\mathcal R} $
and $ {\mathcal C} $ are continuous operators, and the curve-bundle gluing data
satisfies Riemann--Roch theorem iff $ {\mathcal R}{\mathcal C}-1 $ is a Fredholm operator of index
0. In particular, this holds if 1 is not in the essential spectrum of the
operator $ {\mathcal R}\circ{\mathcal C} $. \end{corollary}

The property of 1 not being in the essential spectrum looks very
fragile w.r.t.~changes in the curve and bundle data. However, to
construct the Jacobian of a curve, it is useful to know that {\em all\/} the line
bundles which correspond to points of the Jacobian satisfy Riemann--Roch
theorem. Remark~\ref{rem6.60} suggests another, much more robust criterion:

\begin{definition} Suppose that excess spaces of the curve-bundle data have
uniformly bounded distortions. Say that the curve-bundle gluing data
{\em satisfies the strong Riemann\/}--{\em Roch theorem\/} if $ {\mathcal R}\circ{\mathcal C} $ is defined everywhere,
and is a compact operator. \end{definition}

\begin{theorem} \label{th28.40}\myLabel{th28.40}\relax  Suppose that the curve-bundle gluing data satisfies
the strong Riemann--Roch theorem. Consider another curve-bundle gluing data
which differs from the initial data only by replacing functions $ \psi_{j} $ by
functions $ \psi_{j}' $. If functions $ \psi'_{j}/\psi_{j} $ are uniformly bounded, then the
modified curve-bundle data also satisfies the strong Riemann--Roch theorem \end{theorem}

\begin{proof} As Section~\ref{s27} shows, the modified operator $ {\mathcal R}' $ differs from $ {\mathcal R} $
by a multiplication by a bounded operator. \end{proof}

Note that the last past of Remark~\ref{rem6.60} cannot be applied to
reverse Theorem~\ref{th28.40}. Indeed, the operators $ {\mathcal R} $ have a very special
form only. Moreover, one can construct counter-examples (see Section~%
\ref{s6.70}) when $ {\mathcal R} $ and $ {\mathcal C} $ are bounded, the composition $ {\mathcal R}\circ{\mathcal C} $ is not compact, but
has the essential spectrum $ \left\{0\right\} $ (e.g., $ \left({\mathcal R}\circ{\mathcal C}\right)^{2} $ may be compact). However,
all the counter-examples we know are not stable w.r.t.~appropriate
modifications of {\em curve\/} gluing data: if two components $ C_{k} $, $ C_{\widetilde{k}} $ are glued
together by identification of $ \gamma_{j}\subset C_{k} $, $ \gamma_{j'}\subset C_{\widetilde{k}} $, we can replace this pair by
one curve $ C^{*} $ which is result of this identification; similarly, one can
break a component into two, or chain several such operations.

It is natural to conjecture that the ``interesting'' curves, those
which correspond to points of moduli spaces lying in the support of
``interesting'' measures, can be described by gluing data satisfying the
strong Riemann--Roch theorem. In what follows we consider such
gluing data only.

\subsection{Operator $ {\protect \mathcal C} $ and bar-projectors }\label{s5.7}\myLabel{s5.7}\relax  Obviously, the operator $ {\mathcal C} $ depends
only on the domains $ D_{k}\subset C_{k} $ (and the corresponding excess spaces), but not
on the gluing data; moreover, it breaks into a direct sum of the
corresponding operators for the individual curves $ C_{k} $. Restrict our
attention to the operator $ {\mathcal C} $ for one such curve $ C $; in other words, we may
assume that $ K $ contains one element 1 only, and $ C_{1}=C $.

Consider now the finer block structure of the operator $ {\mathcal C} $ due to the
decomposition of $ W_{\pm} $ into a direct sum over indices $ j $ of curves $ \gamma_{j} $. Let
$ {\mathcal C}_{jl} $ be the component sending $ H_{-}^{1/2}\left(\gamma_{l}\right)/\operatorname{const} $ into $ H_{+}^{1/2}\left(\gamma_{j}\right)/\operatorname{const} $. It is
clear that $ {\mathcal C}_{jj}=0 $.

From now on assume that the curves $ \gamma_{j} $ do not intersect. In
particular, the sets $ D_{k} $ have no open subsets of dimension 1. Under this
condition we can find convenient bounds for the blocks $ {\mathcal C}_{jk} $, $ j\not=k $.

First of all, recall that a choice of the excess space $ E $ on $ C $ allows
one to construct the operator $ \bar{\partial}_{E}^{-1} $ mapping global sections of $ \bar{\omega} $ to global
sections of $ {\mathcal O} $ modulo constant. Thus the operator $ \partial\circ\bar{\partial}_{E}^{-1} $ maps sections of
$ \bar{\omega} $ to sections of $ \omega $.

\begin{definition} Consider a complex compact curve $ C $ with an excess space $ E $,
and two non-intersecting open subsets $ R,R'\subset C $. Define the {\em bar-projector\/}
$ {\mathcal P}_{RR'}\colon H^{0}\left(R',\bar{\omega}\right) \to H^{0}\left(R,\omega\right) $ as the composition $ \rho_{R}\circ\partial\circ\bar{\partial}_{E}^{-1}\circ\varepsilon_{R'} $; here $ \varepsilon_{R'}f $ is
the extension of $ f $ by 0 from $ R' $, and $ \rho_{R}g $ is the restriction of $ g $ to $ R $. \end{definition}

Let $ D' = C\smallsetminus R' $. Obviously, the bar-projector vanishes on $ \bar{\partial}\mathring{H}^{1}\left(R'\right) $;
call the induced operator from $ H^{0}\left(R',\bar{\omega}\right)/\bar{\partial}\mathring{H}^{1}\left(R'\right)={\mathcal R}_{D'\subset C,E}={\mathcal H}_{E}^{1}\left(D'\right)/\operatorname{const} $ by
the same term. When acting from $ {\mathcal H}_{E}^{1}\left(D'\right)/\operatorname{const} $, the bar-projector is
identified with $ \rho_{R}\circ\partial $, or, if $ E $ is $ C\smallsetminus R $-supported, with $ \rho_{R}\circ d_{\text{deRham}} $. On the
other hand, $ d_{\text{deRham}} $ is a unitary operator from $ H_{\text{int}}^{1}\left(R\right) $ to its image in
$ L_{2}\left(R,\Omega^{1}\right) $. In other words, the bar-projector may be also modeled as the
restriction operator $ {\mathcal H}_{E}^{1}\left(D'\right)/\operatorname{const} \to H_{\text{int}}^{1}\left(R\right)/\operatorname{const} $.

The arguments above show that $ {\mathcal C}_{jl} $ differs from the bar-projector
from $ R_{l} $ to $ R_{j} $ only by appropriate approximately-unitary transformation of
the image and the preimage. However, the bar-projector is explicitly written
as an off-diagonal block of a pseudo-differential operator $ \partial\circ\bar{\partial}^{-1} $; in
particular, the bar-projector is an operator with a smooth Schwartz
kernel, so is compact. We obtain

\begin{proposition} The block $ {\mathcal C}_{jk} $ is compact, and can be written as
$ \iota_{1}\circ{\mathcal P}_{R_{j}R_{l}}\circ\iota_{2} $; here $ \iota_{1,2} $ are invertible operators, and
$ \|\iota_{1,2}\|,\|\iota_{1,2}^{-1}\|\leq\Delta\left(E\right) $. \end{proposition}

\subsection{Conformal distance and estimates of bar-projectors }

\begin{definition} \label{def31.40}\myLabel{def31.40}\relax  Consider two non-intersecting regions $ R $, $ R' $ in a
complex curve $ C $ such that $ C\smallsetminus\left(R\cup R'\right) $ is conformally equivalent to a
cylinder $ S^{1}\times\left[0,\lambda\right] $ of radius 1 and length $ \lambda\geq0 $. If $ g\left(C\right)>0 $, suppose that the
support of the excess space $ E $ is in $ R\cup R' $. Call $ \lambda $ the {\em conformal distance\/}
between $ R $ and $ R' $.

For general non-intersecting regions $ R $, $ R' $ the {\em conformal distance\/} is
at least $ \lambda $ if after increasing one of the regions the condition above
holds. \end{definition}

The specific form of the last part of the definition is chosen to
simplify the proof of Proposition~\ref{prop31.60}, while allowing the
construction of Section~\ref{s6.3} to remain simple.

\begin{lemma} \label{lm31.50}\myLabel{lm31.50}\relax  Suppose $ C={\mathbb C}{\mathbb P}^{1} $, and both $ R $ and $ R' $ are disks on $ C $ with
conformal distance $ \lambda $. Then the norm of the bar-projector from $ R' $
to $ R $ is equal to $ e^{-\lambda} $. \end{lemma}

\begin{proof} Since $ E=0 $, the bar-projector is canonically defined. Due to its
conformal invariance, we may suppose $ R $ is given by $ |z|<e^{-\lambda} $, $ R' $ is given
by $ |z|>1 $. Due to rotation-invariance, we may consider each irreducible
component of the action of $ U\left(1\right) $ on $ {\mathcal H}^{1}\left(D'\right) $, $ D'=C\smallsetminus R' $, one-by-one. Clearly,
$ d\left(z^{k}\right)|_{R} $ has the squared $ L^{2} $-norm $ \pi ke^{-2k\lambda} $. The continuation $ f $ of $ z^{k}|_{\partial R'} $
into $ R' $ with the minimal $ L^{2} $-norm of $ \bar{\partial}f $ is $ \bar{z}^{-k} $; its squared norm on $ R' $ is
$ \pi k $. \end{proof}

Since the conformal distance between two disks in $ {\mathbb C} $ of radii $ r_{1,2} $
and the distance between centers $ d $ is $ \operatorname{ch}^{-1}\frac{d^{2}-r_{1}^{2}-r_{2}^{2}}{2r_{1}r_{2}} $, the norm of
the bar-projector is comparable with $ \frac{r_{1}r_{2}}{d\left(d-r_{1}-r_{2}\right)} $.

\begin{proposition} \label{prop31.60}\myLabel{prop31.60}\relax  There is a constant $ c\geq1 $ such that if the conformal
distance between simply-connected domains $ R,R'\subset C $ is at least $ \lambda $, then the
norm of the bar-projector from $ R' $ to $ R $ is less or equal to $ c\Delta\left(E\right) e^{-\lambda} $;
here we suppose that the excess space $ E $ for $ C $ is $ C\smallsetminus\left(R\cup R'\right) $-supported. \end{proposition}

\begin{proof} The operator $ \partial\bar{\partial}_{E}^{-1} $ is bounded by $ \Delta\left(E\right) $. Since the bar-projector
is a block of this operator, its norm is bounded by the norm of $ \partial\bar{\partial}_{E}^{-1} $.
Thus we may assume that $ \lambda\geq2 $. The norm of the bar-projector is monotonous
w.r.t.~the domains, increase one of the domains so that the the first
part of the definition of conformal distance holds. After this, one of
the regions $ R $ or $ R' $ contains the support of $ E $. At first, suppose that
this domain is $ R' $.

Our aim is to reduce the statement to one of Lemma~\ref{lm31.50}. The
obstacles are the presence of the excess space $ E $, that $ R' $ is not
necessarily simply connected, and even if $ R $ and $ R' $ were homeomorphic to
disks, their boundaries may be glued to the boundary of the annulus
$ {\mathcal A}=C\smallsetminus\left(R\cup R'\right) $ in a ``non-standard''way. Note that $ R $ is homeomorphic to a disk;
our first target is to show that we may suppose that $ R' $ is a disk and it
is glued to $ {\mathcal A} $ as in Lemma~\ref{lm31.50}.

Identify the annulus $ \left\{e^{-\lambda}\leq|z|\leq1\right\} $ with $ {\mathcal A}=C\smallsetminus\left(R\cup R'\right) $, glue the disk
$ \left\{|z|\geq1\right\}\subset{\mathbb C}{\mathbb P}^{1} $ to $ C\smallsetminus R' $ via this identification; denote the resulting curve
$ C_{00} $. Show that one may replace $ C $ by $ C_{00} $.

Let $ \widetilde{R}' $ be the image of the annulus $ \left\{e^{-1}\leq|z|\leq1\right\} $ in $ C $, $ \widetilde{R}'_{00} $ be the
image of this annulus in $ C_{00} $.

\begin{lemma} Given $ \alpha\in\mathring{H}^{0}\left(R',\bar{\omega}\right) $, one can find $ \widetilde{\alpha}\in\mathring{H}^{0}\left(\widetilde{R}',\bar{\omega}\right) $ such that $ \bar{\partial}_{E}^{-1}\left(\alpha-\widetilde{\alpha}\right) $
is constant on $ R $, $ \widetilde{\alpha}\in\operatorname{Im}\bar{\partial} $ and $ \bar{\partial}^{-1} \widetilde{\alpha} $ is constant on $ R' $, and $ \|\widetilde{\alpha}\|=O\left(\Delta\left(E\right)\|\alpha\|\right) $.
\end{lemma}

\begin{proof} Find a representative $ f\in H^{1}\left(C\right) $ of $ \bar{\partial}_{E}^{-1}\alpha\in H^{1}\left(C\right)/\operatorname{const} $. Normalize $ f $
by requiring $ \int_{\widetilde{R}'}f\,dz\,d\bar{z}=0 $. Let $ \Psi' $ be a cut-off function which is 0 on $ R' $,
1 outside of $ R'\cup\widetilde{R}' $. Then $ \Psi'f $ is holomorphic outside of $ \widetilde{R}' $. Since the norm
of $ f $ in $ H^{1}\left(C\right)/\operatorname{const} $ is $ O\left(\Delta\left(E\right)\|\alpha\|\right) $, to show a similar estimate for $ \Psi'f $, it
is enough to estimate $ \|f\|_{L_{2}\left(\widetilde{R}'\right)} $; this estimate follows from the estimate
of $ \|df\|_{L_{2}} $. Now take $ \widetilde{\alpha}=\bar{\partial}\left(\Psi'f\right) $. \end{proof}

Now replace $ R' $ by $ R'\cup\widetilde{R}' $, and $ \alpha $ by $ \widetilde{\alpha} $. This decreases $ \lambda $ by 1; thus it
is enough to show that the norm of the bar-projector reduced to $ \operatorname{Im}\bar{\partial} $ can
be estimated as $ O\left(e^{-\lambda}\right) $ (no factor $ \Delta\left(E\right) $). Moreover, since $ \operatorname{Supp}\widetilde{\alpha}\subset{\mathcal A} $, one can
identify $ \widetilde{\alpha} $ with a $ 1 $-form on $ C_{00} $; since $ \bar{\partial}^{-1} \widetilde{\alpha} $ is constant on $ R' $, one can
transfer $ \bar{\partial}^{-1} \widetilde{\alpha} $ to $ C_{00} $ as well. Thus we may assume $ C=C_{00} $, $ \alpha=\widetilde{\alpha} $; since
$ g\left(C_{00}\right)=0 $, there is no $ E $ to care about.

Glue the disk $ \left\{|z|\leq e^{-\lambda}\right\} $ to $ C\smallsetminus R $ via the identification of $ {\mathcal A} $ with
$ \left\{e^{-\lambda}\leq|z|\leq1\right\} $; denote the resulting curve $ C_{0} $. Let $ R_{0}\subset C_{0} $ be the image of the
disk $ \left\{|z|<e^{-\lambda}\right\} $, $ \widetilde{R}\subset C $ and $ \widetilde{R}_{0}\subset C_{0} $ be the images of the annulus
$ \left\{e^{-\lambda}\leq|z|\leq e^{-\lambda+1}\right\} $. By Lemma~\ref{lm31.50}, the norm of $ \partial_{C_{0}}\bar{\partial}_{C_{0}}^{-1}\alpha $ on $ R_{0}\cup\widetilde{R}_{0} $ is
$ O\left(e^{-\lambda}\|\alpha\|\right) $. Using again the trick with a cut-off function, we can find
$ \widetilde{\alpha}\in\mathring{H}^{0}\left(\widetilde{R}_{0},\bar{\omega}\right) $ such that $ \bar{\partial}^{-1}\left(\alpha-\widetilde{\alpha}\right) $ is constant on $ R_{0} $, and $ \|\widetilde{\alpha}\|=O\left(e^{-\lambda}\|\alpha\|\right) $.
Again, since $ \bar{\partial}^{-1}\left(\alpha-\widetilde{\alpha}\right) $ is constant on $ R_{0} $, we can transfer this function to
$ C $. Denote the resulting function $ f $; transfer similarly $ \widetilde{\alpha} $ to $ C $. Then
$ \bar{\partial}^{-1} \widetilde{\alpha}|_{R} $ differs from $ \bar{\partial}^{-1}\alpha|_{R} $ by a constant only; moreover,
$ \|\bar{\partial}^{-1} \widetilde{\alpha}|_{R}\|_{H^{1}/\operatorname{const}}=O\left(e^{-\lambda}\|\alpha\|\right) $.

What remains to prove is the other case, when $ \operatorname{Supp} E\subset R $. Construct $ C_{0} $
as above, replacing $ R $ by a disk $ R_{0} $; then $ g\left(C_{0}\right)=0 $. Given $ \alpha $ as above, apply
the already proved case $ g=0 $; there are a function $ f $ on $ C_{0} $ and a $ 1 $-form $ \widetilde{\alpha} $
with support on $ \widetilde{R}_{0} $ such that $ f|_{R_{0}}=\operatorname{const} $, $ \bar{\partial}f=\alpha-\alpha_{0} $, and $ \|\widetilde{\alpha}\|=O\left(e^{-\lambda}\right) $. One can
transfer $ f $ and $ \widetilde{\alpha} $ to $ C $; since one can replace $ \alpha $ by $ \widetilde{\alpha} $ and $ \lambda $ by 0, the
estimate $ O\left(\Delta\left(E\right)\right) $ for $ \partial\bar{\partial}_{E}^{-1} $ we already used finishes the proof. \end{proof}

This proof in fact implies a much stronger result:

\begin{amplification} The same estimate holds if one considers the
bar-projector acting into the space $ H^{1}\left(R\subset C\right) $ instead of $ H_{\text{int}}^{1}\left(R\right) $. \end{amplification}

\begin{remark} In fact the norm of the bar-projector may be much smaller than
what is given by Proposition~\ref{prop31.60}. Assume that $ C={\mathbb C}{\mathbb P}^{1} $, $ R'=\left\{|z|>1\right\} $.
Suppose that $ R $ sits inside the disk $ \left\{|z|<\rho\right\} $, $ \rho<1 $. Since we know the
(smooth) kernel $ c\,dz\,dz'/\left(z-z'\right)^{2} $ of the bar-projector, it is easy to
calculate the Hilbert--Schmidt norm of this operator, thus estimate its
norm. In fact, due to the conformal invariance of the bar-projector, the
square of its Hilbert--Schmidt norm is proportional to the area of $ R $ in
the hyperbolic metric of the disk $ \left\{|z|<1\right\} $, thus the bar-projector is
bounded by $ O\left(\sqrt{|R|}/\left(1-\rho^{2}\right)\right) $; here $ |R| $ is the Euclidean area of $ R $.

This argument shows that the norm of the bar-projector can be made
arbitrary small even under the requirement that $ R $ contains the interval
$ \left[-\varepsilon,\varepsilon\right] $, $ \varepsilon<1 $. However, the conformal distance between this interval and $ R' $
is finite.

Using the test-function $ z $ in $ \left\{|z|\leq1\right\} $, one can show that the norm of
the bar-projector is bounded from below by $ \sqrt{|R|} $. \end{remark}

\begin{remark} Note that if $ g\left(C\right)=0 $, the norm of the bar-projector from $ R' $ to $ R $
is equal to the one from $ R $ to $ R' $. Indeed, let $ X=\partial\bar{\partial}^{-1} $. Then $ X $
manifestly satisfies $ \bar{X}=X^{-1} $; moreover, $ X $ is unitary due to Remark~%
\ref{rem1.305}. Thus the operator $ X $ equals its transposed. Since the
bar-projector is a block of the operator $ X $, and the blocks which
correspond to projectors from $ R' $ to $ R $ and from $ R $ to $ R' $ are transposed,
this implies the result. \end{remark}

\begin{remark} Sometimes it is possible to calculate the norm of the
bar-projector explicitly. E.g., let $ C $ be the disk $ \left\{|z|\leq1\right\} $ with $ z $
identified with $ z^{-1} $ for any point $ z $ of the boundary; $ C $ is a rational
curve. Let $ R_{+}=\left\{|z|<a\right\} $, $ R\_ $ be the image of $ \left\{b<|z|\leq1\right\} $ on $ C $; here $ a<b<1 $.

Suppose that $ R\subset C $ has a smooth boundary. Given a function $ f $ in
$ {\mathcal H}^{1}\left(C\smallsetminus R\right) $, let $ \widetilde{f} $ be the harmonic extension of $ f|_{\partial R} $ into $ R $. It is clear that
$ \widehat{f} \buildrel{\text{def}}\over{=} \bar{\partial}\widetilde{f} $ is the representative of $ f $ in $ H^{0}\left(R,\bar{\omega}\right) $ with the minimal norm.
Take $ R=R_{-} $, $ f_{n}=-\left(b^{-4n}-1\right)z^{n} $; since harmonic functions in $ R_{-} $ are linear
combinations of $ z^{l}+z^{-l},\quad \bar{z}^{l}+\bar{z}^{-l} $, we conclude that
$ \widetilde{f}_{n}=\left(z^{n}+z^{-n}\right)-b^{-2n}\left(\bar{z}^{n}+\bar{z}^{-n}\right) $, thus $ \widehat{f}_{n}=-nb^{-2n}\left(\bar{z}^{n}-\bar{z}^{-n}\right)d\bar{z}/\bar{z} $. If $ n\not=n' $, $ f_{n} $ and $ f_{n'} $
are orthogonal in $ {\mathcal H}^{1}\left(C\smallsetminus R_{-}\right) $ (since $ \widehat{f}_{n} $, $ \widehat{f}_{n'} $ are orthogonal in $ H^{0}\left(R_{-},\bar{\omega}\right) $).
Since $ R_{+} $ is rotation-invariant, thus $ \partial f_{n} $ and $ \partial f_{n'} $ are orthogonal in
$ H^{0}\left(R_{+},\omega\right) $, one may consider the action of the projector on functions $ f_{n} $
one-by-one.

Since
\begin{equation}
\|\widehat{f}_{n}\|^{2} = 2\pi n^{2}b^{-4n}\int_{b}^{1}\left(r^{2n}+r^{-2n}\right)dr/r,\qquad \|\partial f_{n}|_{R_{+}}\|^{2} = 2\pi n^{2}\left(b^{-4n}-1\right)^{2}\int_{0}^{a}r^{2n}dr/r,
\notag\end{equation}
this component of the projector has the square of the norm $ \left(1-b^{4n}\right)e^{-2n\lambda} $;
here $ \lambda $ is the conformal distance $ \log \left(b/a\right) $ between $ R_{+} $ and $ R_{-} $. Thus the
norm of the projector from $ R_{-} $ to $ R_{+} $ is $ \sqrt{1-b^{4}}e^{-\lambda} $. This suggests that the
constant $ c $ of Proposition~\ref{prop31.60} may be 1. \end{remark}

\subsection{Block matrices }

\begin{definition} Consider Hilbert direct sums $ V=\bigoplus_{l_{2}}V_{l} $, $ V'=\bigoplus_{l_{2}}V_{k}' $ of Hilbert
spaces, and an operator $ A\colon V \to V' $. It induces operators $ A_{kl}\colon V_{l} \to V'_{k} $;
call these operators the {\em blocks\/} of $ A $. \end{definition}

\begin{lemma} \label{lm32.20}\myLabel{lm32.20}\relax  Use the notations of the definition above.
\begin{enumerate}
\item
if the matrix $ \left(\|A_{kl}\|\right) $ corresponds to a bounded operator $ l_{2} \to l_{2} $,
then the operator $ A $ is bounded.
\item
Suppose that all the blocks of $ A $ are compact operators. If the
matrix $ \left(\|A_{kl}\|\right) $ corresponds to a compact operator $ l_{2} \to l_{2} $, then the
operator $ A $ is compact.
\end{enumerate}
\end{lemma}

\begin{proof} If the matrix $ \left(\|A_{kl}\|\right) $ corresponds to a continuous operator $ l_{2}
\to l_{2} $ with a norm $ \leq M $, then $ \|A\|\leq M $. Indeed, otherwise one could find $ v\in V $,
$ v'\in V' $ with $ |v|=|v'|=1 $, $ |\left(Av,v'\right)|>M $. Let $ v=\left(v_{l}\right) $, $ v'=\left(v'_{k}\right) $; then $ \widetilde{v}_{l}=|v_{l}| $,
$ \widetilde{v}'_{k}=|v'_{k}| $ are in $ l_{2} $, and $ |\widetilde{v}_{\bullet}|=|\widetilde{v}'_{\bullet}|=1 $. Now $ |\left(Av,v'\right)|\leq\sum_{kl}\|A_{kl}\|\widetilde{v}_{l}\widetilde{v}_{k}' $
leads to a contradiction.

If $ \left(\|A_{kl}\|\right) $ gives a compact operator, it may be
approximated with arbitrary precision by replacing all but a finite number
of entries by 0. Thus $ A $ may be approximated by replacing all but a
finite number of blocks by 0 blocks. What remains is a finite sum of
compact operators, thus compact. Since $ A $ may be approximated by compact
operators, $ A $ is compact itself.

In the other direction, if $ A $ is compact, so are its blocks. To show
that the operator $ l_{2} \to l_{2} $ given by the matrix $ \left(\|A_{kl}\|\right) $ can be
approximated by finite-dimensional operators, approximate $ A $ by a
finite-dimensional operator; write down this operator as a sum of
one-dimensional operators $ \sum_{n}w'_{n}\otimes w_{n} $, here $ w'\otimes w $ sends $ v $ into $ \left(v,w\right)w' $. One
can approximate each of $ w_{n} $, $ w'_{n} $ by a finite sum of vectors in $ V_{l} $, $ V_{k}' $;
hence $ A $ may be approximated by an operator with only a finite number of
non-0 blocks; consequently, $ A $ may be approximated by a finite sum $ \widehat{A} $ of
its blocks. Then $ \|A-\widehat{A}\| $ can be made arbitrarily small , thus
$ \left|\left(\left(A-\widehat{A}\right)\sum\varepsilon_{l}v_{l},\sum\varepsilon'_{k}v_{k}'\right)\right| $ can be made arbitrarily small for any $ \left(\varepsilon_{l}\right) $,
$ \left(\varepsilon'_{k}\right) $ with $ |\left(\varepsilon_{l}\right)|_{l_{2}}=|\left(\varepsilon_{k}'\right)|_{l_{2}}=1 $. \end{proof}

\subsection{Practical criteria for Riemann--Roch theorem }\label{s54}\myLabel{s54}\relax  Consider the estimates
for the norms of blocks of the operator $ {\mathcal C} $ to obtain easy-to-check
criteria that a curve satisfies a strong Riemann--Roch theorem. Assume
that the gluing data is bounded and quasi-Schottky quasi-circular, so that
the operator $ {\mathcal R} $ is bounded and invertible. In such a case all we need to
show is that the operator $ {\mathcal C} $ is compact.

Given $ j,\widetilde{j}\in J $, $ e^{-\lambda_{j\widetilde{j}}} $ is well-defined if $ \gamma_{j} $ and $ \gamma_{\widetilde{j}} $ are on the same
curve $ C_{k} $;
here $ \lambda_{j\widetilde{j}} $ is the conformal distance between $ \gamma_{j} $ and $ \gamma_{\widetilde{j}} $. Put
$ e^{-\lambda_{j\widetilde{j}}}\buildrel{\text{def}}\over{=}0 $ otherwise.

\begin{theorem} \label{th33.35}\myLabel{th33.35}\relax  Consider a bounded quasi-Schottky quasi-circular curve
and bundle gluing data with uniformly bounded distortions $ \Delta\left(E_{k}\right) $ of the
excess spaces. Suppose that the excess spaces are compatible with
identifications $ \varphi_{j} $ of smooth boundary components, and that the smooth
components of boundaries of subsets $ D_{k} $ do not intersect. These data
satisfies the strong Riemann--Roch condition if the operator $ l_{2} \to l_{2} $
defined by the matrix $ \left(e^{-\lambda_{j\widetilde{j}}}\right) $ is bounded and compact. \end{theorem}

\begin{proof} By Proposition~\ref{prop31.60}, the blocks $ {\mathcal C}_{j\widetilde{j}} $ of the operator $ {\mathcal C} $ are
compact and have norms $ O\left(e^{-\lambda_{j\widetilde{j}}}\right) $. By Lemma~\ref{lm32.20}, in the conditions of
the theorem the operator $ {\mathcal C} $ is compact. Now Corollary~\ref{cor28.35} implies
the theorem. \end{proof}

\begin{corollary} \label{cor33.40}\myLabel{cor33.40}\relax  Consider a bounded quasi-Schottky quasi-circular curve
and bundle gluing data. Suppose that the excess spaces are compatible
with identifications $ \varphi_{j} $ of smooth boundary components, and that the
smooth components of boundaries of subsets $ D_{k} $ do not intersect. These
data satisfies the strong Riemann--Roch condition if $ \sum_{m\not=j}e^{-2\lambda_{mj}}<\infty $. \end{corollary}

\begin{proof} The conditions ensure that the operator with the matrix $ \left(e^{-\lambda_{j\widetilde{j}}}\right) $
is Hilbert--Schmidt, thus compact. \end{proof}

Later, in Section~\ref{s3.55}, we will see that the restriction of
quasi-circularity may be dropped.

\section{Geometry of quasi-algebraic curves }\label{h50}\myLabel{h50}\relax 

\subsection{Quasi-algebraic curves and sheaves } Curve gluing data defines what are
the points of the curve. In this way we obtain a topological space.
Bundle gluing data defines what are fibers over the points, what are
germs of the sections of the bundle, and what are global sections of the
bundle; this defines a sheaf. Call the resulting sheaf on the resulting
topological space a {\em quasi-algebraic sheaf\/} if the gluing data satisfies
the Riemann--Roch theorem.

Taking the gluing functions $ \psi_{j} $ to be 1, we obtain the sheaf $ {\mathcal O} $. Call
the topological space obtained from a curve gluing data a {\em quasi-algebraic
curve\/} if the sheaf $ {\mathcal O} $ satisfies the strong form of the Riemann--Roch
theorem.

\subsection{Ideal points } Let us use the criterion of Section~\ref{s54} to show how far
can be a curve from the classical picture of a complex algebraic curve
without losing the nice properties usually associated with algebraic
curves.\footnote{Of course, the results discussed in this paper do not completely confirm
similarity of these curves with algebraic curves. However, the results of
\cite{Zakh97Qua} suggest that most of others results of algebraic geometry are
applicable to these curves as well: that paper discusses duality for line
bundles (including the case of degree$ \approx $genus), the description of
Jacobians as quotients of vector spaces by lattices, and
conditions on periods of global $ 1 $-forms.} To simplify the exposition, we consider only circular Schottky
curve gluing data with bounded bundle gluing data; we require that these
data satisfy the strong form of the Riemann--Roch theorem; such examples
already exhibit most of the peculiarities of the general theory.

Restrict our attention yet more, to gluing data satisfying
conditions of Corollary~\ref{cor33.40}; for such curves the matching
decomposition and the gluing mappings play no role (as far as they remain
Schottky and bounded), all what matters are conformal distances between
the disks bounding the curve. Consequently, the question boils down to
the following one: how ``bad'' may be a collection of non-intersecting
disks on $ {\mathbb C}{\mathbb P}^{1} $ such that the complement satisfies the strong form of the
Riemann--Roch theorem when equipped with bounded Schottky gluing data.

A natural measure of ``badness'' is how big is the {\em dust\/} (accumulation
points of the disks). We say that a point $ z $ is an {\em accumulation point\/} of a
collection of subsets if any neighborhood of $ z $ intersects with an
infinite number of these subsets.

For example, the first investigated case of quasi-algebraic curves
was the the hyper-elliptic case: the spectral curves for the KdV equation
\cite{McKTru76Hil}. When translated to our language, such a curve corresponds
to the disks accumulating to points $ 0,\infty\in{\mathbb C}{\mathbb P}^{1} $; the disks have centers $ c_{n} $,
$ n\in{\mathbb Z} $; $ c_{n}\approx n^{2} $, $ c_{-n}\approx n^{-2} $ (for $ n\gg 0 $); the radii of the disks are rapidly
decreasing. Many other examples of infinite-genus curves so abundant in
the theory of integrable systems share the same property of having a
finite set as the dust; the conditions on radii of the disks are much
less drastic.

\begin{remark} As we discussed it in the introduction (see also \cite{Zakh97Qua}),
the dust corresponds to {\em ideal points\/} of the curve: the theory of divisors
on the curve should include divisors at the points of dust (though these
points are not included in the smooth part of the curve). \end{remark}

\subsection{Foam curves }\label{s6.3}\myLabel{s6.3}\relax  The natural measure of how far is the curve from the
classical theory is the Hausdorff dimension of the dust: it shows how
``heavy'' are ideal points with respect to ``usual'' points of the curve are.
(Note that any non-dust point of the curve is smooth after the gluing of
the boundary circles together. Dust points are ``very non-smooth''---in any
imaginable way.) Note that for any collection of domains $ D_{j} $, the
accumulation points form a closed subset of $ \coprod_{k}C_{k} $. Moreover, this closed
subset has no interior.

It turns out that the strong Riemann--Roch condition adds no
restrictions on the dust:

\begin{theorem} \label{th41.20}\myLabel{th41.20}\relax  Consider a collection $ \left(C_{k}\right)_{k\in K} $ of compact curves such
that $ g\left(C_{k}\right)\not=0 $ only for a finite number of indices $ k $. For any closed subset
$ {\mathcal D} $ of $ \coprod_{k}C_{k} $ without interior there is curve-bundle gluing data which has $ {\mathcal D} $
as the dust and satisfies the strong form of the Riemann--Roch theorem. \end{theorem}

\begin{proof} To simplify arguments, consider the case of $ K=\left\{1\right\} $ and $ C_{1}=C={\mathbb C}{\mathbb P}^{1} $.
(The key idea needed for the general case is in Remark~\ref{rem41.50}.)
Construct a sequence of non-intersecting closed disks $ R_{j} $ with centers in
$ c_{j} $ and radii $ r_{j} $, $ j\in{\mathbb N} $, one by one.

Choose a sequence $ \left(a_{n}\right) $ of points in $ C\smallsetminus{\mathcal D} $ which has $ {\mathcal D} $ as the set of
accumulation points (possible since $ {\mathcal D} $ has no interior). Start with an
empty collection of disks. Given $ R_{j} $ for $ j<j_{0} $, let $ c_{j_{0}} $ be the first point
of the sequence which is in $ C\smallsetminus\bigcup_{j<j_{0}}R_{j} $. Let $ \widetilde{r}=\operatorname{dist}\left(c_{j_{0}},{\mathcal D}\cup\bigcup_{j<j_{0}}R_{j}\right) $. Let
$ r_{j_{0}} $ be the maximal number $ \leq\widetilde{r}/2 $ satisfying the conditions
$ \sum_{j<j_{0}}e^{-2\lambda_{jj_{0}}}\leq2^{-j_{0}} $; here $ \lambda_{jk} $ is the conformal distance between $ R_{j} $ and $ R_{k} $,
$ \operatorname{ch}\lambda_{jk} = \frac{\left(c_{j}-c_{k}\right)^{2}-r_{j}^{2}-r_{k}^{2}}{2r_{j}r_{k}} $. Such $ r_{j_{0}} $ is unique, since each summand
increases when $ r_{j_{0}} $ increases.

It is clear that $ \bigcup_{j\leq j_{0}}R_{j} $ does not intersect $ {\mathcal D} $, the disks do not
intersect, and have no accumulation points outside $ {\mathcal D} $. Moreover, the
collection of disks satisfies Corollary~\ref{cor33.40}. Now any point $ a_{n} $ is
either a center of a disk from the collection, or is inside one of the
disks $ R_{1},\dots ,R_{n-1} $. Thus any accumulation point of the points $ a_{n} $ is either
in one of the disks $ R_{j} $, or an accumulation point of the disks $ R_{n} $. \end{proof}

In other words, there is no restriction on the dust except that it
is a set of accumulation points of non-intersecting open subsets. In
particular, the dust may have a positive measure. Moreover, there may be
nothing but the dust on the curve:

\begin{theorem} There is a closed subset $ D\subset C={\mathbb C}{\mathbb P}^{1} $ such that
\begin{enumerate}
\item
$ {\mathbb C}{\mathbb P}^{1}\smallsetminus D $ is a union of non-intersecting disks;
\item
$ D $ coincides with the set of accumulation points of these disk;
\item
$ D\subset{\mathbb C}{\mathbb P}^{1} $ satisfies the strong form of the Riemann--Roch theorem when
equipped with an arbitrary bounded Schottky gluing data.
\end{enumerate}
\end{theorem}

\begin{proof} Put $ {\mathcal D}=\varnothing $; enumerate points $ z=x+iy $ with rational $ x $, $ y $; let $ \left(a_{n}\right) $
be the corresponding sequence. Construct a sequence of closed disks in
$ {\mathbb C}{\mathbb P}^{1} $ basing on these $ {\mathcal D} $ and $ \left(a_{n}\right) $ as in the proof of Theorem~\ref{th41.20}.

This sequence of non-intersecting disks contains all the points $ x+iy $
with $ x,y\in{\mathbb Q} $. Let $ R_{m} $ be the interior of the $ m $-th disk. Put $ D $ to be the
complement of $ \bigcup R_{m} $. All we need to show is that any point $ z\in D $ is an
accumulation point of the disks $ R_{m} $. Since the closed disks $ \bar{R}_{m} $ contain all
the rational points $ a_{n} $, there is a sequence of points of $ \bigcup R_{m} $ which goes
to $ z $; if $ z\in\partial R_{m_{0}} $, then we can additionally require that no point of this
sequence is in $ R_{m_{0}} $. It is clear that one can choose a subsequence such
that all its points are in different disks $ R_{m} $. \end{proof}

\begin{remark} \label{rem41.50}\myLabel{rem41.50}\relax  In the case of a curve $ C $ of arbitrary genus we need also
to construct a subset $ \varepsilon $ of positive measure to make it into $ \operatorname{Supp} E $. Use
notations of Theorem~\ref{th41.20}. Construct $ \varepsilon $ recursively: let $ \varepsilon_{0}=C $; choose
$ \varepsilon_{j_{0}} $ as $ \varepsilon_{j_{0}-1}\smallsetminus\widetilde{R}_{j_{0}} $. Here $ \widetilde{R}_{j_{0}} $ is a centered at $ c_{j_{0}} $ disk such that
$ \operatorname{measure}\left(\varepsilon_{j_{0}}\right)/ \operatorname{measure}\left(\varepsilon_{j_{0}-1}\right)>1-2^{-j_{0}} $. Now choose $ R_{j_{0}} $ as a disk in $ \widetilde{R}_{j_{0}} $ of a
very small relative radius. \end{remark}

\begin{remark} There is a more customary description of such ``foam'' curves via
a ``limit process'' involving compact curves of finite genus.\footnote{To make this description into a theorem, we need first to describe which
gluing data describe ``the same'' quasi-algebraic curve, and define a
topology on the resulting moduli space. These topics are outside of the scope
of this paper.} Take $ C_{0}={\mathbb C}{\mathbb P}^{1} $.
Consider a metric on $ C_{0} $, $ \varepsilon_{0}>0 $, and a finite subset such that its
$ \varepsilon_{0} $-neighborhood coincides with $ C_{0} $. Remove from $ C_{0} $ the collection of small
disks with centers in this subset. Break this collection into pairs; for
each pair consider two circles which are boundaries of this pair of
removed disks. Glue in a long cylinder along these two circles; this way
we add a handle per a pair of removed disks. Call the resulting curve $ C_{1} $.

Now repeat the same process for $ C_{1} $, choosing $ \varepsilon_{1}\ll \varepsilon_{0} $. Note that one
needs to remove the disks not only from the common part of $ C_{0} $ and $ C_{1} $,
but also from the added cylinders. Call the resulting curve $ C_{2} $ etc.

One can show that if the lengths of added cylinders increase quickly
enough, then the ``limit'' of the curves $ C_{k} $ is a well-defined
quasi-algebraic curve which enjoys all the benefits of the results of
algebraic geometry. \end{remark}

\begin{remark} As the example above shows, it is possible that a
quasi-algebraic curve has no smooth point at all. Moreover, it may happen
that no point of the curve has a neighborhood homeomorphic to a disk.

Keep in mind that our definition of topology on a quasi-algebraic
curve (and of the set of points of a quasi-algebraic curve) has no deep
underlying reason: they come from the construction of the set of the
functions we allow. This somewhat undermines the observation of the
preceding paragraph. On the other hand, it was shown in \cite{Zakh97Qua} that
by strengthening conditions on the pairwise conformal distances, one may
ensure that there is a well-defined power-series asymptotic formula near
each point of the curve. As asymptotic formulae do, it works outside
some ``sparse'' subset of the neighborhood of a point; the complement of
this subset is locally homeomorphic to a subset of $ {\mathbb C} $. \end{remark}

\begin{remark} One of the incentives to consider quasi-algebraic curves is
that they might describe non-perturbative amplitudes of string
propagation. Recall that the genus $ g $ curves are responsible for $ g $-loops
Feynman diagrams for string propagations; this provides a measure on the
finite-genus moduli spaces. There is a hope that these measures may be
just ``residues'' of a measure on the moduli space of curves of
arbitrary (including infinite) genus.

To investigate the details of this measure one needs much more
advanced knowledge of the theory of quasi-algebraic curves. Even the
question about the {\em support\/} of this measure is not clear at all. Consider
the analogy with the Wiener measure on functions $ f\left(x\right) $: ``the density'' of
this measure $ \exp \left(-\int_{0}^{1}f'\left(x\right)^{2}dx\right) $ is well defined for $ f\in H^{1}\left(\left[0,1\right]\right) $; however,
the support of this measure involves much more singular functions.
Similarly, one cannot a priori expect that ``foam'' curves miss the support
of the above measure. A natural path to investigate this measure is to
consider as general theory as possible;\footnote{This is the reason for why we try to make every statement as general as
possible.} if one of the corollaries of this
theory is that interesting measures do not involve the foam curves, then
one may want to consider simplifications of the theory which work with
much easier ``manageable'' curves only. \end{remark}

\begin{remark} Current approaches to quantum gravity suggest that
consideration of a foam space-time may be necessary. On the other hand,
strings have much less singular distribution of the mass (``along a curve''
instead of ``at a point''), which alleviates the need to drastically change
the topology of space-time. It may be a poetic justice if the
non-perturbative approach to string propagation leads indeed to
consideration of ``foam strings'', thus moving the dust out from under the
carpet into a different place. However, this place promises to be much
more convenient to deal with dust than considerations of ``foam''
space-time. \end{remark}

\begin{remark} Recall why string theory has a chance to be easier to deal with
comparing to the particle theory. Indeed, a particle concentrates its
mass in an arbitrary small region of space. Such a concentration causes a
cloud of virtual particle-antiparticle pairs around this point, as
well as the singularity of the metric (which leads to foam models of
space-time). On the other hand, a string concentrates vanishingly small
mass in small regions of space, thus has a better chance to avoid the
serious singularities.

Judging by the apparent failure of the initial rosy expectations
of the string theory boom, concentrating the mass along space-dimension
one is not enough to avoid singularities. However, if the ``string measure
on the infinite-genus moduli space'' is wild enough to be concentrated on
the set of foam curves, the corresponding strings in space are fractal.
Being fractal, they may have a Hausdorff dimension larger than 1 (or of
$ 1+0 $ type), thus may cause milder singularities than the smooth curves. \end{remark}

\subsection{Why closed operators? }\label{s3.45}\myLabel{s3.45}\relax  Up to this place we applied Theorem~%
\ref{th4.90} only in the case of bounded bundle gluing data, when both the
operators $ A_{1,2} $ are bounded. To deal with this
case one does not need the machinery of Section~\ref{h4}, one could apply the
usual Fredholm theorem. However, simple heuristics show that to have a
good relation between divisors and line bundles the bounded gluing data
is not enough.

Indeed, consider the effect of removing a point $ z_{0} $ from a divisor on
the corresponding bundle gluing data. We want the global sections $ f\left(z\right) $ of
the new bundle correspond to sections of the old bundle which vanish at
this point. Suppose that the curve is glued of a subset $ D\subset{\mathbb C}\subset{\mathbb C}{\mathbb P}^{1} $, then $ f $
is a section of the new bundle if $ \left(z-z_{0}\right)f\left(z\right) $ is a section of the old
bundle. Thus the new bundle gluing data $ \widetilde{\psi} $ may be obtained from the old
data $ \psi $ via\footnote{These gluing data do not satisfy the conditions that all but a finite
number of functions $ \psi_{j} $ are constant. To preserve this condition one needs
to replace $ z-z_{0} $ by a slightly different function (but still vanishing at
$ z_{0} $).
Such a change would not influence the following arguments, so we skip it.} $ \widetilde{\psi}_{j}\left(z\right)=\frac{\varphi_{j}\left(z\right)-z_{0}}{z-z_{0}}\psi_{j}\left(z\right) $.

Suppose that a curve has a very long handle; in the language of
gluing data this corresponds to $ D $ containing an annulus with a very large
conformal distance between two components of the boundary. The robustness
argument from the introduction suggests that moving a point of a divisor
along such a handle has a very little effect on the sections of the
corresponding line bundle (at least at the part of $ D $ ``far'' from this
annulus). On the other hand, the formula above shows that the effect of
such a movement on $ \psi_{j} $ is very large if $ \varphi_{j} $ send a contour ``inside'' the
annulus to one ``outside'' it.

To obtain a ``model'' of the given curve as a subset of $ {\mathbb C}{\mathbb P}^{1} $ with some
gluings of boundaries, one needs to make a meridianal cut across each
``handle'' of the curve. To apply the preceding argument, the handle we
consider should have two marked subsets: a long tube which will become
the annulus, and the loop we made the cut across. They should not
intersect; consequently, two sides of the cut will form two curves, one
inside the annulus, another outside. Moreover, all curves $ \gamma_{j} $ but one are
either all inside or all outside of the annulus on $ {\mathbb C}{\mathbb P}^{1} $. In other words,
there is going to be exactly one gluing function $ \psi_{j} $ which is strongly
affected by the movement of the point inside the annulus.

Thus minor changes to a line bundle may correspond to a giant change
of the gluing functions. For curves which have infinitely many handles
which get progressively longer, one should expect that one can combine
infinitely many such changes (each corresponding to moving a point of the
divisor inside one handle) so that the change to a line bundle is still
negligible, but the gluing functions become unbounded. Thus to have a
nice theory of divisors, one {\em needs\/} to consider unbounded operators $ {\mathcal R} $ too.

Apply the non-bounded case of Theorem~\ref{th4.90} to the settings of
Corollary~\ref{cor33.40}:

\begin{amplification} \label{amp3.450}\myLabel{amp3.450}\relax  If the distortions $ \Delta\left(E_{k,k'},\gamma_{j,j'},\varphi_{j,j'},\psi\equiv 1\right) $ are
uniformly bounded, then in the conditions of Corollary~\ref{cor33.40} one can
replace the boundness condition on the bundle gluing data together with
the condition $ \sum_{m\not=j}e^{-2\lambda_{mj}}<\infty $ by the condition $ \sum_{m\not=j}e^{-2\lambda_{mj}}|\psi_{j'}|^{2}<\infty $. \end{amplification}

\begin{remark} We do not know whether the line bundles which satisfy the
strong form of the Riemann--Roch theorem provide a nice
divisor-to-line-bundle correspondence. This question needs a separate
investigation. \end{remark}

\subsection{Non quasi-circular case }\label{s3.55}\myLabel{s3.55}\relax  Every curve of genus 1 has a circular
Schottky model; one can construct examples of curves of genus 2 which have
no circular Schottky model. On the other hand, any curve of finite genus
has a Schottky model; it is natural to expect that a similar property
should hold for ``reasonable'' curves of infinite genus too. However, it is
not readily obvious that these ``reasonable'' curves would satisfy the
quasi-circularity condition.

Thus it may be vital to be able to drop the condition of
quasi-circularity; it is especially useful since this comes in
essentially no cost. The preceding arguments needed the condition of
quasi-circularity since the operator $ {\mathcal R} $ involves the operator $ \iota_{j}^{\text{skew}} $ of
Section~\ref{s27}. Moreover, this operator is needed since we need to skew the
norm on $ H_{+}^{1/2}\left(\gamma_{j}\right) $ so that Theorem~\ref{th35.157} holds. Recall that the only
role the skew norm plays in the proof of this theorem is the fact that
the image of the direct product of the restriction mappings to an
infinite collection of disjoint subsets is in fact in $ \bigoplus_{l_{2}} $ of their
individual images, as opposed to $ \prod $ (in other words, the ``correlation'' of
the restrictions is negligible).

Avoid the need to consider $ + $-skewed norms by encoding the above fact
into a definition:

\begin{definition} \label{def45.10}\myLabel{def45.10}\relax  We say that a quasi-smooth subset $ D\subset C $ is {\em fat\/} if the
restriction mapping $ H^{1}\left(C\right)/\operatorname{const} \to \prod H^{1}\left(R_{j}\subset C\right)/\operatorname{const} $ has its image in
$ \bigoplus_{l_{2}}H^{1}\left(R_{j}\subset C\right)/\operatorname{const} $. Here $ R_{j} $, $ j\in J $, are connected components of $ C\smallsetminus D $. The
{\em restriction-tolerance\/} of $ D $ is the norm of the corresponding mapping $ H^{1}\left(C\right)
\to \bigoplus_{l_{2}}H^{1}\left(R_{j}\subset C\right)/\operatorname{const} $.

A collection of fat subsets $ D_{k}\subset C_{k} $ is {\em uniformly fat\/} if their
restriction-tolerances are uniformly bounded. \end{definition}

Obviously, if the subsets $ D_{k}\subset C_{k} $ of the gluing data are uniformly
fat, one can consider the (non-skewed) embedding norms on $ H^{1/2}\left(\gamma_{j}\right) $
whenever we considered the skewed norms before. With such an approach
there is no mapping $ \iota_{j}^{\text{skew}} $ in the context of Section~\ref{s27}, thus we may
drop the restriction of quasi-circularity from all the statements as far
as fatness conditions hold.

On the other hand, the condition of fatness is very close to the
conditions of Theorem~\ref{th33.35}. To show this, we start with a technical
statement:

\begin{proposition} \label{prop45.20}\myLabel{prop45.20}\relax  Consider subspaces $ H_{i} $, $ i\in I $, of a Hilbert space $ H $;
denote their orthogonal complements by $ H_{i}^{\perp} $. Then the following conditions
are equivalent:
\begin{enumerate}
\item
the product of projection mappings $ H \to \prod_{i}H/H_{i} $ has it image in
$ \bigoplus_{l_{2},i}H/H_{i} $;
\item
the sum of embedding mappings $ \bigoplus_{i}H_{i}^{\perp} \xrightarrow[]{\iota} H $ extends continuously to a
mapping $ \bigoplus_{l_{2},i}H_{i}^{\perp} \to H $;
\item
the block-Gram matrix $ {\mathcal P} $ (with the blocks $ {\mathcal P}_{ij} $ being the orthogonal
projectors $ H_{j}^{\perp} \to H_{i}^{\perp} $) defines a continuous operator in $ \bigoplus_{l_{2},i}H_{i}^{\perp} $.
\end{enumerate}
\end{proposition}

\begin{proof} The first two statements are dual to each other. The second
implies the third, since $ {\mathcal P} = \iota^{*}\iota $.

To show that the third implies the second, note that $ {\mathcal P} $ is Hermitian,
and non-negative definite. Thus it induces a Hilbert structure on a
completion $ H' $ of a quotient of $ \bigoplus_{l_{2},i}H_{i}^{\perp} $. The Hilbert spaces $ H_{i}^{\perp} $ are
naturally isometrically embedded into $ H' $; it is clear that the natural
mappings $ i $, $ i' $ of $ \bigoplus_{i}H_{i}^{\perp} $ into $ H $ and $ H' $ satisfy $ \|i\left(h\right)\|=\|i'\left(h\right)\| $. This
defines an isometric identification of $ H' $ with a subspace of $ H $; thus a
continuous mapping $ \bigoplus_{l_{2},i}H_{i}^{\perp} \to H $. \end{proof}

\begin{remark} If Hilbert spaces $ H $ and $ \widetilde{H} $ are Hilbert-dual to each other, then
one can consider the subspaces $ H_{i}^{\perp} $ as being subspaces of $ \widetilde{H} $. Since
$ H^{s}\left(R\subset C\right)=H^{s}\left(C\right)/\mathring{H}^{s}\left(C\smallsetminus R\right) $, the natural duality between $ H^{s}\left(C\right) $ and $ H^{-s}\left(C,\omega\otimes\bar{\omega}\right) $
makes the calculation of the orthogonal complement to $ \mathring{H}^{s}\left(C\smallsetminus R\right) $ especially
simple: it is $ \mathring{H}^{-s}\left(R,\omega\otimes\bar{\omega}\right)\subset H^{-s}\left(C,\omega\otimes\bar{\omega}\right) $. Thus to check fatness it is enough
to consider the block-Gram matrix formed by the orthogonal projectors
$ \mathring{H}_{\int=0}^{-1}\left(R_{i},\omega\otimes\bar{\omega}\right) \to \mathring{H}_{\int=0}^{-1}\left(R_{j},\omega\otimes\bar{\omega}\right) $; here $ H_{\int=0}^{s}\subset H^{s} $ consists of sections with
the integral being 0.

Obviously, $ H_{\int=0}^{-1}\left(C,\omega\otimes\bar{\omega}\right) $ has a natural Hilbert norm, since it is
dual to $ H^{1}\left(C\right)/\operatorname{const} $. This norm is invariant with respect to automorphisms
of $ C $. \end{remark}

Call the orthogonal projector $ \mathring{H}_{\int=0}^{-1}\left(R,\omega\otimes\bar{\omega}\right) \to \mathring{H}_{\int=0}^{-1}\left(R',\omega\otimes\bar{\omega}\right) $ the
$ \Delta $-{\em projector\/} between $ R $ and $ R' $. Consider a quasi-smooth subset $ D\subset C $ of a
compact curve $ C $ with connected components $ D_{j} $, $ j\in J $, of the complement.
Associate to it the infinite block-matrix $ \left({\mathcal P}_{jk}\right) $, $ j,k\in J $, with block $ {\mathcal P}_{jk} $
being the $ \Delta $-projector between $ R_{k} $ and $ R_{j} $. One can consider $ \left({\mathcal P}_{jk}\right) $ as a
matrix of an operator $ \widetilde{{\mathcal P}}\colon \bigoplus_{j}\mathring{H}_{\int=0}^{-1}\left(R_{j},\omega\otimes\bar{\omega}\right) \to \prod_{j}\mathring{H}_{\int=0}^{-1}\left(R_{j},\omega\otimes\bar{\omega}\right) $.

\begin{lemma} \label{lm45.40}\myLabel{lm45.40}\relax  If $ D $ is fat, then the operator $ {\mathcal C} $ of Section~\ref{s352} (with
the image in $ + $-skewed spaces $ H_{+}^{1/2}\left(\partial_{j}\right) $) is bounded. \end{lemma}

\begin{proof} Let $ H_{j}=\mathring{H}_{\int=0}^{-1}\left(R_{j},\omega\otimes\bar{\omega}\right) $, $ H=\bigoplus_{l_{2},j}H_{j} $. Let $ H_{j+}=\operatorname{Im}\bar{\partial}|_{\mathring{H}^{0}\left(R_{j},\omega\right)} $,
$ H_{j-}=\operatorname{Im}\partial|_{\mathring{H}^{0}\left(R_{j},\bar{\omega}\right)} $. First, the component of $ {\mathcal P}_{jk} $ acting from $ H_{k-} $ to $ H_{j-} $
vanishes, similarly for the component acting from $ H_{k+} $ to $ H_{j+} $.

Indeed, consider the operator $ \Delta=i\partial\bar{\partial}\colon H^{1}\left(C\right)/\operatorname{const} \to H_{\int=0}^{-1}\left(C,\omega\otimes\bar{\omega}\right) $.
Being an elliptic self-dual operator with no kernel, it is invertible.
The norm in $ H_{\int=0}^{-1}\left(C,\omega\otimes\bar{\omega}\right) $ is defined via the duality with $ H^{1}\left(C\right)/\operatorname{const} $;
$ \|\alpha\|_{H^{-1}}=\max _{f}\frac{|\int\bar{f}\alpha|}{\|f\|_{H^{1}}^{1/2}} $; here $ f\in H^{1}\left(C\right)/\operatorname{const} $. Let $ g=\Delta^{-1}\alpha $; then
$ \int\bar{f}\alpha=-i\int\bar{f}\bar{\partial}\partial g=i\int\bar{\partial}\bar{f}\partial g=i\int\overline{\partial f}\partial g $. Since $ \|\partial f\|_{L_{2}}=\|f\|_{H^{1}/\operatorname{const}} $, we can see that
$ \|\alpha\|_{H^{-1}}=\|\partial g\|_{L_{2}}=\|g\|_{H^{1}/\operatorname{const}} $. Thus $ \Delta $ is a unitary operator, and
$ \left(\alpha_{1},\alpha_{2}\right)_{H^{-1}}=\left(\bar{\partial}\Delta^{-1}\alpha_{1},\bar{\partial}\Delta^{-1}\alpha_{2}\right)_{L_{2}}= \left(\partial^{-1}\alpha_{1},\partial^{-1}\alpha_{2}\right)_{L_{2}} $ (here we choose $ \partial^{-1}:
H_{\int=0}^{-1}\left(C,\omega\otimes\bar{\omega}\right) \to H^{0}\left(C,\bar{\omega}\right) $ so that the image is orthogonal to global
holomorphic forms). This shows that for $ g\left(C\right)=0 $ the component of
$ \Delta $-projector acting from $ H_{k-} $ to $ H_{j-} $ vanishes; similarly for $ + $-components.
For general $ g\left(C\right) $ these components do not vanish; the corresponding norm
on $ \bigoplus_{l_{2}}H_{k-} $ is not the direct sum norm, but its correction via the
projection on the orthogonal complement to the image of $ \partial $. Thus this norm
differs from the direct sum norm by a continuous finite-rank operator.

Similarly, if $ \alpha_{1}=\bar{\partial}\beta_{1} $, $ \alpha_{2}=\partial\beta_{2} $, and $ g\left(C\right)=0 $, then $ \left(\alpha_{1},\alpha_{2}\right)_{H^{-1}} $ is
proportional to $ \int\bar{\beta}_{1}\partial\bar{\partial}^{-1}\beta_{2} $. Thus the identifications above identify the
$ \Delta $-projector from $ H_{k+} $ to $ H_{j-} $ with the bar-projector. Again, for general $ C $
the projector from $ \bigoplus H_{k+} $ to $ \bigoplus H_{k-} $ differs from the operator $ {\mathcal C} $ by a
continuous finite-rank operator. \end{proof}

\begin{remark} If $ {\mathcal C} $ is bounded, then all the the components of $ {\mathcal P} $ between
$ \bigoplus H_{k+} $, $ \bigoplus H_{k-} $ are bounded too. Moreover, $ H_{k+}+H_{k-}=H_{k} $. Indeed, the
orthogonal complements to $ H_{k+} $ and $ H_{k-} $ coincide with $ \operatorname{Ker}\bar{\partial} $ and $ \operatorname{Ker}\partial $ in
$ H^{1}\left(R_{k}\subset C\right)/\operatorname{const} $; these subspaces do not intersect.

However, this does not yet imply that $ {\mathcal P} $ is bounded. Indeed, the
angle between $ H_{k+} $ and $ H_{k-} $ can be arbitrary small (as examples of long
ellipses show). However, in the context of ``practical criteria'' of
Section~\ref{s54}, $ {\mathcal C} $ and $ {\mathcal P} $ behave the same way: \end{remark}

\begin{proposition} Consider two regions $ R,R'\subset C $ of conformal distance $ l $.
\begin{enumerate}
\item
There is a constant $ c $ such that if $ C={\mathbb C}{\mathbb P}^{1} $ and $ R,R' $ are disks, then the
norm of the $ \Delta $-projector between $ R $ and $ R' $ is $ c e^{-l} $.
\item
There is a constant $ c' $ such that the norm of the $ \Delta $-projector between
$ R $ and $ R' $ is less than $ c' e^{-l} $.
\end{enumerate}
\end{proposition}

\begin{proof} The first part follows from Lemma~\ref{lm31.50} and the proof
of Lemma~\ref{lm45.40}; indeed, for circles the subspaces $ H_{k+} $ and $ H_{k-} $ of this
proof are orthogonal.

The second part of the statement can be proved similarly to
Proposition~\ref{prop31.60}. \end{proof}

\begin{corollary} In Theorem~\ref{th33.35} and in Corollary~\ref{cor33.40} one can drop
the condition of quasi-circularity. \end{corollary}

\begin{remark} One can combine this corollary and Amplification~\ref{amp3.450}.
However, to ensure that conditions of Theorem~\ref{th4.90} hold, one needs to
require that both $ \sum_{m\not=j}e^{-2\lambda_{mj}}|\psi_{j'}|^{2}<\infty $ and $ \sum_{m\not=j}e^{-2\lambda_{mj}}<\infty $ are finite. The
first part gives compactness of $ {\mathcal R}\circ{\mathcal C} $, the second boundness of $ {\mathcal C} $. \end{remark}

\subsection{Non-pseudo-smooth case }\label{s6.6}\myLabel{s6.6}\relax  Recall that so far we assumed that the
connected components $ R_{i} $ of the complement to the model domain $ D\subset C $ is a
union of regions with smooth boundary. However, the only place when this
assumption is crucial is Lemma~\ref{lm131.40}; in turn, this lemma is needed
for the introduction of $ + $-skewed norms.

As the previous section shows, the consideration of $ + $-skewed norms
can be avoided in many cases. Consequently, in these cases one can weaken
the assumptions on the components $ R_{i} $. For example, one can assume that
the curves $ \gamma_{i}=\partial R_{i} $ are Jordan curves (what is crucial, due to Lemma~%
\ref{lm131.44}, is that $ \partial R_{i} $ has measure 0). In such a case one considers
Lemma~\ref{lm131.32} as the {\em definition\/} of $ H^{1/2}\left(\gamma\right) $.

However, we need the gluing functions $ \varphi_{i} $ and $ \psi_{i} $ induce a mapping
between $ H^{1/2}\left(\gamma_{i}\right) $ and $ H^{1/2}\left(\gamma_{i'}\right) $. The simplest modification to allow this
is to require $ \varphi_{i} $, $ \psi_{i} $ be defined not on $ \gamma_{i} $, but on a neighborhood of $ \gamma_{i} $,
and require that the corresponding gluing operators act in $ H^{1} $ of these
neighborhoods.\footnote{Recall that the arguments of \cite{CourHurv} suggest that any reasonable
curve should have a model with $ \varphi_{\bullet} $ being fraction-linear. As shown in
\cite{Zakh97Qua}, line bundles should be representable by constant functions $ \psi_{\bullet} $
(with an exception of one pair of curves $ \gamma $, $ \gamma' $ to allow non-zero degree).} (For the latter condition, it is enough to require that
these functions are Lipschitz.)

\begin{remark} This modification is very welcome, since as Section~\ref{s6.3}
shows, our approach allows consideration of gluing data of fractal
nature. It does not make a lot of sense to allow $ \partial D $ to be fractal, while
requiring that all the connected components of $ C\smallsetminus D $ have smooth
boundaries.

Moreover, the typical description of a complex curve via the
Schottky model leads to an invariantly defined {\em Schottky group\/}; our
language requires a choice of a fundamental domain for this group. One
should hope that weakening the requirements on $ \partial R_{i} $ allows treating most
of the ``natural'' Schottky groups using our language. \end{remark}

\subsection{Black-white curves }\label{s6.70}\myLabel{s6.70}\relax  The abstract form of the Riemann--Roch
theorem~\ref{th4.90} we use to establish the Riemann--Roch theorem needs the
condition that 1 is not in the essential spectrum of an appropriate
operator $ B $ (which is $ {\mathcal C}\circ{\mathcal R} $ or $ {\mathcal R}\circ{\mathcal C} $). However, up to now we used a much
weakened form where the operator $ B $ is assumed to be compact. By analogy
with Remark~\ref{rem6.60}, one could expect that the compactness condition
should be close to the {\em strict\/} Riemann--Roch condition. However, the very
restrictive form of the operator $ {\mathcal R} $ allows curves to have a non-compact
operator $ {\mathcal C} $ such that the Riemann--Roch theorem still holds.

Indeed, if $ {\mathcal R} $ is bounded (so $ {\mathcal C}\circ{\mathcal R} $ is bounded), and $ \left({\mathcal C}\circ{\mathcal R}\right)^{2} $ is compact,
then $ {\mathcal C}\circ{\mathcal R} $ has only 0 in the essential spectrum. This may be achieved if $ {\mathcal C} $
can be made $ 2\times2 $ block-diagonal with one block being compact, another
bounded, and $ {\mathcal R} $ interchanging these blocks.

The following example of such curves is similar to the curves
studied in \cite{GiesKnoTru93Geo,FelKnoTru96Inf}, and \cite{Mer00Rie}; this is
why we discuss it in details. A curve gluing data defines a {\em black-white\/}
curve, if the set of indices $ K $ enumerating compact curves $ C_{k} $ is broken
into two parts, $ K_{b} $ and $ K_{w} $, and no boundary components of the $ w $-parts are
glued together. In other words, the glued curve is colored in two colors,
and the connected components of the white part consist of one domain $ D_{k} $,
$ k\in K_{w} $.

It is easy to see that for such curves to satisfy the strong form of
Riemann--Roch theorem with bounded bundle gluing data, it is enough to
require that blocks $ {\mathcal C}_{jl} $ corresponding to {\em black curves only\/} form a compact
matrix (with our ``usual assumptions'' about white curves, so that the
corresponding white blocks $ {\mathcal C}_{jl} $ are bounded; one can also drop the ``usual
assumptions'' on the white blocks if one simultaneously replaces the
conditions on the black curves by stronger restrictions than just
compactness). In other words, there is no restriction on the white
domains $ D_{k}\subset C_{k} $, $ k\in K_{w} $ (except the usual restrictions on boundness of
distortions of excess spaces).

In particular, given an arbitrary curve-bundle gluing data
$ \left(D_{k}\subset C_{k},',\varphi_{\bullet},\psi_{\bullet},V_{\bullet}\right) $, $ k\in K $, one can convert it to a black-white gluing data.
To do so, let $ K_{w}=K $, $ K_{b} $ consists of $ ' $-orbits. Let $ D_{k} $ for $ k\in K_{b} $ be an
annulus $ \left\{1<|z|<N_{k}\right\} $ with $ N_{k}\gg 1 $ appropriately embedded into $ C_{k}\simeq{\mathbb C}{\mathbb P}^{1} $. Glue
this wide annulus (same as a long tube) $ D_{k} $ between the corresponding
curves $ \gamma_{j} $, $ \gamma_{j'} $. It is easy to see that the new gluing data is a
black-white curve, and that with an appropriate choice of $ N_{k} $, $ k\in K_{b} $, this
gluing data satisfies the Riemann--Roch theorem.

In fact the only thing to prove is that one can embed $ D_{k} $ in $ {\mathbb C}{\mathbb P}^{1} $ and
glue $ \partial D_{k} $, $ k\in K_{b} $, to $ \gamma_{j} $, $ \gamma_{j'} $ in such a way that the distortion
$ \Delta\left(E_{k,k'}\gamma_{j,j'},\varphi_{j,j'},\psi\equiv 1\right) $ is bounded. The simplest way to do this is the
following one: by definition, the real curve $ \gamma_{j} $ is a boundary of the
region $ R_{j} $ inside the corresponding complex curve $ C_{k_{j}} $. Choose $ R_{j}'\subset R_{j} $ such
that the conformal distance between $ \partial R_{j}' $ and $ \partial R_{j} $ is $ N_{k} $; proceed similarly
for $ \gamma_{j'} $. Now identify the annuli $ R_{j}\smallsetminus R_{j}' $ and $ R_{j'}\smallsetminus R_{j'}' $ so that $ \partial R_{j} $ is
identified with $ \partial R'_{j'} $; gluing $ R_{j} $ and $ R_{j'} $ using this identification gives
a rational curve $ C_{k} $ with two contractible regions $ R_{j}' $, $ R_{j'}' $ with
conformal distance $ N_{j} $. Moreover, the boundary of $ D_{k}\buildrel{\text{def}}\over{=}R_{j}\smallsetminus R_{j}'=R_{j'}\smallsetminus R_{j'}' $ is
naturally glued to $ \gamma_{j} $ and $ \gamma_{j'} $. Now the check that the distortion is
bounded is tautological; thus the distortion
$ \Delta\left(E_{k_{j}},0_{'}\gamma_{j},\partial R_{j},\varphi_{j,j'},\psi\equiv 1\right)=1 $.

To take into account $ \psi_{j,j'} $, it is enough to consider the case of
constant $ \psi_{j} $. Make $ \psi_{\bullet} $ corresponding to the gluing of $ \gamma_{j} $ and $ \partial R_{j} $ to be 1,
and $ \psi_{\bullet} $ corresponding to the gluing of $ \gamma_{j'} $ and $ \partial R_{j'} $ to be $ \psi_{j} $; this makes
the corresponding $ \Delta $ to become $ |\psi_{j}| $. Note that we do not need boundness of
$ \left\{\psi_{j}\right\} $, an increase of $ \psi_{j} $ can be compensated by an increase of the
corresponding $ N_{k\left(j\right)} $.

\subsection{The bundle $ {\protect \mathcal O} $ }\label{s3.60}\myLabel{s3.60}\relax  When all the functions $ \psi_{j} $, $ j\in J $, are 1, the bundle
gluing data describes the bundle $ {\mathcal O} $; one can take the allowance spaces $ V_{j} $
to be spanned by 1. In this important special case the mismatch operator
$ \mu $ of Section~\ref{s10.80} not only has index 0, but also is an isomorphism:

\begin{theorem} \label{th8.311}\myLabel{th8.311}\relax  Suppose that $ g\left(C_{k}\right)=0 $, $ k\in K $, for all the curves $ C_{k} $ of the
curve gluing data. Suppose also that the result of gluing is connected.
If the bundle $ {\mathcal O} $ over this curve gluing data satisfies the Riemann--Roch
theorem, then the mismatch operator is a bijection. \end{theorem}

\begin{proof} Since all the excess spaces vanish, the Riemann--Roch theorem
states that $ \mu $ has index 0. Thus it is enough to show that $ \operatorname{Ker}\mu=\left\{0\right\} $, or to
show that for a function $ f $ modulo $ \operatorname{const} $ with $ \bar{\partial}f=0 $ and the vanishing
mismatch one has $ \|\partial f\|_{L_{2}}=0 $. In the case of domains $ D_{k}\subset C_{k} $ with smooth
boundary and smooth $ f $, the following argument works:
$ \int_{D}\partial f\wedge\bar{\partial}\bar{f}=\int_{D}\partial f\wedge d\bar{f}=\int_{D}df\wedge d\bar{f}=\int_{D}d\left(f\,d\bar{f}\right)=\int_{\partial D}f\,d\bar{f} $; here $ D=\coprod_{k}D_{k} $. Now $ \partial D $ is broken
into pairs $ \gamma_{j} $, $ \gamma_{j'} $ with identifications between them, and the pull-backs
of $ d\bar{f} $ to these pairs are compatible with the identifications. Thus
$ \int_{\gamma_{j}\cup\gamma_{j'}}f\,d\bar{f}=\int_{\gamma_{j}}\Delta_{j}f\,d\bar{f}=\Delta_{j}f\int_{\gamma_{j}}d\bar{f}=0 $; here $ \Delta_{j}f $ is the (constant) jump of $ f $
when $ \gamma_{j} $ is identified with $ \gamma_{j'} $.

In our, more general, situation $ f $ may be extended to become an $ H_{\text{loc}}^{1} $
function on $ C=\coprod_{k}C_{k} $, so $ df,d\bar{f} $ are $ L_{2}=H^{0} $-section of $ \Omega^{1}C $. In particular,
$ df\wedge d\bar{f} $ is a well-defined $ L_{1} $-section of $ \Omega^{2}C $. Thus $ \int_{D}df\wedge d\bar{f} $ may be written as
$ \sum_{k}\int_{C_{k}}df\wedge d\bar{f}-\sum_{j}\int_{R_{j}}df\wedge d\bar{f}. $ What we achieved so far is to reduce the question
to integration over compact smooth manifolds with boundaries, as above.
However, the differential forms we need to consider are not smooth.

Given a surface $ S $, a smooth curve $ \gamma\subset S $, and $ H^{1} $-functions $ f $ and
$ g $ on $ C $ {\em define\/} $ \int_{\gamma}f\,dg $ as the result of the pairing of $ f|_{\gamma}\in H^{1/2}\left(\gamma,{\mathcal O}\right) $ and
$ dg\in H^{-1/2}\left(\gamma,\Omega^{1}\gamma\right) $. Now for a subset $ T\subset S $ with a smooth boundary the
expressions $ \int_{T}df\wedge dg $ and $ \int_{\partial T}f\,dg $ are both well-defined and continuous in $ f $
(or $ g $). Since these expressions coincide for smooth $ f $ and $ g $, they
coincide everywhere; this covers the case when the curves $ \gamma_{j} $ are smooth.

In general, when the curves $ \gamma_{j} $ may be Jordan curves, we need an
extra argument. It is easy to reduce what we need to the following
statement:

\begin{lemma} Consider a surface $ S $, a domain $ D\subset S $ such that $ \Gamma=\partial D $ is of measure
0, and a sequence of smooth embedded curves $ \left(\Gamma_{n}\right) $ in $ S $ such that $ \Gamma_{n} $ bounds
a domain $ D_{n} $, and for any neighborhood $ U $ of $ \Gamma $ the symmetric difference of
$ D $ and $ D_{n} $ is inside $ U $ for large enough $ n $. Then bilinear functionals $ A_{n}:
\left(f,g\right) \to \int_{\Gamma_{n}}f\,dg $ on $ H^{1}\left(S\right) $ have a limit $ A $ when $ n \to \infty $. Moreover, $ A\left(f,g\right)=0 $
if the image of $ f $ in $ H^{1}\left(\Gamma\subset S\right) $ vanishes. \end{lemma}

\begin{proof} Since the area between $ \Gamma $ and $ \Gamma_{n} $ goes to 0, the arguments above
imply that the sequence $ A_{n}\left(f,g\right) $ is fundamental for any fixed pair $ \left(f,g\right) $,
thus has a limit. Similarly, the limit is continuous. By definition, if
the image of $ f $ in $ H^{1}\left(\Gamma\subset S\right) $ vanishes, one can approximate $ f $ by a function $ \widetilde{f} $
which vanish near $ \Gamma $; then $ A_{n}\left(\widetilde{f},g\right)=0 $ for large $ n $, so $ A\left(\widetilde{f},g\right)=0 $; by
continuity, $ A\left(f,g\right)=0 $. \end{proof}

This finishes the proof of the theorem. \end{proof}

This theorem is the central statement in the description of the
geometry of the Jacobian of the curve \cite{Zakh97Qua}.

\section{Appendix: Fredholm theorem }\label{h4}\myLabel{h4}\relax 

In this section we discuss details and motivations for Theorem~%
\ref{th4.90}.

\subsection{Quasi-complementary subspaces } Continue using notations of Section~%
\ref{s35.24}. Note that if $ V_{1} $, $ V_{2} $ are quasi-complementary, then the natural
mapping $ V_{1} \to V/V_{2} $ is a Fredholm mapping with the index being the excess
of $ V_{1} $, $ V_{2} $. Recall that arguments of Section~\ref{s4.11} required
finite-dimensional ``corrections'' to graphs of linear operators. First,
show that relative dimension is invariant w.r.t.~finite-dimensional
variations.

\begin{proposition} Suppose that $ V_{1} $, $ V_{2} $ are quasi-complementary with the
excess $ d $, and $ V_{k}' $ is comparable with $ V_{k} $ with the relative dimension $ d_{k} $,
$ k=1,2 $. Then $ V'_{1} $, $ V'_{2} $ are quasi-complementary with the excess $ d+d_{1}+d_{2} $. \end{proposition}

\begin{proof} It is enough to consider the case when $ V_{2}=V_{2}' $, and $ V_{1}'\supset V_{1} $ with
codimension 1. Let $ v\in V_{1}'\smallsetminus V_{1} $.

If $ v\in V_{1}+V_{2} $, we need to show that the codimension of $ V_{1}\cap V_{2} $ in $ V_{1}'\cap V_{2} $
is 1. Let $ v=v_{1}+v_{2} $, $ v_{k}\in V_{k} $, $ k=1,2 $. Then $ v_{2}\notin V_{1}\cap V_{2} $, but $ v_{2}=v-v_{1}\in V_{1}'\cap V_{2} $, thus
$ \operatorname{codim}\left(V_{1}\cap V_{2}\subset V_{1}'\cap V_{2}\right)\geq1 $. If $ w\in V_{1}'\cap V_{2} $, then $ w-\tau v\in V_{1} $, thus $ w-\tau v_{2}\in V_{1}\cap V_{2} $, thus
$ \operatorname{codim}\left(V_{1}\cap V_{2}\subset V_{1}'\cap V_{2}\right)\leq1 $.

Similarly, if $ v\notin V_{1}+V_{2}=\overline{V_{1}+V_{2}} $, then $ V_{1}'\cap V_{2}=V_{1}\cap V_{2} $, $ V_{1}'+V_{2}=V_{1}+V_{2}+{\mathbb C}v $,
thus is closed, and $ \operatorname{codim}\left(V_{1}+V_{2}\subset V_{1}'+V_{2}\right)=1 $. \end{proof}

\begin{remark} Call subspaces $ V_{1} $, $ V_{2} $ {\em weakly quasi-complementary\/} if they are
closed, $ \dim  V_{1}\cap V_{2}<\infty $, and $ \operatorname{codim}\left(\overline{V_{1}+V_{2}}\right)<\infty $; similarly, define
$ \operatorname{reldim}\left(V_{1},V_{2}\right) $. The law of a change to comparable subspaces fails
spectacularly for weakly quasi-complementary subspaces. Indeed, if in the
proof above $ v\notin V_{1}+V_{2} $, but $ v\in\overline{V_{1}+V_{2}} $, then the excesses of $ \left(V_{1},V_{2}\right) $ and
$ \left(V_{1}',V_{2}\right) $ coincide. Thus the weak quasi-complementarity is preserved by
changing spaces to comparable, but there is no way to control the excess.
This is why in what follows we are not interested in weak
quasi-complementarity. \end{remark}

\subsection{Closed operators }\label{s40.22}\myLabel{s40.22}\relax  As explained in Section~\ref{s3.45}, to get a
satisfactory divisor-bundle correspondence, one needs to allow unbounded
bundle gluing functions. This would lead to unbounded operator $ {\mathcal R} $. Recall
the settings of closed operators (see \cite{RudFun} for details).

A {\em partial operator\/} from $ V $ to $ V' $ is a linear operator $ V_{1} \to V' $ with
$ V_{1}\subset V $. Such a partial operator $ A $ is {\em densely defined\/} if $ V_{1} $ is dense in $ V $,
and it has a {\em closed graph\/} if the graph of $ A $ (which is a vector subspace
of $ V_{1}\oplus V'\subset V\oplus V' $) is closed in $ V\oplus V' $. A partial operator is {\em closed\/} if it is
densely defined and it has a closed graph. Obviously, any continuous
operator $ V \to V' $ is closed. Such closed operators are called {\em bounded}.

Two partial operators $ A\colon V_{1} \to W $ (with $ V_{1}\subset V $) and $ B\colon W_{1} \to V^{*} $ (with
$ W_{1}\subset W^{*} $) are {\em in duality\/} if $ w\in W_{1} $ iff $ \left< Av,w \right>=\left< v,v' \right> $ for an
appropriate $ v'\in V^{*} $, and $ Bw=v' $. If $ A $ is densely defined, there is a unique
operator which is in duality with $ A $, it is called the {\em dual operator}. The
dual operator automatically has a closed graph. If $ A $ is closed, then the
dual operator $ A^{*} $ is closed, $ A^{**}=A $, and the graph of $ -A^{*} $ (which is a
vector subspace in $ W^{*}\oplus V^{*}\simeq V^{*}\oplus W^{*} $) is the orthogonal complement to the graph
of $ A $ in $ V\oplus W $.

A {\em composition\/} $ A'\circ A $ of two partial operators $ A\colon V_{1} \to V' $ (with $ V_{1}\subset V $)
and $ A'\colon V'_{1} \to V'' $ (with $ V_{1}'\subset V' $) is defined on $ A^{-1}\left(V'_{1}\right) $ as
$ A'\circ A|_{A^{-1}\left(V_{1}'\right)} $. The composition of two closed operators is not necessarily
densely defined, and the closure of the graph of the composition is not
necessarily a graph of a partial mapping (in other words, this closure may
intersect $ 0\oplus V''\subset V\oplus V'' $). However, if $ A' $ is bounded, then $ A'\circ A $ is densely
defined; if $ A $ is bounded, then the closure of the graph of $ A'\circ A $ is a
graph of a partial operator.

If the composition $ A'\circ A $ is densely defined, then $ \left(A'\circ A\right)^{*} $ is defined
on $ v\in V''{}^{*} $ if $ A^{*}\circ A'{}^{*} $ is defined on $ v $, and $ \left(A'\circ A\right)^{*}v=\left(A^{*}\circ A'{}^{*}\right)v $. Note that if
$ A' $ is not bounded, the domain of $ A^{*}\circ A'{}^{*} $ may be strictly smaller than the
domain of $ \left(A'\circ A\right)^{*} $.

We say that two partial operators $ A' $ and $ A $ have {\em a compact
composition\/} $ A'\circ A $ if the graph of $ A'\circ A $ is a vector subspace of a graph
of a compact operator $ B\colon V \to V'' $. Similarly, say that 1 is not in significant
spectrum of $ A'\circ A $ if there is a bounded ``extension'' $ B\colon V \to V'' $ with 1
not in the $ \operatorname{SpecEss}\left(B\right) $.

\subsection{Complementarity criterion } By Theorem~\ref{th52.40}, to get a sufficiently
general version of Riemann--Roch theorem, one should be able to
characterize a sufficiently large subset of the set of
quasi-complementary pairs. Consider two closed vector subspaces $ V_{1},V_{2}\subset H $
such that the projection of $ V_{i} $ on $ H_{i} $ has no null-space and a dense image.
This means that one can consider $ V_{1} $ as a graph of a closed mapping $ A_{1}\colon H_{1}
\to H_{2} $, similarly $ V_{2} $ is a graph of a closed mapping $ A_{2}\colon H_{2} \to H_{1} $.

\begin{lemma}[abstract finiteness] If $ A_{1} $ and $ A_{2} $ have a compact composition
$ A_{1}\circ A_{2} $, then $ V_{1}\cap V_{2} $ is finite dimensional. If $ A_{2}^{*} $ and $ A_{1}^{*} $ have a compact
composition $ A_{2}^{*}\circ A_{1}^{*} $, then $ \overline{V_{1}+V_{2}} $ has a finite codimension. \end{lemma}

\begin{proof} The projection of $ V_{1}\cap V_{2} $ to $ H_{2} $ is a subspace of $ \operatorname{Ker}\left(A_{1}\circ A_{2}-\boldsymbol1\right) $,
thus is finite-dimensional. The second part can be proven by taking
orthogonal complements to $ V_{1} $ and $ V_{2} $. \end{proof}

\begin{corollary}[weak form] If $ A_{1} $ and $ A_{2} $ have a compact densely defined
composition $ A_{1}\circ A_{2} $, then $ V_{1} $ and $ V_{2} $ are weakly quasi-complementary. \end{corollary}

\begin{proposition}[strong bounded form] If $ A_{1} $, $ A_{2} $ are bounded, and $ A_{1}\circ A_{2} $ is
compact, then $ V_{1} $ and $ V_{2} $ are quasi-complementary with the excess 0. \end{proposition}

\begin{proof} It is enough to prove is that $ V_{1}+V_{2} $ is closed.

Let $ v_{l}+v'_{l} \xrightarrow[l\to\infty]{} v\in H_{1}\oplus H_{2} $, here $ v_{l}\in V_{1} $, $ v_{l}'\in V_{2} $, $ l\in{\mathbb N} $. Then
$ v_{l}=\left(h_{l},A_{1}h_{l}\right) $, $ v'_{l}=\left(A_{2}h'_{l},h'_{l}\right) $. Let $ v=\left(x,x'\right) $, $ x\in H_{1} $, $ x'\in H_{2} $. Then
$ h_{l}+A_{2}h'_{l} \to x $, $ h'_{l}+A_{1}h_{l} \to x' $. We need to show that
\begin{equation}
h+A_{2}h' = x,\qquad h'+A_{1}h = x'
\label{equ138.20}\end{equation}\myLabel{equ138.20,}\relax 
has a solution.

One may substitute $ x-A_{2}h'_{l} $ instead of $ h_{l} $, thus it is enough to
show that if $ \left(1-A_{1}A_{2}\right)h'_{l} \to x' $, then $ x'=\left(1-A_{1}A_{2}\right)h' $ for an appropriate
$ h' $, which follows from closeness of $ \operatorname{Im}\left(1-A_{1}A_{2}\right) $. \end{proof}

\begin{remark} Similarly, one may require that $ A_{2}\circ A_{1} $ is compact. This is not
equivalent to $ A_{1}\circ A_{2} $ being compact. Indeed, let $ K $ be a compact operator.
Consider block matrices $ A_{1}= \left(
\begin{matrix}
\boldsymbol1 & 0 \\ 0 & K
\end{matrix}
\right) $, $ A_{2}=\left(
\begin{matrix}
0 & \boldsymbol1 \\ 0 & 0
\end{matrix}
\right) $. \end{remark}

\begin{amplification} One can require instead that 1 is not in the essential
spectrum of $ A_{1}\circ A_{2} $ or $ A_{2}\circ A_{1} $. \end{amplification}

The proposition is an insignificant generalization of Fredholm
theorem, however, since our target is to get as large a class of pairs of
closed subspaces as possible, we need a much stronger result which allows
$ A_{1} $ to be closed instead of bounded. In applications $ A_{2} $ is bounded, we
have a possibility to control $ A_{2} $, and can make it arbitrarily small.
Investigate which conditions on $ A_{1} $ imply quasi-complementarity.

First, consider the case when $ A_{2}=0 $. Then $ A_{1}A_{2} $ is automatically
defined everywhere and compact. However, $ V_{1} $ and $ V_{2} $ being
quasi-complementary is {\em equivalent\/} to $ A_{1} $ being bounded! Thus there is no
hope to generalize the above proposition literally. However, this example
shows that to have $ V_{1}+V_{2} $ closed it is enough to strengthen the topology
on $ H_{2} $, so that $ A_{1} $ becomes bounded. Let $ \|\bullet\|_{A_{1}} $ be the norm on the domain of
$ A_{1} $ defined by identification of it with $ V_{1} $, $ \|v\|_{A_{1}}^{2}=\|v\|_{H_{1}}^{2}+\|A_{1}v\|_{H_{2}}^{2} $, let
$ H_{1}^{\left(A_{1}\right)} $ be the domain of $ A_{1} $ considered with this norm. Obviously, $ H_{1}^{\left(A_{1}\right)} $
is complete, and $ A_{1} $ can be extended to a bounded operator $ H_{1}^{\left(A_{1}\right)} \to H_{2} $,
and $ V_{1} $ is closed in $ H_{1}^{\left(A_{1}\right)}\oplus H_{2} $. Moreover, $ \left(V_{1},H_{2}\right) $ are quasi-complementary
in $ H_{1}^{\left(A_{1}\right)}\oplus H_{2} $ with the excess 0.

How to generalize this to the case $ A_{2}\not=0? $ It is clear that it is
enough to require that $ \operatorname{Im} A_{2}\subset H_{1}^{\left(A_{1}\right)} $, or that $ A_{1}\circ A_{2} $ is defined everywhere.
This finishes the proof of Theorem~\ref{th4.90}.

\begin{remark} There is another argument why it is not possible to consider
closed operators $ A_{1} $ without making some corrections. Indeed, $ V_{2} $ being
quasi-complementary with $ V_{1} $ with the excess 0 heuristically implies that
$ \dim  V_{2} =\dim  H_{2} $. A change of $ V_{2} $ to a comparable subspace of $ H_{1}\oplus H_{2} $ of non-0
relative dimension would break this coincidence. On the other hand, if $ A_{2} $
is a closed unbounded operator, there there is another closed operator
$ A_{2}' $ such that $ \operatorname{Graph}\left(A_{2}\right) $ is comparable with $ \operatorname{Graph}\left(A_{2}'\right) $ with an arbitrary
relative dimension. Indeed, since considering $ A_{2}^{*} $ inverts the relative
dimension, it is enough to build $ A_{2} $ with $ \operatorname{Graph}\left(A_{2}'\right)\subset\operatorname{Graph}\left(A_{2}\right) $ with
codimension 1. Take a linear functional $ \alpha $ on $ H_{1}\oplus H_{2} $ with $ H_{2}\subset\operatorname{Ker}\alpha $. Let
$ V_{2}'=V_{2}\cap\operatorname{Ker}\alpha $. Then $ V_{2} $ is a closed subspace of $ H_{1}\oplus H_{2} $, and is a graph of a
partial operator. Obviously, the domain of this operator is dense unless
$ \alpha|_{H_{1}} $ is in the domain of $ A_{2}^{*} $. \end{remark}

\begin{remark} \label{rem6.60}\myLabel{rem6.60}\relax  It is easy to amplify this theorem by replacing compact
operators by some larger set of operators $ K $ such that $ 1-K $ is Fredholm.
The strongest such result (which is almost tautological) corresponds to
the class of operators of the form $ K_{1}+K_{0} $, with $ 1\notin\operatorname{Spec} K_{1} $, and $ K_{0} $ compact.
Another amplification, with the condition which is easier to check,
corresponds to operators of the form $ K_{1}+K_{0} $, with $ K_{0} $ compact, and the
spectral radius $ R $ of $ K_{1} $ satisfying $ R<1 $.

However, the stated form has the advantage of being invariant
w.r.t.~the change $ A_{2}'=A_{2}\circ B $. In fact, this is the strongest form which is
invariant w.r.t.~such changes. Indeed, it is enough to show that for any
non-compact operator $ A $ one can find a bounded operator $ B $ such that 1 is
in the essential spectrum of $ A\circ B $. In turn, by the polar decomposition
one can assume $ A $ to be self-adjoint; now the statement follows from the
spectral theorem and the characterization of compact self-adjoint
operators by their spectrum. \end{remark}

\bibliography{ref,outref,mathsci}

\def\cprime{$'$} \def\cprime{$'$} \def\cprime{$'$} \def\cprime{$'$}
  \def\cprime{$'$}
\providecommand{\bysame}{\leavevmode\hbox to3em{\hrulefill}\thinspace}
\providecommand{\MR}{\relax\ifhmode\unskip\space\fi MR }
\providecommand{\MRhref}[2]{%
  \href{http://www.ams.org/mathscinet-getitem?mr=#1}{#2}
}
\providecommand{\href}[2]{#2}
\begin{thebibliography}{10}

\bibitem{Bik00Rie}
I.~A. Bikchantaev, \emph{The {R}iemann problem on a finite-sheeted {R}iemann
  surface of infinite genus}, Mat. Zametki \textbf{67} (2000), no.~1, 25--35.
  \MR{2001d:30072}

\bibitem{FelKnoTru96Inf}
J.~Feldman, H.~Kn\"orrer, and Trubowitz E., \emph{Infinite genus riemann
  surfaces}, Canadian {M}athematical {S}ociety 1945-1995 (James~B. Carrell and
  Ram Murty, eds.), no.~3, Canadian Mathematical Society, Ottawa, 1996,
  pp.~91--112.

\bibitem{Fried82Int}
F.~G. Friedlander, \emph{Introduction to the theory of distributions},
  Cambridge University Press, Cambridge, 1982. \MR{86h:46002}

\bibitem{GelNai50Uni}
I.~M. Gel{\cprime}fand and M.~A. Na{\u\i}mark, \emph{Unitarnye predstavleniya
  klassi\v ceskih grupp}, Izdat. Nauk SSSR, Moscow-Leningrad, 1950.
  \MR{13,722f}

\bibitem{GesHol00Dar}
Fritz Gesztesy and Helge Holden, \emph{Darboux-type transformations and
  hyperelliptic curves}, J. Reine Angew. Math. \textbf{527} (2000), 151--183.
  \MR{2002b:37108}

\bibitem{GiesKnoTru93Geo}
D.~Gieseker, H.~Kn{\"o}rrer, and E.~Trubowitz, \emph{The geometry of algebraic
  {F}ermi curves}, Perspectives in Mathematics, vol.~14, Academic Press Inc.,
  Boston, MA, 1993.

\bibitem{CourHurv}
Adolf Hurwitz and R.~Courant, \emph{Vorlesungen \"uber allgemeine
  {F}unktionentheorie und elliptische {F}unktionen}, Interscience Publishers,
  Inc., New York, 1944. \MR{6,148e}

\bibitem{Ichi97Schot}
Takashi Ichikawa, \emph{Schottky uniformization theory on {R}iemann surfaces
  and {M}umford curves of infinite genus}, J. Reine Angew. Math. \textbf{486}
  (1997), 45--68.

\bibitem{Jin96Max}
Naondo Jin, \emph{On maximal {R}iemann surfaces}, Hiroshima Math. J.
  \textbf{26} (1996), no.~2, 385--404.

\bibitem{Kosh89Cur}
G.~A. Koshevo{\u\i}, \emph{Curves of infinite genus that can be uniformized in
  the sense of {S}chottky}, Uspekhi Mat. Nauk \textbf{44} (1989), no.~4(268),
  239--240.

\bibitem{McKTru76Hil}
H.~P. McKean and C.~Trubowitz, \emph{Hill's operator and hyperelliptic function
  theory in the presence of infinitely many branch points}, Comm. Pure Appl.
  Math. \textbf{29} (1976), no.~1, 143--226.

\bibitem{McKVan97Act}
H.~P. McKean and K.~L. Vaninsky, \emph{Action-angle variables for the cubic
  {S}chr\"odinger equation}, Comm. Pure Appl. Math. \textbf{50} (1997), no.~6,
  489--562. \MR{98b:35183}

\bibitem{Mer99Asy}
Franz Merkl, \emph{An asymptotic expansion for {B}loch functions on {R}iemann
  surfaces of infinite genus and almost periodicity of the
  {K}adomcev-{P}etviashvilli flow}, Math. Phys. Anal. Geom. \textbf{2} (1999),
  no.~3, 245--278. \MR{2000m:14037}

\bibitem{Mer00Rie}
\bysame, \emph{A {R}iemann {R}och theorem for infinite genus {R}iemann
  surfaces}, Invent. Math. \textbf{139} (2000), no.~2, 391--437.
  \MR{2001e:32030}

\bibitem{MulSchmiSchra98Hyp}
W.~M{\"u}ller, M.~Schmidt, and R.~Schrader, \emph{Hyperelliptic {R}iemann
  surfaces of infinite genus and solutions of the {K}d{V} equation}, Duke Math.
  J. \textbf{91} (1998), no.~2, 315--352. \MR{98m:58060}

\bibitem{PreSeg86Loo}
Andrew Pressley and Graeme Segal, \emph{Loop groups}, Oxford Mathematical
  Monographs, The Clarendon Press, Oxford University Press, New York, 1986,
  Oxford Science Publications.

\bibitem{Rau91Par}
Jeffrey Rauch, \emph{Partial differential equations}, Springer-Verlag, New
  York, 1991. \MR{94e:35002}

\bibitem{RudFun}
Walter Rudin, \emph{Functional analysis}, second ed., McGraw-Hill Inc., New
  York, 1991. \MR{92k:46001}

\bibitem{SatSat82Sol}
Mikio Sato and Yasuko Sato, \emph{Soliton equations as dynamical systems on
  infinite-dimensional {G}rassmann manifold}, Nonlinear partial differential
  equations in applied science (Tokyo, 1982), North-Holland, Amsterdam, 1983,
  pp.~259--271. \MR{86m:58072}

\bibitem{Schmi96Int}
Martin~U. Schmidt, \emph{Integrable systems and {R}iemann surfaces of infinite
  genus}, Mem. Amer. Math. Soc. \textbf{122} (1996), no.~581, viii+111.

\bibitem{SegWil85Loo}
Graeme Segal and George Wilson, \emph{Loop groups and equations of {K}d{V}
  type}, Inst. Hautes \'Etudes Sci. Publ. Math. (1985), no.~61, 5--65.
  \MR{87b:58039}

\bibitem{Vog87Uni}
David~A. Vogan, Jr., \emph{Unitary representations of reductive {L}ie groups},
  Princeton University Press, Princeton, NJ, 1987. \MR{89g:22024}

\bibitem{Zakh97Qua}
Ilya Zakharevich, \emph{Quasi-algebraic geometry of curves {I}, {T}he
  {R}iemann--{R}och theorem and {J}acobian},
  \url{http://xxx.lanl.gov/e-print/alg-geom/9710013}, August 1997.

\end{thebibliography}
\end{document}